\documentclass[11pt,reqno]{amsart}
\usepackage[T1]{fontenc}
\usepackage[utf8]{inputenc}
\usepackage{lmodern}
\usepackage{amsmath,amssymb,amsthm,mathtools}
\usepackage{microtype}
\usepackage[margin=1.2in]{geometry}
\usepackage{enumitem}
\usepackage{booktabs,array}
\usepackage{xcolor}
\usepackage{tikz}
\usetikzlibrary{arrows.meta,calc,decorations.markings,decorations.pathmorphing,patterns,positioning}
\usepackage[hidelinks]{hyperref}%

\theoremstyle{plain}
\newtheorem{theorem}{Theorem}[section]
\newtheorem{proposition}[theorem]{Proposition}
\newtheorem{lemma}[theorem]{Lemma}
\newtheorem{corollary}[theorem]{Corollary}
\newtheorem*{mainthm}{Main Theorem}
\newtheorem*{thmA}{Theorem A}
\newtheorem*{thmB}{Theorem B}
\newtheorem*{thmC}{Theorem C}
\theoremstyle{definition}
\newtheorem{definition}[theorem]{Definition}
\newtheorem{example}[theorem]{Example}
\newtheorem{question}[theorem]{Question}
\theoremstyle{remark}
\newtheorem{remark}[theorem]{Remark}
\numberwithin{equation}{section}

\newcommand{\PP}{\mathbb{P}}
\newcommand{\CC}{\mathbb{C}}
\newcommand{\ZZ}{\mathbb{Z}}
\newcommand{\QQ}{\mathbb{Q}}
\newcommand{\RR}{\mathbb{R}}
\newcommand{\DD}{\mathbb{D}}
\newcommand{\Qbar}{\overline{\mathbb{Q}}}
\DeclareMathOperator{\Sing}{Sing}
\DeclareMathOperator{\Reg}{Reg}
\DeclareMathOperator{\Br}{Br}
\DeclareMathOperator{\PC}{PC}
\DeclareMathOperator{\Supp}{Supp}
\DeclareMathOperator{\Fix}{Fix}
\DeclareMathOperator{\Aut}{Aut}
\DeclareMathOperator{\Lie}{Lie}
\newcommand{\Id}{\mathrm{Id}}
\newcommand{\cO}{\mathcal{O}}
\newcommand{\fm}{\mathfrak{m}}

\DeclareMathOperator{\Sym}{Sym}
\DeclareMathOperator{\ord}{ord}
\DeclareMathOperator{\gr}{gr}
\DeclareMathOperator{\pr}{pr}
\newcommand{\PCdiv}{\PC^{(1)}}           %
\newcommand{\Brdiv}{\Br^{(1)}}           %
\newcommand{\cG}{\mathcal{G}}            %
\newcommand{\cV}{\mathcal{V}}
\newcommand{\cW}{\mathcal{W}}
\newcommand{\cZ}{\mathcal{Z}}
\newcommand{\cR}{\mathcal{R}}            %
\newcommand{\tR}{\widetilde{\mathcal R}} %
\newcommand{\tX}{\widetilde{X}}
\newcommand{\hX}{\widehat{X}}
\newcommand{\hB}{\widehat{B}}
\newcommand{\hPhi}{\widehat{\Phi}}
\newcommand{\hvarphi}{\hat\varphi}
\newcommand{\Xo}{X^\circ}
\newcommand{\hXo}{\hX^\circ}
\newcommand{\hcO}{\widehat{\cO}}
\newcommand{\cOa}{\cO^{\mathrm{alg}}}
\newcommand{\fn}{\mathfrak{n}}
\newcommand{\fp}{\mathfrak{p}}
\newcommand{\fP}{\mathfrak{P}}
\newcommand{\BB}{\mathbb{B}}             %

\title[Super-saddle or repelling cycles]{Periodic cycles of post-critically finite endomorphisms are either super-saddle or repelling}
\author{Van Tu Le}
\address{Faculty of Mathematics and Informatics, Hanoi University of Science and Technology, Hanoi, Vietnam}
\email{tu.levan@hust.edu.vn}

\begin{document}
\begin{abstract}
We prove that every multiplier of a periodic point of a post-critically finite endomorphism of $\mathbb{P}^N$ of degree at least two is zero or has modulus greater than one. More generally, we extend the concept of post-critical finiteness to polarized endomorphisms of normal projective varieties (possibly singular), namely to the class of polarized endomorphisms admitting a branch boundary, i.e. a proper invariant algebraic set containing the singular locus and the branch loci of all iterates. We generalize the result on multipliers of periodic points to such maps, where the associated eigenvectors are considered on the Zariski tangent space. The core analytic ingredient is a boundary theorem for rotation discs. A maximal irrational rotation disc whose punctured image avoids the post-critical set has almost every radial limit in it. The centre of the rotation disc may be a singular point of the ambient space or in the post-critical set. This gives an analogue of the boundary property of Siegel discs in dimension one. When the map admits a branch boundary, it excludes rotation discs altogether. Independently, we prove a multiplier-transfer theorem for finite holomorphic semiconjugacies. It allows the inequalities of multipliers to be transferred through normalizations and finite quotients, including at singular points.
\end{abstract}

\keywords{Post-critically finite endomorphisms, periodic points, multipliers, rotation discs, polarized endomorphisms}
\subjclass[2020]{Primary 32H50, 37F80; Secondary 37F50, 37P55}

\maketitle

\section{Introduction}\label{sec:intro}

Let $f\colon\PP^N\to\PP^N$ be a holomorphic endomorphism of algebraic degree $d\geq2$, and let $C(f)$ be its critical set. The \emph{post-critical set} $\PC(f)$ is the closure of $\bigcup_{n\geq1}f^n(C(f))$, and $f$ is \emph{post-critically finite}, or PCF, if $\PC(f)$ is a proper algebraic subset of $\PP^N$. The \emph{multipliers} of $f$ at a periodic point $p$ of period $a$ are the eigenvalues of $D_pf^a$. A multiplier $\lambda$ is \emph{superattracting} if $\lambda=0$, \emph{attracting} if $0<|\lambda|<1$, \emph{neutral} if $|\lambda|=1$ and \emph{repelling} if $|\lambda|>1$. A periodic point is said to be {\it super-saddle} if each of its multipliers is superattracting or repelling, and {\it repelling} if all of its multipliers are repelling (Definition~\ref{def:multiplier}).

When $N=1$, every multiplier of a PCF map is superattracting or repelling \cite[Corollary~14.5]{Mil06}. We asked in \cite[Question~A]{Le22} whether the same holds in every dimension and we obtained a partial positive answer in \cite[Conjecture~A]{Le23}.  In \cite{Le22} we showed the positive answer in dimension two and, in every dimension, at periodic points outside $\PC(f)$; at a periodic point of $\Reg\PC(f)$, the multiplier transverse to $\PC(f)$ is zero or repelling \cite[Theorems~1.2 and~1.3, Proposition~4.1]{Le22}. Following Astorg, we obtained the positive answer for the class of maps which are weakly PCF all the way down \cite[Theorem~1.1]{Le23}, among them the maps constructed by Koch \cite{Koc13}. In this article, we give a positive answer without any condition.

\begin{mainthm}\label{thm:main}
Let $f\colon\PP^N\to\PP^N$ be a post-critically finite endomorphism of degree $d\geq2$. If $f^a(p)=p$ for some $a\geq1$, then every eigenvalue of $D_pf^a$ is zero or has modulus greater than one.
\end{mainthm}

A periodic point may have both zero and repelling multipliers. Consequently, every multiplier is algebraic, every conjugate of a nonzero multiplier has modulus greater than one (Corollary~\ref{cor:conjugates}), and connected holomorphic families of PCF endomorphisms are isospectral (Proposition~\ref{prop:families}).

Post-critically finite rational maps are central objects in one-dimensional complex dynamics. The higher-dimensional version was first considered in \cite{FS92}, where Forn\ae ss and Sibony coined the term {\it   critically finite endomorphisms} of $\PP^2$. As observed in \cite{Le22}, a theorem of Forn\ae ss and Sibony \cite[Theorem~6.1]{FS94} implies that the multipliers at periodic points outside $\PC(f)$ have modulus at least one when the complement of $\PC(f)$ is Kobayashi hyperbolic and hyperbolically embedded. Ueda \cite{Ued98} proved that the inverse branches of the iterates form normal families, and that a PCF endomorphism has no rotation domain, that is, no Fatou component on which a sequence of iterates converges to the identity \cite[Definition~4.1 and Proposition~5.1]{Ued98}. Jonsson \cite{Jon98} obtained repelling multipliers for a class of PCF endomorphisms of $\PP^2$ whose critical components are strictly preperiodic. Astorg \cite[Theorem~B]{Ast20} excluded attracting multipliers and roots of unity under a transversality assumption on the periodic components of $\PC(f)$, and then Ji \cite{Ji22} described the Julia sets of PCF endomorphisms of $\PP^2$. More recently, Dupuy, Levy, Tucker and Wang \cite{DLTW26} proved that PCF endomorphisms remain PCF, in a divisorial sense, on the normalizations of periodic divisors.

Recently, Sun and Zhao \cite[Theorem~A]{SZ26b} proved that every Fatou component of a PCF endomorphism of $\PP^N$ is preperiodic to an immediate basin of a superattracting cycle, which Rong had proved for $N=2$ \cite[Theorem~1.1]{Rong08}. They attach to a periodic point an invariant subvariety, its \emph{algebraic core}, whose normalization has a regular fixed point carrying the nonzero multipliers \cite[Theorem~3.2.1]{SZ26b}. They deduce that a periodic point has a repelling multiplier unless its differential is nilpotent \cite[Proposition~3.3.1]{SZ26b}. Their inverse branches \cite[Lemma~2.3.3]{SZ26b} give $|\lambda|\geq1$ for every nonzero multiplier $\lambda$, by the Cauchy estimates. The existence of neutral multipliers was not excluded entirely.

\subsection*{The setting}
By \cite{Le22}, we need to consider the periodic cycles in the post-critical set; the proofs have to pass to normalizations of invariant subvarieties, which may be singular. A \emph{polarized endomorphism} of a normal projective variety $X$ is a morphism $g\colon X\to X$ with $g^*H\simeq H^{\otimes M}$ for an ample line bundle $H$ and an integer $M\geq2$; it is finite and surjective. The \emph{branch locus} $\Br(g)$ is the set of points over which $g$ is not locally biholomorphic at some preimage. The \emph{post-critical set} $\PC(g)$ is the closure of $\Sing X\cup\bigcup_{n\geq1}\Br(g^n)$, and a \emph{branch boundary} for $g$ is a closed algebraic subset $B\subsetneq X$ such that
\[
g(B)\subset B,\qquad\Sing X\subset B\qquad\text{and}\qquad\Br(g^n)\subset B\quad\text{for every }n\geq1.
\]
Sun and Zhao use this hypothesis for maps induced on normalizations of invariant subvarieties of a PCF endomorphism of projective space \cite[Lemmas~2.3.2 and~4.2.1]{SZ26b}; we introduce the term \emph{branch boundary} to extend the concept of post-critical finiteness to polarized endomorphisms of projective varieties (possibly singular). On a normal variety, its existence is equivalent to divisorial post-critical finiteness (Proposition~\ref{prop:divisorial}), and the branch boundary need not have pure dimension. The precise comparison with their definitions is given in Section~\ref{ssec:branch-boundary}. For a PCF endomorphism $f$ of $\PP^N$, polarized by $\cO(1)$ with $M=d$, the set $\PC(f)$ is a branch boundary.

\subsection*{Rotation discs}
Let $D(0,r)\subset\CC$ be the disc of radius $r$ about $0$, $D^*(0,r):=D(0,r)\setminus\{0\}$, and $S^1$ the unit circle. A \emph{rotation disc} of $g$ for $\lambda\in S^1$ is a nonconstant holomorphic map $\varphi\colon D(0,\varepsilon)\to X$, with $\varepsilon>0$, such that $g(\varphi(z))=\varphi(\lambda z)$ for $z\in D(0,\varepsilon)$. On both routes to the Main Theorem, a neutral multiplier gives rise to such a disc. The boundary of a Siegel disc of a rational map lies in the closure of its post-critical set \cite[Theorem~11.17]{Mil06}; our first result is its analogue for rotation discs, whose centre $\varphi(0)$ may be a critical point of $g$ or a singular point of $X$.

\begin{thmA}\label{thm:A}
Let $g$ be a polarized endomorphism of a normal projective variety $X$, and let $E\subset X$ be a closed set such that
\[
g(E)\subset E\qquad\text{and}\qquad\Sing X\cup\bigcup_{n\geq1}\Br(g^n)\subset E.
\]
Let $\lambda  = e^{2 \pi i t}\in S^1, t \notin \QQ$ and $\varphi\colon D(0,\varepsilon)\to X$ be a rotation disc for $\lambda$ with $\varphi(D^*(0,\varepsilon))\cap E=\emptyset$. Let $R\in[\varepsilon,\infty]$ be the supremum of the radii $r$ such that $\varphi$ extends holomorphically to $D(0,r)$ with \mbox{$\varphi(D^*(0,r))\cap E=\emptyset$}.
\begin{enumerate}[label=\textup{(\arabic*)}]
\item We have $R<\infty$, and $\varphi$ extends to a rotation disc $\varphi\colon D(0,R)\to X$ for $\lambda$.
\item For almost every $\zeta\in S^1$, the radial limit $\varphi^*(\zeta):=\lim_{r\to R}\varphi(r\zeta)$ exists and lies in $E$. Whenever $\varphi^*(\zeta)$ exists, so does $\varphi^*(\lambda\zeta)$, and $\varphi^*(\lambda\zeta)=g(\varphi^*(\zeta))$.
\item If $A\subset X$ is a closed algebraic subset which does not contain $\varphi(D(0,\varepsilon))$, then the set of points $\zeta\in S^1$ such that $\varphi^*(\zeta)$ exists and lies in $A$ has measure zero.
\end{enumerate}
\end{thmA}

Theorem~A applies to the closed post-critical set without an algebraicity assumption: the smallest admissible $E$ is $\PC(g)$ (Lemma~\ref{lem:pc-admissible}). Its proof constructs a bounded equivariant lift of the maximal disc to the affine cone (Lemma~\ref{lem:lift-extension}). For an algebraic branch boundary $E=B$, the boundary concentration in~(2) and the boundary uniqueness in~(3) exclude the existence of such a disc. A disc contained in $B$ is treated on the normalization of its Zariski closure; hence we obtain the non-existence of rotation discs in the presence of an algebraic branch boundary.

\begin{corollary}\label{cor:no-rotation-disc}
Let $g$ be a polarized endomorphism of a normal projective variety $X$ which admits a branch boundary, and let $\lambda\in\CC$ with $|\lambda|=1$. There is no nonconstant holomorphic map $\varphi\colon D(0,\varepsilon)\to X$, with $\varepsilon>0$, such that $g(\varphi(z))=\varphi(\lambda z)$ for $z\in D(0,\varepsilon)$.
\end{corollary}

Corollary~\ref{cor:no-rotation-disc} is the analogue for discs of Ueda's theorem that a PCF endomorphism of $\PP^N$ has no rotation domain \cite[Proposition~5.1]{Ued98}; Examples~\ref{ex:sharp} and~\ref{ex:discs} discuss its hypotheses.

\subsection*{Multiplier transfer}
Our second result concerns holomorphic maps between complex spaces and multipliers of tangent maps on the Zariski tangent spaces. For a point $x$ of a reduced complex space $X$, with local ring $(\cO_{X,x},\fm_x)$, let $T_xX:=(\fm_x/\fm_x^2)^*$ be the Zariski tangent space. For a finite map $f \colon X \to Y$, denote by $D_x f \colon T_x X \to T_{f(x)}Y$ the differential of $f$ at $x$. A \emph{monomial of degree $k$} in a set of complex numbers is a product of $k$ of them, with repetitions allowed.

\begin{thmB}\label{thm:B}
Let $\pi\colon Y\to X$ be a finite surjective holomorphic map between irreducible reduced complex spaces, and let $f\colon X\to X$ and $g\colon Y\to Y$ be holomorphic maps with $\pi\circ g=f\circ\pi$. Let $p$ be a fixed point of $f$. There is an integer $N\geq1$, depending only on the germ of $\pi$ along $\pi^{-1}(p)$, with the following properties.
\begin{enumerate}[label=\textup{(\arabic*)}]
\item Assume that every periodic point of $g$ in $\pi^{-1}(p)$ is fixed. Then every nonzero eigenvalue of $D_pf$ is a monomial of degree between $1$ and $N$ in the eigenvalues of $D_yg$, for some fixed point $y\in\pi^{-1}(p)$ of $g$. If $g$ fixes every point of $\pi^{-1}(p)$, this holds for every eigenvalue of $D_pf$.
\item If $X$ is locally irreducible at $p$, then for every fixed point $y\in\pi^{-1}(p)$ of $g$, every eigenvalue of $D_pf$ is a monomial of degree between $1$ and $N$ in the eigenvalues of $D_yg$.
\end{enumerate}
\end{thmB}

Explicit values $N=N(\pi,p)$ will be given in Section~\ref{sec:transfer}. Theorem~B extends \cite[Proposition~2.2]{Le23}, where $X$ is the union of two smooth invariant hypersurfaces tangent at a fixed point. For symmetric products of a rational map, it expresses powers of their multipliers as monomials in those of the original map (Proposition~\ref{prop:symmetric}). We use it to carry multiplier inequalities through normalizations and finite quotients (Corollary~\ref{cor:descent}).

Our third result holds at singular points and contains the Main Theorem.

\begin{thmC}\label{thm:C}
Let $g$ be a polarized endomorphism of a normal projective variety $X$ which admits a branch boundary, and let $x$ be a periodic point of $g$ of period $a$. Then every eigenvalue of the differential $D_xg^a$ on the Zariski tangent space $T_xX$ is zero or has modulus greater than one.
\end{thmC}

The point $x$ may be singular and may lie in every branch boundary. By Theorem~B, the conclusion extends to a polarized endomorphism of a non-normal projective variety, provided that the map induced on the normalization admits a branch boundary (Corollary~\ref{cor:nonnormal}). However, the polarization assumption cannot be omitted (Example~\ref{ex:polarized}).

\subsection*{Two routes}
\begin{figure}[ht]
\centering
\begin{tikzpicture}[>={Stealth[length=2mm]},font=\small,
  box/.style={draw,rounded corners=2pt,align=center,inner sep=4pt,text width=5.75cm},
  shared/.style={draw,rounded corners=2pt,align=center,inner sep=4pt,text width=12.1cm,fill=gray!22},
  goal/.style={draw,line width=0.9pt,align=center,inner sep=4pt,text width=3.0cm}]
  \node[box,anchor=north,minimum height=33mm] (L1) at (-3.18,0) {\textbf{First route} (Section~\ref{sec:route-one})\\[2pt]
    $p$ periodic for $f$ on $\PP^N$\\ algebraic core of $p$ \cite{SZ26b}:\\ a regular fixed point $q$\\ on its normalization;\\ inverse branches of \cite{SZ26b}};
  \node[box,anchor=north,minimum height=33mm] (R1) at (3.18,0) {\textbf{Second route} (Section~\ref{sec:route-two})\\[2pt]
    $x$ periodic for $g$ on a normal $X$\\ completed local ring at $x$:\\ a fixed point $q$ on a normalized\\ invariant subvariety;\\ arcs with centred lifts \cite{SZ26b}\\ (Appendix~\ref{app:arcs})};
  \node[shared,anchor=north] (S1) at ($(L1.south)!0.5!(R1.south)+(0,-0.5)$) {bounded inverse branches at $q$ (Proposition~\ref{prop:local-retract}):\\ multipliers of modulus at least one, an idempotent limit $\rho$};
  \node[box,anchor=north] (L3) at ($(S1.south -| L1.south)+(0,-0.5)$) {the image of $\rho$ is a manifold:\\ linearization, a coordinate axis};
  \node[box,anchor=north] (R3) at ($(S1.south -| R1.south)+(0,-0.5)$) {the image of $\rho$ is a singular germ:\\ a torus action};
  \node[shared,anchor=north] (S2) at ($(L3.south)!0.5!(R3.south)+(0,-0.5)$) {a rotation disc, which contradicts Corollary~\ref{cor:no-rotation-disc}};
  \node[goal,anchor=north] (L4) at ($(S2.south -| L1.south)+(0,-0.5)$) {Main Theorem};
  \node[goal,anchor=north] (R4) at ($(S2.south -| R1.south)+(0,-0.5)$) {Theorem~C};
  \draw[->] (L1.south) -- (L1.south |- S1.north);
  \draw[->] (R1.south) -- (R1.south |- S1.north);
  \draw[->] (L3.north |- S1.south) -- (L3.north);
  \draw[->] (R3.north |- S1.south) -- (R3.north);
  \draw[->] (L3.south) -- (L3.south |- S2.north);
  \draw[->] (R3.south) -- (R3.south |- S2.north);
  \draw[->] (L4.north |- S2.south) -- (L4.north);
  \draw[->] (R4.north |- S2.south) -- (R4.north);
  \draw[->,dashed] (R4.west) -- node[above=1.5pt,inner sep=1pt] {$X=\PP^N$} (L4.east);
\end{tikzpicture}
\caption{The two routes to the Main Theorem. The grey steps are common to both.}
\label{fig:routes}
\end{figure}
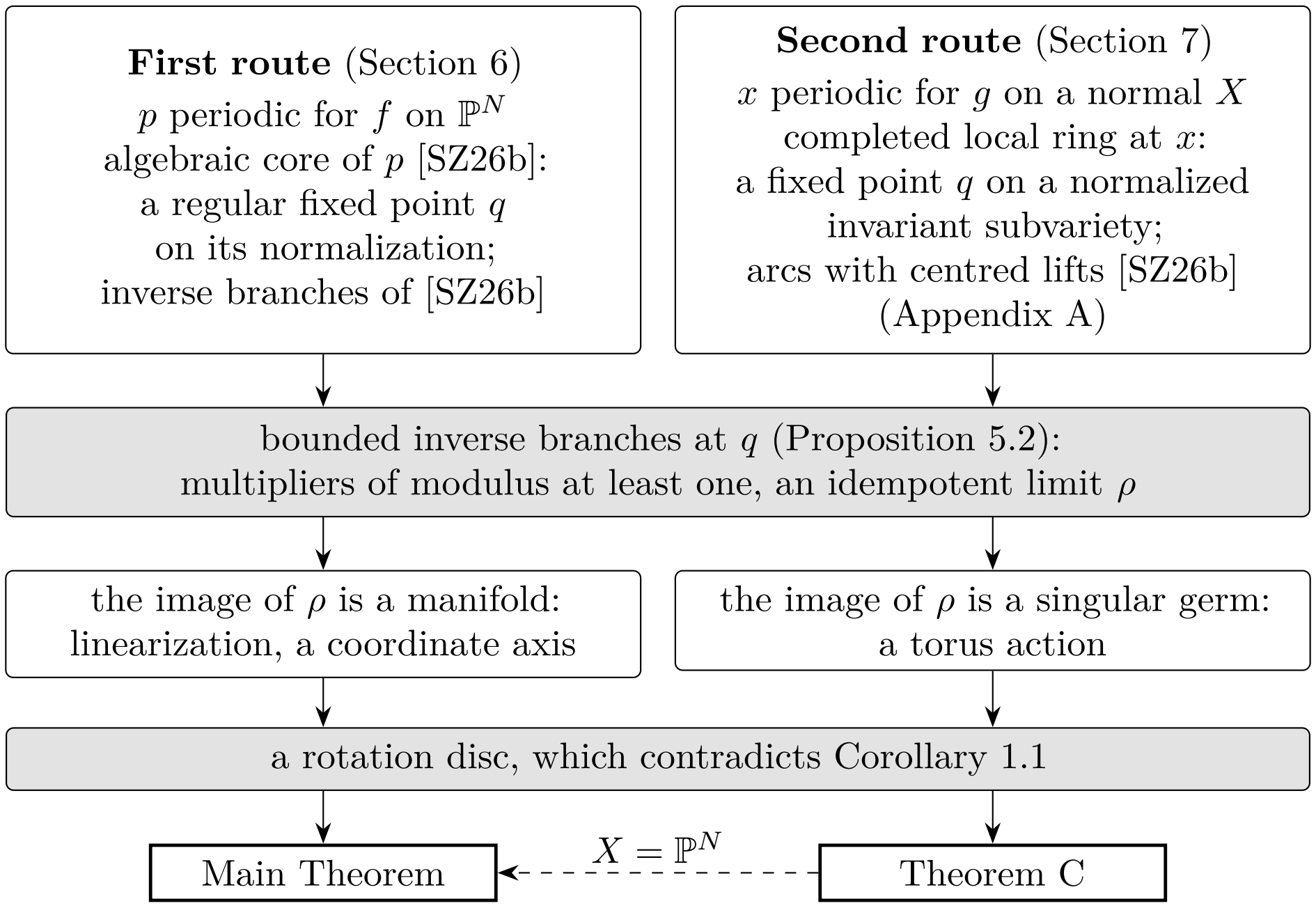

We prove the Main Theorem along two routes, which share two steps (Figure~\ref{fig:routes}). Both reach a fixed point $q$ of a map induced by an iterate on the normalization of an invariant subvariety, at which this map is a local biholomorphism with bounded inverse branches. At such a point the multipliers have modulus at least one, and a limit of inverse branches is an idempotent whose image carries the neutral directions (Proposition~\ref{prop:local-retract}). Both then find a rotation disc in this image, which contradicts Corollary~\ref{cor:no-rotation-disc}.

The first route (Section~\ref{sec:route-one}) is geometric and extends the method of \cite{Le22}. The point $q$ is the regular point given by the algebraic core, with the inverse branches of Sun and Zhao \cite[Theorem~3.2.1 and Lemma~2.3.3]{SZ26b}. The image of the idempotent is a manifold on which the map is linearizable, and a coordinate axis is a rotation disc.

The second route (Section~\ref{sec:route-two}) is algebraic. It does not use the algebraic core, but it uses the two tools with which Sun and Zhao construct the algebraic core at a periodic point of a normal variety: the part of the completed local ring on which the pull-back is invertible \cite[Section~3.1]{SZ26b} (Lemma~\ref{lem:fitting}), and their families of arcs which lift through every iterate \cite[Lemmas~3.2.4, 3.2.5 and~4.2.1]{SZ26b} (Appendix~\ref{app:arcs}). They yield an invariant subvariety and a point $q$ of its normalization, possibly singular, at which the induced map is a local biholomorphism. The image of the idempotent is then a singular germ, and the disc comes from a torus action. 

\subsection*{Relation with earlier work}
The proof of the rotation-disc exclusion develops the method of \cite[Section~3]{Le22}. We retain the lift to the affine cone, the Green-function bound, the lifted inverse branch, and the closed retract obtained from a limit of its iterates. The final contradiction again uses radial limits and the theorem of F.~and M.~Riesz.

The covering is now determined by the punctured disc, so its centre may lie in the branch boundary. Limits of inverse branches may meet this boundary in higher codimension; we lift each limit on the complement of its inverse image. On the resulting retract, a character of the closed group generated by the lifted inverse branch produces the rotation, and a radial flow constructs the maximal disc. This replaces the fixed-point linearization used in \cite{Le22} and develops the automorphism-group method of Ueda \cite[Appendix]{Ued98}.

If a rotation disc is contained in the boundary, we lift it to the normalization of its Zariski closure. The existence of a branch boundary there was proved for invariant curves of $\PP^2$ in \cite[Lemma~5.2]{Le22}, for invariant subvarieties of projective space by Sun and Zhao \cite[Lemmas~2.3.2 and~A.0.2]{SZ26b}, and for invariant prime divisors by Dupuy, Levy, Tucker and Wang \cite[Theorem~1.1]{DLTW26}. We give in Appendix~\ref{app:whitney} another proof for normal ambient varieties.

\subsection*{Organization}
Section~\ref{sec:prelim} contains the preliminaries. Section~\ref{sec:transfer} proves Theorem~B and can be read independently. Section~\ref{sec:rotation} proves Corollary~\ref{cor:no-rotation-disc}, then Theorem~A. After Section~\ref{sec:bib} on bounded inverse branches, Section~\ref{sec:route-one} proves the Main Theorem along the first route, and Section~\ref{sec:route-two} proves Theorem~C along the second. Consequences, examples and questions follow in Sections~\ref{sec:consequences}, \ref{sec:examples} and~\ref{sec:questions}.
\subsection*{AI usage} The author used multiple LLMs (ChatGPT, Claude) via a locally built harness to assist at several steps of the development of the manuscript, including revising the presentation and the reasoning steps, symbolic verification of examples, and drawing the illustration figures. The author is responsible for the mathematical contents and the final text.

\section{Preliminaries}\label{sec:prelim}

A variety is an irreducible reduced complex algebraic variety. We use the Euclidean topology unless the Zariski topology is mentioned. With the notation of Section~\ref{sec:intro} for discs, we set $\DD:=D(0,1)$, so that $S^1=\partial\DD$.

\subsection{Zariski tangent spaces and multipliers}\label{ssec:zariski}

\begin{definition}\label{def:multiplier}
Let $q$ be a point of a reduced complex space $X$, with analytic local ring $(\cO_{X,q},\fm_q)$. Set $V_q:=\fm_q/\fm_q^2$; its dual is the Zariski tangent space $T_qX$. A holomorphic germ $h\colon(X,q)\to(X,q)$ induces a linear map $h^*$ of $V_q$, whose transpose is the differential $D_qh$ of $T_qX$; the two maps have the same eigenvalues.

Let $x$ be a periodic point of a holomorphic self-map $h$ of $X$, of period $\ell$. The \emph{multipliers} of $h$ at $x$ are the eigenvalues of $D_xh^\ell$. The point $x$ is super-saddle if every multiplier $\lambda$ of $h$ at $x$ satisfies $\lambda=0$ or $|\lambda|>1$; it is called repelling if every multiplier has modulus strictly greater than $1$. A periodic cycle is super-saddle, respectively repelling, if each of its points is.
\end{definition}

If $X$ is algebraic, we denote by $\cOa_{X,q}\subset\cO_{X,q}$ its algebraic local ring. The two rings have the same completion $\hcO_{X,q}$ \cite[(36.7)]{Abh64}, hence the same cotangent space $V_q=\widehat\fm_q/\widehat\fm_q^2$.

A normal complex space $Y$ is locally irreducible \cite[Chapter~6, \S4.2]{GrRe84}. We also use three standard facts. (N1)~A continuous function on an open subset $W\subset Y$ which is holomorphic outside a nowhere dense analytic subset of $W$ is holomorphic on $W$ \cite[Chapter~7, \S4.2]{GrRe84}. (N2)~A finite holomorphic map between complex spaces of the same pure dimension, whose target is locally irreducible, is open \cite[Chapter~5, \S4.3]{GrRe84}. (N3)~If $\eta\colon(A,\fm_A)\to(C,\fm_C)$ is a local homomorphism of analytic algebras such that $C/\eta(\fm_A)C$ is finite dimensional, then $C$ is a finite $A$-module; this is the Weierstrass finiteness theorem \cite[Chapter~I, Theorem~1.10]{GLS07}.

We write $\BB(\rho)$ for the open ball of radius $\rho$ about $0$ in $\CC^\varkappa$. Let $W$ be an open subset of a normal complex space, identified with a closed analytic subset of $\BB(\rho)$. With the topology of locally uniform convergence, $\cO(\BB(\rho))$ and $\cO(W)$ are Fr\'echet spaces, by (N1) for $\cO(W)$ (see also \cite[Theorem~7.4.9]{Hor90}). The restriction map $r\colon\cO(\BB(\rho))\to\cO(W)$ is continuous, and surjective by Cartan's Theorem~B in the form \cite[Theorem~7.4.8]{Hor90}, as the ball is Stein \cite[Corollary~2.5.8 and p.~116]{Hor90}. By the open mapping theorem for Fr\'echet spaces \cite{BS17}, $r$ is open. Applied to the open set $\{H:\sup_{\BB(\rho/2)}|H|<1\}$, this gives a compact set $K\subset W$ and $\varepsilon>0$ with
\begin{equation}\label{eq:extension}
\bigl\{h\in\cO(W):\sup_K|h|<\varepsilon\bigr\}\subset r\bigl(\bigl\{H\in\cO(\BB(\rho)):\sup_{\BB(\rho/2)}|H|<1\bigr\}\bigr).
\end{equation}
We give here two lemmas which are versions of Montel's theorem for maps on complex spaces. We include the proofs for completeness.
\begin{lemma}\label{lem:montel}
Let $Y$ be a normal complex space whose topology has a countable basis. Every uniformly bounded sequence of holomorphic functions on $Y$ has a subsequence which converges locally uniformly on $Y$ to a holomorphic function.
\end{lemma}

\begin{proof}
Choose $C>0$ with $|f_j|\leq C$ on $Y$. Let $y\in Y$, and let $W\subset Y$ be an open neighbourhood of $y$ identified with a closed analytic subset of a ball $\BB(\rho)$, with $y=0$; let $K$ and $\varepsilon$ be as in \eqref{eq:extension}. Since $\sup_K|\varepsilon f_j/(2C)|<\varepsilon$, we can write $f_j|_W=(2C/\varepsilon)H_j|_W$ with $H_j\in\cO(\BB(\rho))$ and $\sup_{\BB(\rho/2)}|H_j|<1$. By Montel's theorem, a subsequence of $(H_j)$ converges locally uniformly on $\BB(\rho/2)$, so the corresponding subsequence of $(f_j)$ converges locally uniformly on the neighbourhood $W\cap\BB(\rho/2)$ of $y$, to a holomorphic function. As $Y$ has a countable basis, it is covered by countably many such neighbourhoods, and a diagonal extraction gives the required subsequence.
\end{proof}

\begin{lemma}\label{lem:jet}
Let $X$ be a normal complex space, $q\in X$, and let $X_\rho$ be an open neighbourhood of $q$ identified with a closed analytic subset of the ball $\BB(\rho)\subset\CC^\varkappa$, with $q=0$. Fix a norm on the cotangent space $V_q$. There are a compact set $K\subset X_\rho$ and $C>0$ such that the class $[h]\in V_q$ of every $h\in\cO(X_\rho)$ with $h(q)=0$ satisfies $\|[h]\|\leq C\sup_K|h|$. In particular, if $h_j\to h$ locally uniformly on $X_\rho$ and $h_j(q)=h(q)=0$, then $[h_j]\to[h]$.
\end{lemma}

\begin{proof}
Let $K$ and $\varepsilon$ be as in \eqref{eq:extension}, for $W=X_\rho$. Let $h\in\cO(X_\rho)$ with $h(q)=0$ and $\sup_K|h|<\varepsilon$, and let $H$ be an extension of $h$ with $\sup_{\BB(\rho/2)}|H|<1$. Then $H(0)=0$, and $H-\sum_k\partial_kH(0)z_k$ vanishes to order two at $0$. So $[h]=\sum_{k=1}^\varkappa\partial_kH(0)[z_k|_{X_\rho}]$, where $|\partial_kH(0)|\leq2/\rho$ by the Cauchy estimates, and $\|[h]\|\leq C_1:=(2/\rho)\sum_k\|[z_k]\|$. Applying this to $th$ for all $t>0$ if $\sup_K|h|=0$, and to $th$ with $t:=\varepsilon/(2\sup_K|h|)$ otherwise, gives $\|[h]\|\leq 2C_1\varepsilon^{-1}\sup_K|h|$ for every $h\in\cO(X_\rho)$ with $h(q)=0$.
\end{proof}

\subsection{Polarized endomorphisms admitting a branch boundary}\label{ssec:branch-boundary}

For a finite holomorphic endomorphism $g \colon X  \to X$  of a complex space, we denote by $C(g) \subset X$ the set of points at which $g$ is not locally biholomorphic. Its image $\Br(g):=g(C(g)) \subset X$ is the \emph{branch locus}. If $g$ is a finite morphism of varieties, it is locally biholomorphic at a point if and only if it is \'etale there \cite[Expos\'e~XII, Proposition~3.1 and its proof, Remarque~3.3(a)]{SGA1}, so that $C(g)$ and $\Br(g)$ are closed algebraic subsets. Denote by
\[
\PC(g) = \overline{\Sing X\cup\bigcup_{n\geq1}\Br(g^n)},
\]
where the closure is taken in the Euclidean topology.

\begin{definition}\label{def:branch-boundary}
Let $g\colon X\to X$ be a finite endomorphism of a projective variety.
\begin{enumerate}[label=\textup{(\arabic*)}]
\item The \emph{post-critical set} of $g$ is $\PC(g)$.
\item The map $g$ is \emph{polarized} if $g^*H\simeq H^{\otimes M}$ for an ample line bundle $H$ and an integer $M\geq2$.
\item A \emph{branch boundary} for $g$ is a closed algebraic subset $B\subsetneq X$ such that
\[
g(B)\subset B,\qquad\Sing X\subset B\qquad\text{and}\qquad\Br(g^n)\subset B\quad\text{for every }n\geq1.
\]
\item An endomorphism $f$ of $\PP^N$ of algebraic degree $d\geq2$, which is polarized by $\cO(1)$ with $M=d$, is \emph{post-critically finite}, or PCF, if $\PC(f)$ is a proper algebraic subset of $\PP^N$.
\end{enumerate}
\end{definition}

As $\Br(g^n)\subset\bigcup_{k<n}g^k(\Br(g))$, the third condition in~(3) is equivalent, given $g(B)\subset B$, to $\Br(g)\subset B$. A polarized endomorphism is proper and contracts no curve, as $g^*H$ is ample; hence it is finite, and surjective. Every iterate $g^k$ is polarized, with $(g^k)^*H\simeq H^{\otimes M^k}$, and every branch boundary for $g$ is one for $g^k$.

\begin{lemma}\label{lem:finite-fix}
Let $g$ be a polarized endomorphism of a projective variety. Then the set $\Fix(g^k)$ is finite for every $k\geq1$.
\end{lemma}

\begin{proof}
The set $\Fix(g^k)$ is closed and algebraic. An irreducible curve $C\subset\Fix(g^k)$ would satisfy $H|_C\simeq((g^k)^*H)|_C\simeq H^{\otimes M^k}|_C$, so that $(M^k-1)\deg(H|_C)=0$, which contradicts the ampleness of $H$.
\end{proof}

The name \emph{branch boundary}, and its use as a standing hypothesis on polarized endomorphisms of normal varieties, are essentially due to Sun and Zhao. They prove that the map induced by a PCF endomorphism of $\PP^N$ on the normalization of an invariant subvariety has such a set \cite[Lemma~2.3.2]{SZ26b}, which they describe as a fixed ramification boundary for every iterate \cite[Section~1.1]{SZ26b}; the set constructed in their proof is the minimal one of Lemma~\ref{lem:boundary}(3), and the inclusions (1) and~(2) of that lemma are stated there. In \cite[Lemma~4.2.1]{SZ26b} they take as a hypothesis a proper algebraic subset containing $\Sing X$ and the branch loci of all the iterates. For divisorially PCF maps the set $\Sing X\cup\PCdiv(g)$ plays this role \cite[Definition~3.2.1, Lemmas~3.2.3 and~3.2.4]{SZ26a}. We will discuss the relation between these notions below. 

When $X$ is smooth, $C(g)$ is the critical set, and $\Br(g^n)=\bigcup_{k=1}^ng^k(C(g))$ by the chain rule; so $\PC(g)$ is the closure of the usual post-critical set $\bigcup_{k\geq1}g^k(C(g))$, a countable union of hypersurfaces. Countably many hypersurfaces contained in a proper algebraic subset are finitely many. Hence, by Lemma~\ref{lem:boundary}(3) below, if $X$ is smooth, then $g$ admits a branch boundary if and only if the usual post-critical set is algebraic. In particular, if $f$ is a PCF endomorphism of $\PP^N$, then $\PC(f)=\bigcup_{k\geq1}f^k(C(f))$ is a hypersurface, and it is a branch boundary for $f$.

\begin{lemma}\label{lem:boundary}
Let $g$ be a finite endomorphism of a projective variety $X$.
\begin{enumerate}[label=\textup{(\arabic*)}]
\item We have $g(\Sing X)\subset\Sing X\cup\Br(g)$.
\item If $y\in\Br(g^n)$ and $X$ is locally irreducible at $y$, then $g(y)\in\Br(g^{n+1})$.
\item The map $g$ admits a branch boundary if and only if $\bigcup_{n\geq1}\Br(g^n)$, equivalently $\PC(g)$, is contained in a proper closed algebraic subset of $X$. Then the union of $\Sing X$ and of the Zariski closure of $\bigcup_{n\geq1}\Br(g^n)$ is the smallest branch boundary for $g$.
\end{enumerate}
\end{lemma}

\begin{proof}
(1) If $x\in\Sing X$ and $g(x)\notin\Br(g)$, then $g$ is locally biholomorphic at $x$, so $X$ is singular at $g(x)$.

(2) Let $x\in g^{-n}(y)$ be a point at which $g^n$ is not locally biholomorphic, and suppose that $g^{n+1}$ maps a neighbourhood $U$ of $x$ biholomorphically onto an open set. Then $g^n|_U$ is injective. The image of the germ of $g^n$ at $x$ is an analytic germ of dimension $\dim X$ \cite[Chapter~5, \S4.1]{GrRe84} contained in the irreducible germ $(X,y)$, hence equal to it \cite[Chapter~5, \S3.2]{GrRe84}. As $g^n$ is finite, after shrinking $U$, the map $g^n|_U$ is proper and injective onto a neighbourhood of $y$, hence a homeomorphism onto it. So $g$ maps the open set $g^n(U)$ bijectively onto the open set $g^{n+1}(U)$, with holomorphic inverse $g^n\circ(g^{n+1}|_U)^{-1}$. Thus $g|_{g^n(U)}$ is a biholomorphism, and $g^n|_U=(g|_{g^n(U)})^{-1}\circ g^{n+1}|_U$ is biholomorphic onto $g^n(U)$, a contradiction. Hence $g(y)=g^{n+1}(x)\in\Br(g^{n+1})$.

(3) Points at which $X$ is not locally irreducible are singular. So (1) and~(2) give $g(A)\subset A$ for $A:=\Sing X\cup\bigcup_{n\geq1}\Br(g^n)$. Since $g$ is continuous and closed for the Zariski topology, the Zariski closure $B$ of $A$ satisfies $g(B)\subset B$. It is a branch boundary when it is proper, and every branch boundary contains it. The two conditions of~(3) are equivalent, since $A'\cup\Sing X$ is a proper closed algebraic subset containing $\PC(g)$ whenever $A'$ is one containing $\bigcup_{n\geq1}\Br(g^n)$.
\end{proof}

\begin{lemma}\label{lem:sing-pc}
Let $g$ be a polarized endomorphism of a projective variety $X$. Then $\Sing X\subset\bigcup_{n\geq1}\Br(g^n)$.
\end{lemma}

\begin{proof}
Set $m:=\dim X$, and suppose that $s\in\Sing X$ lies in no $\Br(g^n)$. Fix $n\geq1$. As $g^n$ is finite and locally biholomorphic at every point of $g^{-n}(s)$, a neighbourhood $V$ of $s$ has a preimage which is a disjoint union of open sets $U_x$, $x\in g^{-n}(s)$, each mapped biholomorphically onto $V$. By the projection formula and the positivity of $H^m$ \cite[Tags~0BET and~0BEV]{Stacks}, $\deg g^n=M^{nm}$, and, over a dense Zariski open subset of the normal locus of $X$, the fibres of $g^n$ have $\deg g^n$ points \cite[Chapter~2, \S6.3, Theorem~2.29]{Sha13}. Hence $g^{-n}(s)$ has $M^{nm}$ points, all of them singular.

Let $S_1,\dots,S_r$ be the irreducible components of $\Sing X$, of dimensions $k_j\leq m-1$, and let $b$ be the number of irreducible components of the germ of $\Sing X$ at $s$. Put $W_j:=g^n(S_j)$ and $\delta_j:=\deg(g^n|_{S_j}\colon S_j\to W_j)$. By the projection formula and the positivity of $H^{k_j}\cdot W_j$ \cite[Tags~0BET and~0BEV]{Stacks}, $\delta_j\leq\delta_j\,(H^{k_j}\cdot W_j)=M^{nk_j}\,(H^{k_j}\cdot S_j)$.

Let $x\in S_j\cap g^{-n}(s)$. The germ of $S_j$ at $x$ is a union of irreducible components of the germ of $\Sing X$ at $x$. The biholomorphism $g^n|_{U_x}$ preserves the singular locus, so it maps this germ onto a union of components of the germ of $\Sing X$ at $s$, each of which is a component of the germ of $W_j$ at $s$; choose one of them, $\Gamma_{i(x)}$. The points of $W_j$ at which $W_j$ is normal and over which $g^n|_{S_j}$ has exactly $\delta_j$ preimages contain a dense Zariski open subset of $W_j$ \cite[Chapter~2, \S6.3, Theorem~2.29]{Sha13}, and this subset meets $\Gamma_i\cap V$, which has dimension $k_j=\dim W_j$. After shrinking $V$ so that $\Gamma_i\cap V\subset g^n(S_j\cap U_x)$ whenever $i(x)=i$, such a point $w\in\Gamma_i\cap V$ has a preimage in $S_j\cap U_x$ for each such $x$, so that at most $\delta_j$ points $x$ satisfy $i(x)=i$. Hence $S_j\cap g^{-n}(s)$ has at most $b\,\delta_j$ points, and
\[
M^{nm}\leq b\sum_{j=1}^rM^{nk_j}\,(H^{k_j}\cdot S_j)\leq C\,M^{n(m-1)},
\]
with $C$ independent of $n$. This is impossible for $n$ large, as $M\geq2$.
\end{proof}

On a normal variety, branch boundaries are determined by ramification in codimension one. Let $h\colon Y\to X$ be a finite surjective morphism of normal varieties. For a prime divisor $E\subset Y$, the local rings of $Y$ at $E$ and of $X$ at $h(E)$ are discrete valuation rings; let $e_h(E)$ be the ramification index of this extension. The ramification divisor $R_h:=\sum_E(e_h(E)-1)E$ is a finite sum, and we set $\Brdiv(h):=h(\Supp R_h)$. For an endomorphism $g$ of $X$, the \emph{divisorial post-critical set} is $\PCdiv(g):=\bigcup_{n\geq1}g^n(\Supp R_g)$, and $g$ is \emph{divisorially PCF} if $\PCdiv(g)$ is a finite union of prime divisors \cite[Definition~3.2.1]{SZ26a}. When $X$ is smooth, $\Supp R_g=C(g)$ and this is the usual notion \cite[Remark~3.2.2]{SZ26a}.

\begin{proposition}\label{prop:divisorial}
Let $g$ be a finite endomorphism of a normal projective variety $X$.
\begin{enumerate}[label=\textup{(\arabic*)}]
\item We have $\PCdiv(g)\subset\bigcup_{n\geq1}\Br(g^n)$ and $\bigcup_{n\geq1}\Br(g^n)\setminus\Sing X=\PCdiv(g)\setminus\Sing X$.
\item The map $g$ admits a branch boundary if and only if it is divisorially PCF. In that case $B_g:=\Sing X\cup\PCdiv(g)$ is the smallest branch boundary for $g$.
\item If moreover $g$ is polarized and admits a branch boundary, then $\PC(g)=\Br(g^n)=B_g$ for every large $n$.
\end{enumerate}
\end{proposition}

That $B_g$ is $g$-invariant and contains every $\Br(g^n)$ when $g$ is divisorially PCF is due to Sun and Zhao \cite[Lemmas~3.2.3 and~3.2.4]{SZ26a}; the equivalence in~(2) and the minimality of $B_g$ are proved here. Both proofs rest on the purity theorem of Zariski and Nagata.

\begin{proof}
(1) Let $x\in C(g)$ with $g(x)\in\Reg X$. By the purity theorem of Zariski and Nagata \cite[Expos\'e~X, Th\'eor\`eme~3.1]{SGA1}, applied to $g$ over $\Reg X$, the set of points at which $g$ is not \'etale has codimension one at $x$. So $x$ lies on a prime divisor $E$ at whose generic point $g$ is not \'etale, that is, $e_g(E)\geq2$. We apply this to $g^n$. As $e_{g^n}(E)=\prod_{i<n}e_g(g^i(E))$, we obtain $\Br(g^n)\setminus\Sing X\subset\bigcup_{k=1}^ng^k(\Supp R_g)$. Conversely, if $E$ is a prime divisor with $e_g(E)\geq2$, then $e_{g^k}(E)\geq e_g(E)\geq2$ by the same multiplicativity, so $g^k$ is \'etale at no point of $E$, and $g^k(E)\subset\Br(g^k)$ for every $k\geq1$.

(2) If $B$ is a branch boundary, each $g^k(E)$ with $E\subset\Supp R_g$ is a prime divisor contained in $B$ by~(1), hence one of the finitely many irreducible components of $B$. Conversely, if $g$ is divisorially PCF, then $B_g$ is a proper closed algebraic subset, it contains every $\Br(g^n)$ by~(1), and $g(B_g)\subset B_g$ by Lemma~\ref{lem:boundary}(1). Every branch boundary contains $\Sing X$ and $\PCdiv(g)$.

(3) By~(1) and Lemma~\ref{lem:sing-pc}, $\bigcup_{n\geq1}\Br(g^n)=B_g$, which is closed; hence $\PC(g)=B_g$. We claim that $\Br(g^n)\subset\Br(g^{n+1})$. Let $x\in C(g^n)$, and choose $x'\in g^{-1}(x)$, as $g$ is surjective. Suppose that $g^{n+1}$ maps a neighbourhood $U'$ of $x'$ biholomorphically onto an open set. Then $g|_{U'}$ is injective. As $X$ is locally irreducible at $x$ \cite[Chapter~6, \S4.2]{GrRe84}, the image of the germ of $g$ at $x'$ is the germ $(X,x)$ \cite[Chapter~5, \S\S3.2 and~4.1]{GrRe84}, so after shrinking $U'$ the finite map $g$ sends $U'$ homeomorphically onto a neighbourhood $U$ of $x$. The inverse $(g^{n+1}|_{U'})^{-1}\circ g^n|_U$ of $g|_{U'}$ is holomorphic, so $g^n|_U=g^{n+1}|_{U'}\circ(g|_{U'})^{-1}$ is biholomorphic onto an open set, a contradiction. Hence $x'\in C(g^{n+1})$ and $g^n(x)=g^{n+1}(x')\in\Br(g^{n+1})$. Each irreducible component of $B_g$ is the increasing union of its intersections with the closed algebraic subsets $\Br(g^n)$, so it is contained in one of them by Baire's theorem.
\end{proof}

On a normal variety, the minimal branch boundary $B_g=\Sing X\cup\PCdiv(g)$ consists of the divisorial part $\PCdiv(g)$ and of the irreducible components of $\Sing X$ not contained in it, which have codimension at least two \cite[Chapter~2, \S5.1, Theorem~2.19]{Sha13}. Such components occur (Examples~\ref{ex:kummer}, \ref{ex:P112} and~\ref{ex:mixed-dimension}).

The class of polarized endomorphisms of normal projective varieties admitting a branch boundary contains the maps induced by PCF endomorphisms of $\PP^N$ on the normalizations of invariant subvarieties (Proposition~\ref{prop:normalization-boundary}). It is closed under such restrictions and under products of maps with the same $M$: if $g_i\colon X_i\to X_i$ is polarized by $H_i$ and admits a branch boundary $B_i$, then $g_1\times g_2$ is polarized by $p_1^*H_1\otimes p_2^*H_2$, with $p_i$ the projections, and $B_1\times X_2\cup X_1\times B_2$ is a branch boundary for it, since a product of germs is smooth, and a product of maps is locally biholomorphic, exactly where both factors are. Lemma~\ref{lem:pcf-transport}(e) transports branch boundaries through finite quotient maps when the quotient branch locus has proper forward orbit closure; the equivariance hypothesis stated after that lemma ensures this condition.
\subsection{The cone and the Green function}\label{ssec:cone}

\begin{proposition}[{\cite[Corollary~2.2]{Fak03}}]\label{prop:fakhruddin}
Let $g$ be a polarized endomorphism of a normal projective variety $X$, with $g^*H\simeq H^{\otimes M}$. There are an embedding $\iota\colon X\to\PP^L$ and an endomorphism $\Phi$ of $\PP^L$ such that $\Phi\circ\iota=\iota\circ g$ and $\iota^*\cO(1)\simeq H^{\otimes k}$ for some $k\geq1$.
\end{proposition}

If $\dim X\geq1$, then $\Phi$ has algebraic degree $M$: if $D$ is its degree, then $H^{\otimes kD}\simeq\iota^*\Phi^*\cO(1)\simeq g^*\iota^*\cO(1)\simeq H^{\otimes kM}$, and $H$ is ample.

Let $g$, $H$ and $M$ be as above, with $\dim X\geq1$. We fix $\iota$ and $\Phi$ as in Proposition~\ref{prop:fakhruddin}, identify $X$ with $\iota(X)\subset\PP^L$, and choose homogeneous polynomials $\hPhi=(\hPhi_0,\dots,\hPhi_L)$ of degree $M$ representing $\Phi$; it is homogeneous and \emph{non-degenerate}:
\[
\hPhi(tv)=t^M\hPhi(v)\qquad\text{and}\qquad\hPhi^{-1}(0)=\{0\}.
\]
The affine cone $\hX:=\{0\}\cup\{v\in\CC^{L+1}\setminus\{0\}:[v]\in X\}$ is closed in $\CC^{L+1}$. The map $\varpi\colon\hX\setminus\{0\}\to X$, $v\mapsto[v]$, is the restriction of the tautological $\CC^*$-bundle over $\PP^L$, hence a holomorphic principal $\CC^*$-bundle over all of $X$, including $\Sing X$. Since $\Phi(X)\subset X$, the map $\hat g:=\hPhi|_{\hX}$ sends $\hX$ to itself, and
\[
\varpi\circ\hat g=g\circ\varpi\quad\text{on }\hX\setminus\{0\},\qquad\hat g^{-1}(0)=\{0\}\qquad\text{and}\qquad\hat g(tv)=t^M\hat g(v).
\]

Assume that $g$ admits a branch boundary $B$. There are homogeneous polynomials $Q_1,\dots,Q_l$ of positive degree with $B=\{[v]:Q_1(v)=\dots=Q_l(v)=0\}$ (the coordinate functions if $B=\emptyset$). We set $\Xo:=X\setminus B$,
\[
\hB:=\varpi^{-1}(B)\cup\{0\}=\hX\cap\{Q_1=\dots=Q_l=0\}\qquad\text{and}\qquad\hXo:=\hX\setminus\hB=\varpi^{-1}(\Xo).
\]

\begin{lemma}\label{lem:covering}
\begin{enumerate}[label=\textup{(\arabic*)}]
\item The sets $\Xo$ and $\hXo$ are connected complex manifolds. Moreover $\hat g(\hB)\subset\hB$ and $\hat g^{-n}(\hXo)\subset\hXo$ for every $n\geq1$.
\item For every $n\geq1$, the maps $g^n\colon g^{-n}(\Xo)\to\Xo$ and $\hat g^n\colon\hat g^{-n}(\hXo)\to\hXo$ are finite coverings.
\end{enumerate}
\end{lemma}

\begin{proof}
(1) As $\Sing X\subset B$, we have $\Xo\subset\Reg X$. The manifold $\Reg X$ is connected, as the set of complex points of an irreducible variety \cite[Chapter~7, \S2, Theorem~7.1]{Sha13II}, and $B\cap\Reg X$ is a proper analytic subset of it.  Hence, $\Xo$ is a connected manifold, and so is the $\CC^*$-bundle $\hXo$ over it. The inclusions follow from $g(B)\subset B$, $\hat g^{-1}(0)=\{0\}$ and $\hat g^{-n}(\hXo)=\varpi^{-1}(g^{-n}(\Xo))$.

(2) The set $Y_n:=g^{-n}(\Xo)$ is open and contained in $\Xo$. The restriction $g^n\colon Y_n\to\Xo$ of the finite map $g^n$ over an open set is proper with finite fibres, and it is a local homeomorphism because $g^n(Y_n)\subset\Xo$ and $\Br(g^n)\subset B$. Its image is open, closed and nonempty, as $g$ is surjective, hence equal to the connected set $\Xo$. A proper surjective local homeomorphism with finite fibres between locally compact Hausdorff spaces is a finite covering.

On the cone, let $s$ and $s'$ be holomorphic sections of $\varpi$ over open sets $V,V'\subset X$ with $g^n(V)\subset V'$. In the coordinates $(t,x)\mapsto ts(x)$ and $(t,x')\mapsto ts'(x')$, homogeneity gives $\hat g^n(t,x)=(t^{M^n}u_n(x),g^n(x))$ with $u_n$ holomorphic and nowhere zero. So $\hat g^n$ is locally biholomorphic on $\varpi^{-1}(Y_n)=\hat g^{-n}(\hXo)$, and it is proper over $\hXo$ since $|\hPhi^n(v)|\geq m_n|v|^{M^n}$ for some $m_n>0$. As for $g^n$, it is a finite covering.
\end{proof}

\begin{lemma}[{Hubbard--Papadopol \cite[Theorem~2.1]{HP94}}]\label{lem:green}
The functions $M^{-k}\log|\hPhi^k|$ converge uniformly on $\CC^{L+1}\setminus\{0\}$ to a continuous function $\cG$, and there is a constant $C_0$ such that
\[
\cG\circ\hPhi=M\,\cG\qquad\text{and}\qquad\bigl|\cG(v)-\log|v|\bigr|\leq C_0\quad(v\neq0).
\]
In particular $\{\cG=0\}$ is a compact subset of $\CC^{L+1}\setminus\{0\}$.
\end{lemma}

\begin{proof}
The map $\hPhi$ is homogeneous of degree $M\geq2$ and non-degenerate. By \cite[Theorem~2.1]{HP94}, there is a neighbourhood $V$ of $0$ on which the functions $M^{-j}\log|\hPhi^j|$ converge, uniformly on compact subsets of $V\setminus\{0\}$, to a limit $\cG$ that satisfies the two relations on $V\setminus\{0\}$. Since $M^{-j}\log|\hPhi^j(tv)|=\log|t|+M^{-j}\log|\hPhi^j(v)|$, the function $M^{-j}\log|\hPhi^j(v)|-\log|v|$ depends only on $[v]$. Hence uniform convergence on a sphere contained in $V$ gives uniform convergence on $\CC^{L+1}\setminus\{0\}$, and both relations extend. The set $\{\cG=0\}$ is closed in $\CC^{L+1}\setminus\{0\}$ and contained in $\{e^{-C_0}\leq|v|\leq e^{C_0}\}$.
\end{proof}

\subsection{Tools from complex analysis}\label{ssec:tools}

A \emph{holomorphic retraction} of a complex manifold $Y$ is a holomorphic map $\rho\colon Y\to Y$ with $\rho\circ\rho=\rho$. Its image is its fixed-point set, and it is connected if $Y$ is connected.

\begin{lemma}[Rossi {\cite[Lemma~2.1.28, p.~119]{Aba89}}]\label{lem:rossi}
The image of a holomorphic retraction of a complex manifold $Y$ is a closed complex submanifold of $Y$.
\end{lemma}

Let $\varrho$ be the Poincar\'e distance on $\DD$. The Kobayashi pseudodistance $d_{\mathcal X}(x,y)$ of a complex manifold $\mathcal X$ is the infimum of $\sum_{i=1}^k\varrho(a_i,b_i)$ over all chains of holomorphic discs $f_i\colon\DD\to\mathcal X$ and points $a_i,b_i\in\DD$ with $f_1(a_1)=x$, $f_i(b_i)=f_{i+1}(a_{i+1})$ and $f_k(b_k)=y$. Holomorphic maps decrease $d$, and $\mathcal X$ is \emph{hyperbolic} if $d_{\mathcal X}$ is a distance \cite{Kob67}.

\begin{lemma}\label{lem:hyperbolic}
Let $\mathcal X$ be a connected complex manifold of positive dimension and $\chi\colon\mathcal X\to\CC^m$ a locally injective holomorphic map whose image is contained in a ball. Then $\mathcal X$ is hyperbolic, and $d_{\mathcal X}$ induces the topology of $\mathcal X$.
\end{lemma}

\begin{theorem}[Kobayashi {\cite[Theorem~5.4.2(1), p.~262]{Kob98}}, {\cite[Theorem~6.2]{Kob67}}]\label{thm:kobayashi}
Let $\mathcal X$ be a connected hyperbolic manifold. Then the group $\Aut(\mathcal X)$ of holomorphic automorphisms of $\mathcal X$, with the compact-open topology, is a real Lie group.
\end{theorem}

We use the following facts on Lie groups. A closed subgroup $\Gamma$ of a Lie group is an embedded Lie subgroup, with the subspace topology \cite[Theorem~20.12]{Lee13}. Its identity component $\Gamma^\circ$ is open in $\Gamma$ and generated by $\exp(\Lie\Gamma)$ \cite[Propositions~7.14, 7.15 and~20.8(f)]{Lee13}. A continuous homomorphism of connected Lie groups is smooth \cite[Chapter~IV, \S XIII, Proposition~1]{Che46}. Lie groups are metrizable, and the closure of an abelian subgroup of a topological group is abelian. If $\mathcal X$ is locally compact and Hausdorff, evaluation $\Aut(\mathcal X)\times\mathcal X\to\mathcal X$ is continuous for the compact-open topology; so a continuous one-parameter subgroup $t\mapsto h_t$ of $\Aut(\mathcal X)$ gives a continuous map $(t,x)\mapsto h_t(x)$ on $\RR\times\mathcal X$.

By Fatou's theorem \cite[Theorem~11.32, p.~249]{Rud87}, a bounded holomorphic function on $\DD$ has radial limits at almost every point of $S^1$. We apply it componentwise to bounded holomorphic maps $\DD\to\CC^{L+1}$, and to discs $D(0,R)$ after rescaling.

\begin{lemma}[F.~and M.~Riesz {\cite[Theorem~A.3, p.~220]{Mil06}}, {\cite[Theorem~15.19, pp.~308--309]{Rud87}}]\label{lem:riesz}
Let $u\colon\DD\to\CC$ be bounded and holomorphic. If the radial limit $\lim_{r\to1}u(r\zeta)$ exists and vanishes for $\zeta$ in a subset of $S^1$ of positive Lebesgue measure, then $u\equiv0$.
\end{lemma}

\begin{theorem}[Ueda {\cite[Theorem~2.1, p.~320]{Ued98}}]\label{thm:ueda-normality}
Let $\Phi$ be an endomorphism of $\PP^m$ of algebraic degree at least $2$, let $Z$ be a complex space, and let $\psi_0\colon Z\to\PP^m$ be a holomorphic map. Let $(\psi_\nu)$ be a family of holomorphic maps $Z\to\PP^m$ such that, for every $\nu$, $\Phi^{j_\nu}\circ\psi_\nu=\psi_0$ for some integer $j_\nu\geq0$. Then the family $(\psi_\nu)$ is normal.
\end{theorem}

We use Ueda's notion of analytic covering \cite[Section~3]{Ued98}. Let $\Omega$ be a connected complex manifold and $Y$ an irreducible normal complex space. A holomorphic map $\eta\colon Y\to\Omega$ is an \emph{analytic covering} if it is proper and surjective with finite fibres; its \emph{order} $\operatorname{ord}(\eta,y)$ at $y\in Y$ is its local sheet number at $y$. For an analytic subset $\Sigma\subset\Omega$ of codimension one, $\eta$ is \emph{$\Sigma$-branched} if $\eta(\{y:\operatorname{ord}(\eta,y)>1\})\subset\Sigma$.

\begin{lemma}[Ueda {\cite[Lemma~3.8, p.~324]{Ued98}}]\label{lem:ueda-local}
Let $\Omega'$ be a complex manifold and $\Sigma$ an analytic subset of $\Omega'$ of codimension one. For every $x_0\in\Omega'$ there is a neighbourhood $\Omega$ of $x_0$ such that, for every $\Sigma\cap\Omega$-branched analytic covering $\eta\colon Y\to\Omega$, the fibre $\eta^{-1}(x_0)$ consists of a single point.
\end{lemma}

This is part~(ii) of \cite[Lemma~3.8]{Ued98}, and $\Omega$ is the neighbourhood given there. The definition of an analytic covering of $\Omega$ requires $\Omega$ to be connected, and it is: by part~(i) of the same lemma, $\Omega$ is the image of an analytic covering, whose source is irreducible, hence connected.

\subsection{Invariant subvarieties and normalization}\label{ssec:invariant}

Let $g$ be a polarized endomorphism of a normal projective variety $X$, with $g^*H\simeq H^{\otimes M}$, let $Z\subset X$ be an irreducible subvariety of dimension at least one with $g(Z)=Z$, and let $\nu\colon X_Z\to Z$ be its normalization. By the universal property of the normalization \cite[Chapter~2, \S5.2, Theorem~2.21(ii)]{Sha13}, the finite surjective morphism $g\circ\nu\colon X_Z\to Z$ factors as $g\circ\nu=\nu\circ g_Z$ for a morphism $g_Z\colon X_Z\to X_Z$. It is unique, since $\nu$ is injective over the normal locus of $Z$, whose preimage under $g\circ\nu$ is dense in $X_Z$, and it is finite and surjective. We call $g_Z$ the map \emph{induced} by $g$. Every iterate $g^k$ satisfies $g^k(Z)=Z$ and induces $g_Z^k$.

\begin{proposition}\label{prop:normalization-boundary}
Assume that $g$ admits a branch boundary $B$, and let $Z$, $\nu\colon X_Z\to Z$ and $g_Z$ be as above. Then $g_Z$ is a polarized endomorphism of $X_Z$, with $g_Z^*\nu^*H\simeq(\nu^*H)^{\otimes M}$, and it admits a branch boundary.
\end{proposition}

No condition is imposed on the position of $Z$, which may lie in $B$ or in $\Sing X$. For $X=\PP^N$ and $g$ a PCF endomorphism, the proposition is due to Sun and Zhao \cite[Lemmas~2.3.2 and~A.0.2]{SZ26b}. When $Z$ is a prime divisor, the existence of a branch boundary on $X_Z$ follows, by Proposition~\ref{prop:divisorial}(2), from the theorem of Dupuy, Levy, Tucker and Wang \cite[Theorem~1.1]{DLTW26}, proved by purity and Abhyankar's lemma. We give another proof in general, which uses Whitney stratifications and is given in Appendix~\ref{app:whitney}.

\begin{lemma}\label{lem:Z}
Let $g$ be a finite endomorphism of a projective variety $X$, let $\mu\in S^1$, and let $\psi\colon D(0,\varepsilon)\to X$ be a nonconstant holomorphic map with $g\circ\psi(z)=\psi(\mu z)$. Let $Z$ be the Zariski closure of $\psi(D(0,\varepsilon))$. Then $Z$ is irreducible, $\dim Z\geq1$ and $g(Z)=Z$. 

Moreover, for every closed algebraic subset $A\subset X$ with $Z\not\subset A$, $\psi^{-1}(A)$ is closed and discrete in $D(0,\varepsilon)$.
\end{lemma}

\begin{proof}
Write $D:=D(0,\varepsilon)$. For a closed algebraic subset $A\subset X$, the analytic subset $\psi^{-1}(A)$ of $D$ is either $D$, and then $Z\subset A$, or closed and discrete. This gives the last assertion, and the irreducibility of $Z$: if $Z=A\cup A'$, then $D=\psi^{-1}(A)\cup\psi^{-1}(A')$ is uncountable, so $Z\subset A$ or $Z\subset A'$. As $\psi(D)$ is infinite, $\dim Z\geq1$. The identity $\psi(z)=g(\psi(\mu^{-1}z))$ gives $\psi(D)\subset g(Z)$. Since $g$ is finite, $g(Z)$ is closed, irreducible and of dimension $\dim Z$; it contains $Z$, hence $g(Z)=Z$.
\end{proof}
Such a map $\psi$ is called a {\it rotation disc} of $g$. The next lemma shows that we can lift a rotation disc via the normalization. 
\begin{lemma}\label{lem:N}
In the setting of Lemma~\ref{lem:Z}, let $\nu\colon X_Z\to Z$ be the normalization, and $g_Z\colon X_Z\to X_Z$ be the lift of $g$. There is a unique holomorphic lift $\psi_Z\colon D(0,\varepsilon)\to X_Z$ of $\psi$ via $\nu$, i.e. $\nu\circ\psi_Z=\psi$. Moreover, $\psi_Z$ is a rotation disc of $g_Z$ and $\psi_Z(D(0,\varepsilon))$ is Zariski dense in $X_Z$.

In addition, if $\mu$ is not a root of unity and $B_Z\subsetneq X_Z$ is a closed algebraic subset with $g_Z(B_Z)\subset B_Z$, then $\psi_Z(D^*(0,\varepsilon))\subset X_Z\setminus B_Z$.
\end{lemma}

\begin{proof}
Write $D:=D(0,\varepsilon)$. Let $\mathcal N\subset Z$ be the non-normal locus, a proper closed algebraic subset; $\nu$ is finite, and biholomorphic over $Z\setminus\mathcal N$. By Lemma~\ref{lem:Z}, $J:=\psi^{-1}(\mathcal N)$ is closed and discrete in $D$. Set $\psi_Z:=(\nu|_{X_Z\setminus\nu^{-1}(\mathcal N)})^{-1}\circ\psi$ on $D\setminus J$.

Let $a\in J$, and choose $s_0>0$ with $\overline D(a,s_0)\cap J=\{a\}$. For $0<s\leq s_0$, the sets $C_s:=\overline{\psi_Z(\overline D(a,s)\setminus\{a\})}$ are connected and lie in the compact set $\nu^{-1}(\psi(\overline D(a,s_0)))$. Their intersection is thus nonempty, connected and contained in the finite fibre $\nu^{-1}(\psi(a))$, hence a single point $x_a$, and $\psi_Z(z)\to x_a$ as $z\to a$. By Riemann's extension theorem in a local embedding of $X_Z$ near $x_a$, $\psi_Z$ is holomorphic at $a$. So $\psi_Z$ extends holomorphically to $D$; it is unique, as two lifts agree off $J$, and nonconstant, as $\psi$ is.

For $z\notin J\cup\mu^{-1}J$, the points $g_Z(\psi_Z(z))$ and $\psi_Z(\mu z)$ both lie over $g(\psi(z))=\psi(\mu z)\notin\mathcal N$, so they coincide; by continuity, $g_Z\circ\psi_Z(z)=\psi_Z(\mu z)$ on $D$. If $\psi_Z(D)\subset A$ with $A$ closed algebraic, then $\psi(D)$ lies in the closed set $\nu(A)$, so $\nu(A)=Z$, and $A=X_Z$ as $X_Z$ is irreducible of dimension $\dim Z$. This proves the density.

Finally, $\psi_Z^{-1}(B_Z)$ is closed and discrete in $D$, by density, and $\psi_Z(\mu z)=g_Z(\psi_Z(z))\in B_Z$ for $z\in\psi_Z^{-1}(B_Z)$. As $\mu$ is not a root of unity, a point $z\neq0$ of $\psi_Z^{-1}(B_Z)$ would give infinitely many points $\mu^kz$ of it on the compact circle $\{|w|=|z|\}$. So $\psi_Z^{-1}(B_Z)\subset\{0\}$.
\end{proof}

\section{Multiplier transfer}\label{sec:transfer}

We prove Theorem~B with an explicit integer $N$, and draw its consequences for multipliers and for post-critical finiteness.

Let $\pi\colon Y\to X$, $f$, $g$ and $p$ be as in Theorem~B, and set $S:=\pi^{-1}(p)$. This fibre is finite, and $g(S)\subset S$. For $y\in S$ and $a\in\cO_{Y,y}\setminus\{0\}$, let $\ord_y(a)$ be the largest integer $k$ with $a\in\fm_y^k$, and set $\ord_y(0):=+\infty$. For $h\in\cO_{X,p}$, let $\ord_y(\pi^*h)$ refer to the germ of $h\circ\pi$ at $y$, and set $\ord_S(\pi^*h):=\min_{y\in S}\ord_y(\pi^*h)$. Define
\[
N_S(\pi,p):=\max_{h\in\fm_p\setminus\fm_p^2}\ord_S(\pi^*h),\qquad N_y(\pi,p):=\max_{h\in\fm_p\setminus\fm_p^2}\ord_y(\pi^*h)\quad(y\in S),
\]
with the convention $\max\emptyset:=1$. Let $N(\pi,p):=\max_{y\in S}N_y(\pi,p)$ if $X$ is locally irreducible at $p$, and $N(\pi,p):=N_S(\pi,p)$ otherwise. These quantities depend only on the germ of $\pi$ along $S$, and $N_S(\pi,p)\leq N_y(\pi,p)$. An individual $N_y(\pi,p)$ may be infinite when $X$ is locally reducible at $p$: a nonzero cotangent class may pull back to zero on one branch of a node. The integer $N(\pi,p)$ is finite by Lemma~\ref{lem:transfer-filtration} and the proof below, and is the bound in Theorem~B.

\subsection{Local rings and the proof of Theorem~B}
A \emph{local ring} is a Noetherian local $\CC$-algebra $(R,\fm)$ with residue field $\CC$, for instance the local ring of a complex space at a point; its cotangent space $\fm/\fm^2$ is finite dimensional. An endomorphism $\sigma$ of $R$ with $\sigma(\fm)\subset\fm$ induces linear maps $\gr^k\sigma$ of $\fm^k/\fm^{k+1}$, and we write $\overline\sigma:=\gr^1\sigma$. The same notation is used for finite products of local rings. The multiplication map $\Sym^k(\fm/\fm^2)\to\fm^k/\fm^{k+1}$ is onto and intertwines $\Sym^k\overline\sigma$ with $\gr^k\sigma$. Triangularizing $\overline\sigma$, we obtain:
\begin{equation}\label{eq:graded-eigenvalues}
\text{every eigenvalue of $\gr^k\sigma$ is a product of $k$ eigenvalues of $\overline\sigma$.}
\end{equation}

\begin{lemma}\label{lem:transfer-filtration}
Let $(R,\fm)$ and $(R_1,\fm_1),\dots,(R_s,\fm_s)$ be local rings, $R':=R_1\times\dots\times R_s$ and $J:=\fm_1\times\dots\times\fm_s$. Let $\iota\colon R\to R'$ be an injective homomorphism such that $R'$ is a finite $R$-module and $\iota(\fm)\subset J$. Then $\iota^{-1}(J^{N+1})\subset\fm^2$ for some integer $N\geq1$; let $N(\iota)$ be the least one. Let $\sigma_R$ and $\sigma_{R'}$ be endomorphisms of $R$ and $R'$ with $\sigma_R(\fm)\subset\fm$, $\sigma_{R'}(J)\subset J$ and $\sigma_{R'}\circ\iota=\iota\circ\sigma_R$. Then every eigenvalue of $\overline{\sigma}_R$ is an eigenvalue of $\gr^k\sigma_{R'}$ on $J^k/J^{k+1}$ for some $1\leq k\leq N(\iota)$.
\end{lemma}

\begin{proof}
We identify $R$ with $\iota(R)$. The maximal ideals of $R'$ are the ideals $\fm_i\times\prod_{j\neq i}R_j$; they contain $\fm R'$, and their intersection is $J$. The decreasing chain of subspaces $(J^k+\fm R')/\fm R'$ of the finite-dimensional space $R'/\fm R'$ is stationary, so that the finite $R'$-module $K:=(J^c+\fm R')/\fm R'$ satisfies $JK=K$ for some $c\geq1$. By Nakayama's lemma \cite[Theorem~2.2]{Mat89}, $K=0$, that is, $J^c\subset\fm R'$. By the Artin--Rees lemma \cite[Theorem~8.5]{Mat89}, applied to the finite $R$-module $R'$ and its submodule $R$, there is $c'\geq1$ such that $\fm^rR'\cap R\subset\fm^{r-c'}$ for $r>c'$. Hence $R\cap J^{c(c'+2)}\subset R\cap\fm^{c'+2}R'\subset\fm^2$.

Set $N:=N(\iota)$ and $I_k:=R\cap J^k$. The ideal $I_1$ is proper and contains $\fm$, so $I_1=\fm$; moreover $I_{k+1}=I_k\cap J^{k+1}$. Since $\sigma_{R'}(J^k)\subset J^k$, we have $\sigma_R(I_k)\subset I_k$. The images $V_k$ of the ideals $I_k$ in $V:=\fm/\fm^2$ form a $\overline\sigma_R$-invariant filtration
\[
V=V_1\supset V_2\supset\dots\supset V_{N+1}=0,
\]
so every eigenvalue of $\overline\sigma_R$ is an eigenvalue of the map induced on some $V_k/V_{k+1}$, with $1\leq k\leq N$. Equivariantly, $V_k/V_{k+1}$ is a quotient of $I_k/I_{k+1}$, which embeds in $J^k/J^{k+1}$.
\end{proof}

In the proof of Theorem~B, the map $\gr^k\sigma_{R'}$ acts on $J^k/J^{k+1}=\bigoplus_i\fm_i^k/\fm_i^{k+1}$ as an operator of the following form.

\begin{lemma}\label{lem:block-operator}
Let $\omega$ be a self-map of a finite set $S$, let $(W_s)_{s\in S}$ be finite-dimensional vector spaces and $T_s\colon W_{\omega(s)}\to W_s$ linear maps, and let $\Psi$ be the endomorphism of $W:=\bigoplus_{s\in S}W_s$ given by $\Psi(w)_s:=T_s(w_{\omega(s)})$. For a periodic point $s$ of $\omega$, of period $l$, set $P_s:=T_sT_{\omega(s)}\cdots T_{\omega^{l-1}(s)}\colon W_s\to W_s$.
\begin{enumerate}[label=\textup{(\roman*)}]
\item A number $\mu\neq0$ is an eigenvalue of $\Psi$ if and only if $\mu^l$ is an eigenvalue of $P_s$ for some periodic point $s$ of period $l$.
\item If every periodic point of $\omega$ is fixed, then
\[
\det(t-\Psi)=t^{d}\prod_{\omega(s)=s}\det(t-T_s),\qquad d:=\sum_{\omega(s)\neq s}\dim W_s.
\]
\end{enumerate}
\end{lemma}

\begin{proof}
The subspace $W''$ of $W$ spanned by the spaces $W_s$ with $s$ not periodic is $\Psi$-invariant, because $\omega(s)$ is periodic when $s$ is. Since $\omega^{|S|}(S)$ consists of periodic points, $\Psi^{|S|}$ vanishes on $W''$. The quotient $W/W''$ is the direct sum of the spaces $W_C:=\bigoplus_{s\in C}W_s$ over the cycles $C$ of $\omega$. The map induced by $\Psi$ preserves each $W_C$, and its $l$-th power on $W_C$ is $\bigoplus_{s\in C}P_s$, where $l:=|C|$. This proves~(ii) and the direct implication of~(i).

Conversely, let $s_0$ be a periodic point of period $l$, set $s_i:=\omega^i(s_0)$, with indices modulo $l$, and let $v\neq0$ satisfy $P_{s_0}v=\mu^lv$, where $\mu\neq0$. Set $w_{s_0}:=v$ and $w_{s_i}:=\mu^{i-l}T_{s_i}T_{s_{i+1}}\cdots T_{s_{l-1}}v$ for $0<i<l$. Then $T_{s_i}w_{s_{i+1}}=\mu w_{s_i}$ for every $i$. So the class of $w:=(w_{s_i})_i\in W_C$ in $W/W''$ is an eigenvector, with eigenvalue $\mu$, of the map induced by $\Psi$.
\end{proof}

\begin{proof}[Proof of Theorem~B]
Set $R:=\cO_{X,p}$, $R':=\prod_{y\in S}\cO_{Y,y}$ and $J:=\prod_{y\in S}\fm_y$. The ring $R'$ is the stalk at $p$ of $\pi_*\cO_Y$ \cite[Chapter~2, \S\S3.2--3.3]{GrRe84}, a finite $R$-module by the Finite Mapping Theorem \cite[Chapter~3, \S1.3]{GrRe84}. The homomorphism $\iota:=\pi^*\colon R\to R'$ is injective. Indeed, if $\pi^*h=0$, then $h$ vanishes on $\pi(U)$ for some neighbourhood $U$ of $S$; as $\pi$ is closed and surjective, $\pi(U)$ is a neighbourhood of $p$, and $h=0$ because $X$ is reduced. For $y\in S$, let $g_y^*\colon\cO_{Y,g(y)}\to\cO_{Y,y}$ be the homomorphism induced by $g$. Set $\sigma_R:=f^*$ and $\sigma_{R'}(b)_y:=g_y^*(b_{g(y)})$. The relation $\pi\circ g=f\circ\pi$ gives $\sigma_{R'}\circ\iota=\iota\circ\sigma_R$, and $D_pf$ is the transpose of $\overline\sigma_R$. As $R\cap J^k$ consists of the germs $h$ with $\ord_S(\pi^*h)\geq k$, the integer $N(\iota)$ of Lemma~\ref{lem:transfer-filtration} is $N_S(\pi,p)$. By this lemma, every eigenvalue $\lambda$ of $D_pf$ is an eigenvalue of $\gr^k\sigma_{R'}$ for some $1\leq k\leq N_S(\pi,p)$. Since $J^k/J^{k+1}=\bigoplus_{y\in S}\fm_y^k/\fm_y^{k+1}$, the map $\gr^k\sigma_{R'}$ is the operator $\Psi$ of Lemma~\ref{lem:block-operator}, with $\omega:=g|_S$, $W_y:=\fm_y^k/\fm_y^{k+1}$ and $T_y:=\gr^kg_y^*$.

(1) Assume that every periodic point of $g$ in $S$ is fixed. By Lemma~\ref{lem:block-operator}(ii), either $\lambda=0$ and some point of $S$ is not fixed by $g$, or $\lambda$ is an eigenvalue of $\gr^kg_y^*$ for a fixed point $y\in S$. In the latter case, $\lambda$ is a product of $k$ eigenvalues of $D_yg$ by~\eqref{eq:graded-eigenvalues}. As $N_S(\pi,p)\leq N(\pi,p)$, this proves~(1).

(2) Assume that $X$ is locally irreducible at $p$, and let $y\in S$ be a fixed point of $g$. The space $Y$ is irreducible, hence pure dimensional \cite[Chapter~9, \S1.3]{GrRe84}, so that an irreducible component $Z$ of the germ $(Y,y)$ has dimension $\dim Y=\dim X$. Its image $\pi(Z)$ is an irreducible germ of dimension $\dim X$ \cite[Chapter~5, \S4.1]{GrRe84}, contained in the irreducible germ $(X,p)$, hence equal to it \cite[Chapter~5, \S3.2]{GrRe84}. As above, it follows that $\iota_y:=\pi^*\colon\cO_{X,p}\to\cO_{Y,y}$ is injective. It is finite, since $\cO_{Y,y}$ is a quotient of $R'$, and $N(\iota_y)=N_y(\pi,p)$. Lemma~\ref{lem:transfer-filtration}, applied to $\iota_y$, $f^*$ and $g_y^*$, and~\eqref{eq:graded-eigenvalues} show that every eigenvalue of $D_pf$ is a product of $k$ eigenvalues of $D_yg$, with $k\leq N_y(\pi,p)\leq N(\pi,p)$.
\end{proof}

Let $V_k$ be the image in $\fm_p/\fm_p^2$ of the ideal of germs $h$ with $\ord_S(\pi^*h)\geq k$, respectively $\ord_y(\pi^*h)\geq k$, and let $\lambda$ be an eigenvalue of the map induced by $f^*$ on $V_k/V_{k+1}$. Under the assumption of~(1), if $\lambda\neq0$ or if $g$ fixes every point of $S$, respectively under that of~(2), $\lambda$ is a monomial of degree exactly $k$. We will see later in Example~\ref{ex:transfer}(a) and~(d) that the degree bound is sharp in general, and in Example~\ref{ex:transfer}(b) that local irreducibility is needed.

\subsection{Multipliers at periodic points}

\begin{corollary}\label{cor:transfer-multipliers}
Let $\pi$, $f$ and $g$ be as in Theorem~B, let $p$ be a periodic point of $f$ of period $a$, and $N:=N(\pi,p)$.
\begin{enumerate}[label=\textup{(\arabic*)}]
\item For every nonzero multiplier $\lambda$ of $f$ at $p$, there is a periodic point $y\in\pi^{-1}(p)$ of $g$, of period $b$, such that $\lambda^{b/a}$ is a monomial of degree between $1$ and $N$ in the multipliers of $g$ at $y$.
\item If $X$ is locally irreducible at $p$, then for every multiplier $\lambda$ of $f$ at $p$ and every periodic point $y\in\pi^{-1}(p)$ of $g$, of period $b$, the number $\lambda^{b/a}$ is a monomial of degree between $1$ and $N$ in the multipliers of $g$ at $y$.
\end{enumerate}
\end{corollary}

\begin{proof}
(1) The multiplier $\lambda$ is an eigenvalue of $D_pf^a$. The first paragraph of the proof of Theorem~B, applied to $f^a$ and $g^a$, shows that $\lambda$ is an eigenvalue of the operator $\Psi$ of Lemma~\ref{lem:block-operator} with $\omega:=g^a|_S$ and $T_y:=\gr^k(g^a)_y^*$, for some $1\leq k\leq N$. By Lemma~\ref{lem:block-operator}(i), $\lambda^l$ is an eigenvalue of $P_y=\gr^k(g^{al})_y^*$ for a point $y\in S$ of period $l$ under $g^a$. If $b$ is the period of $y$ under $g$, then $l=b/a$, and $\lambda^{b/a}$ is a product of $k$ eigenvalues of $D_yg^b$ by~\eqref{eq:graded-eigenvalues}.

(2) Apply Theorem~B(2) to $f^b$ and $g^b$.
\end{proof}

\begin{remark}\label{rem:transfer-general}
Irreducibility is used only in the proof of~(2). Theorem~B(1) and Corollary~\ref{cor:transfer-multipliers}(1) hold, with $N=N_S(\pi,p)$, for every finite holomorphic map $\pi\colon Y\to X$ of reduced complex spaces which is surjective onto a neighbourhood of $p$. Theorem~B(2) then holds, with $N=N_y(\pi,p)$, at every fixed point $y\in\pi^{-1}(p)$ with $\dim_yY=\dim_pX$, if $X$ is locally irreducible at $p$. If $X$ and $Y$ have finitely many irreducible components and $f$ and $g$ are finite and surjective, suitable iterates $f^r$ and $g^r$ fix all components. If moreover $p$ is a fixed point of $f^r$ at which $X$ is locally irreducible, Theorem~B applies to $\pi|_{Y_0}\colon Y_0\to X_0$, $f^r$ and $g^r$, where $X_0$ is the component through $p$, for which $T_pX_0=T_pX$, and $Y_0$ is a component with $\pi(Y_0)=X_0$. When several components of $X$ pass through $p$, Theorem~B(1) is applied to $\pi$ itself, as in the second sentence. Let $\{w=0\}$ and $\{w=u(z)\}$ be smooth hypersurface germs at $p$, with $\ord_0u=d\geq2$, so that they are tangent at $p$, each invariant under a holomorphic germ $f$. The class of $w$ spans $V_d$ and $V_{d+1}=0$, so the proof of Theorem~B(1) gives \cite[Proposition~2.2]{Le23}: the eigenvalue of $D_pf$ in the normal direction is a product of exactly $d$ eigenvalues of $D_pf$ on the common tangent space.
\end{remark}

\begin{corollary}\label{cor:descent}
Let $\pi$, $f$ and $g$ be as in Theorem~B, and let $p$ be a periodic point of $f$.
\begin{enumerate}[label=\textup{(\alph*)}]
\item If every periodic point of $g$ in $\pi^{-1}(p)$ is super-saddle, then so is $p$.
\item If $X$ is locally irreducible at $p$ and some periodic point of $g$ in $\pi^{-1}(p)$ is super-saddle, respectively repelling, then so is $p$.
\end{enumerate}
\end{corollary}

\begin{proof}
Let $a$ be the period of $p$, and $\lambda\neq0$ a multiplier of $f$ at $p$. By Corollary~\ref{cor:transfer-multipliers}, $\lambda^{b/a}$ is a monomial of positive degree in the multipliers of $g$ at a periodic point $y\in\pi^{-1}(p)$ of period $b$; in case~(b), $y$ can be any periodic point of $g$ in $\pi^{-1}(p)$. This monomial is nonzero, so its factors are nonzero, hence of modulus greater than $1$. Thus $|\lambda|>1$. In the repelling case of~(b), every multiplier of $f$ at $p$, zero included, is a monomial of positive degree in numbers of modulus greater than $1$, hence is nonzero.
\end{proof}

\subsection{Transport of post-critical finiteness}
We use the ramification divisor $R_h$, the set $\Brdiv(h)$ and the divisorial post-critical set $\PCdiv(g)$ of Section~\ref{ssec:branch-boundary}. We will now show that post-critical finiteness is preserved under semiconjugacies via finite holomorphic maps under a certain dynamical assumption on the branch locus of the semiconjugacy maps.

\begin{lemma}\label{lem:pcf-transport}
Let $\pi\colon Y\to X$, $f\colon X\to X$ and $g\colon Y\to Y$ be finite surjective morphisms of varieties with $\pi\circ g=f\circ\pi$. If $X$ and $Y$ are normal, then:
\begin{enumerate}[label=\textup{(\alph*)}]
\item $R_g+g^*R_\pi=R_\pi+\pi^*R_f$;
\item we have
\[
\PCdiv(g)\subset\pi^{-1}\Bigl(\PCdiv(f)\cup\bigcup_{n\geq1}f^n(\Brdiv(\pi))\Bigr),\qquad
\PCdiv(f)\subset\pi\bigl(\PCdiv(g)\bigr)\cup\bigcup_{n\geq0}f^n(\Brdiv(\pi));
\]
\item if $\bigcup_{n\geq0}f^n(\Brdiv(\pi))$ is a finite union of prime divisors, then $g$ is divisorially PCF if and only if $f$ is.
\end{enumerate}
If $X$ and $Y$ are projective, then:
\begin{enumerate}[label=\textup{(\alph*)},resume]
\item if $f^*L\simeq L^{\otimes M}$ for an ample line bundle $L$, then $\pi^*L$ is ample and $g^*\pi^*L\simeq(\pi^*L)^{\otimes M}$;
\item let $Q$ be the Zariski closure of $\bigcup_{n\geq0}f^n(\Br(\pi))$, and assume $Q\neq X$. If $B_X$ is a branch boundary for $f$, then $\pi^{-1}(B_X\cup Q)$ is a branch boundary for $g$. If $B_Y$ is a branch boundary for $g$, then $\pi(B_Y)\cup Q$ is a branch boundary for $f$.
\end{enumerate}
\end{lemma}

The assumption of~(c) holds in two cases. If $\pi$ is quasi-\'etale, that is $R_\pi=0$, then $\Brdiv(\pi)$ is empty. Assume instead that $\pi$ is the quotient map by a finite group $\Gamma$ of automorphisms of $Y$, and that there is an automorphism $\chi$ of $\Gamma$ such that $g\circ\gamma=\chi(\gamma)\circ g$ for every $\gamma\in\Gamma$. As $e_\pi(E)$ is the order of the subgroup of $\Gamma$ fixing the prime divisor $E$ pointwise, and $\chi(\gamma)\neq\Id$ fixes $g(E)$ pointwise when $\gamma\neq\Id$ fixes $E$ pointwise, we have $g(\Supp R_\pi)\subset\Supp R_\pi$; so $f(\Brdiv(\pi))\subset\Brdiv(\pi)$, and the union in~(c) is $\Brdiv(\pi)$. Moreover $\pi$ is locally biholomorphic at $y$ if and only if the stabilizer of $y$ in $\Gamma$ is trivial; as $\chi$ maps the stabilizer of $y$ into that of $g(y)$, we get $f(\Br(\pi))\subset\Br(\pi)$ and $Q=\Br(\pi)$ in~(e).

\begin{proof}
(a) Ramification indices are multiplicative in towers of discrete valuation rings. Applied to $g^*\circ\pi^*=\pi^*\circ f^*$ at the generic points of a prime divisor $E\subset Y$ and of its images, this gives
\[
e_g(E)\,e_\pi\bigl(g(E)\bigr)=e_\pi(E)\,e_f\bigl(\pi(E)\bigr),
\]
which is the equality of the coefficients of $E$ in the two sides of~(a), up to adding $1$.

(b) By this identity, a prime divisor $E\subset\Supp R_g$ satisfies $E\subset\Supp R_\pi$ or $\pi(E)\subset\Supp R_f$. As $g^n(\mathcal J)\subset\pi^{-1}(f^n(\pi(\mathcal J)))$ for every $\mathcal J\subset Y$, this gives the first inclusion. Conversely, a prime divisor $D\subset\Supp R_f$ is $\pi(E)$ for a prime divisor $E$, and then $E\subset\Supp R_g$ or $g(E)\subset\Supp R_\pi$. Since $f(D)=\pi(g(E))$, we infer that $f(\Supp R_f)\subset\pi(g(\Supp R_g))\cup\Brdiv(\pi)$, and applying $f^{n-1}$, with $f^{n-1}\circ\pi=\pi\circ g^{n-1}$, gives the second inclusion.

(c) The set $Q_1:=\bigcup_{n\geq0}f^n(\Brdiv(\pi))$ is closed. If $f$ is divisorially PCF, then, by the first inclusion of~(b), every prime divisor $g^n(E)$ with $E\subset\Supp R_g$ lies in the proper closed subset $\pi^{-1}(\PCdiv(f)\cup Q_1)$, hence is one of its finitely many irreducible components. The converse follows in the same way from the second inclusion, $\pi$ being closed.

(d) We have $g^*\pi^*L\simeq\pi^*f^*L\simeq(\pi^*L)^{\otimes M}$, and $\pi^*L$ is ample since $\pi$ is finite.

(e) The set $Q$ is closed and algebraic, as $\Br(\pi)$ is, and $f(Q)\subset Q$. Let $B_X$ be a branch boundary for $f$. The set $B_Y:=\pi^{-1}(B_X\cup Q)$ is a proper closed algebraic subset of $Y$ with $g(B_Y)\subset B_Y$. Let $y\notin B_Y$, $n\geq1$ and $y'\in g^{-n}(y)$. Then $\pi$ is locally biholomorphic at $y$, and $\pi(y)$ is a smooth point, hence so is $y$. Moreover $f^n$ is locally biholomorphic at $\pi(y')$, and $\pi(y')\notin\Br(\pi)$, for otherwise $\pi(y)\in Q$. So $\pi$ is locally biholomorphic at $y'$, and $g^n=(\pi|_V)^{-1}\circ f^n\circ\pi$ near $y'$, for a small neighbourhood $V$ of $y$, and $\Br(g^n)\subset B_Y$. The second assertion is proved in the same way: $B_X:=\pi(B_Y)\cup Q$ is proper, closed and $f$-invariant, every point over $x\notin B_X$ is a smooth point at which $\pi$ is locally biholomorphic, so that $x$ is smooth. For $x'\in f^{-n}(x)$ and $y'\in\pi^{-1}(x')$, the map $\pi$ is locally biholomorphic at $y'$, as $x\notin Q$, and so is $g^n$, as $g^n(y')\notin B_Y$; as $\pi$ is also locally biholomorphic at $g^n(y')$, which lies over $x$, the map $f^n=\pi\circ g^n\circ(\pi|_{V'})^{-1}$ is locally biholomorphic at $x'$.
\end{proof}

\section{Rotation discs}\label{sec:rotation}

We first exclude rotation discs for maps with a branch boundary, then prove the boundary theorem for closed invariant sets.

\begin{proposition}\label{prop:key}
Let $g$ be a polarized endomorphism of a normal projective variety $X$, let $B$ be a branch boundary for $g$, and let $\lambda = e^{2 \pi i t}\in S^1, t \notin \QQ$. There is no nonconstant holomorphic map $\varphi\colon D(0,\varepsilon)\to X$, with $\varepsilon>0$, such that
\[
g\bigl(\varphi(z)\bigr)=\varphi(\lambda z)\quad\text{for }z\in D(0,\varepsilon)\qquad\text{and}\qquad\varphi\bigl(D^*(0,\varepsilon)\bigr)\subset X\setminus B.
\]
\end{proposition}

Throughout the section, $X$, $g$, $B$ and $\lambda$ are as in Proposition~\ref{prop:key}, with $g^*H\simeq H^{\otimes M}$ and $\Xo=X\setminus B$, and we assume, by contradiction, that there is a map $\varphi$ as in the proposition. Then $\dim X\geq1$, and we use the notation of Section~\ref{ssec:cone}.

\subsection{Lifting the disc to the cone}\label{ssec:lift-cone}

\begin{lemma}\label{lem:lift}
There is a nonconstant holomorphic map $\hvarphi\colon D(0,\varepsilon)\to\hX\setminus\{0\}$ with $\varpi\circ\hvarphi=\varphi$, $\hvarphi(D^*(0,\varepsilon))\subset\hXo$ and
\begin{equation}\label{eq:equivariant}
\hat g\bigl(\hvarphi(z)\bigr)=\hvarphi(\lambda z)\qquad(z\in D(0,\varepsilon)).
\end{equation}
\end{lemma}

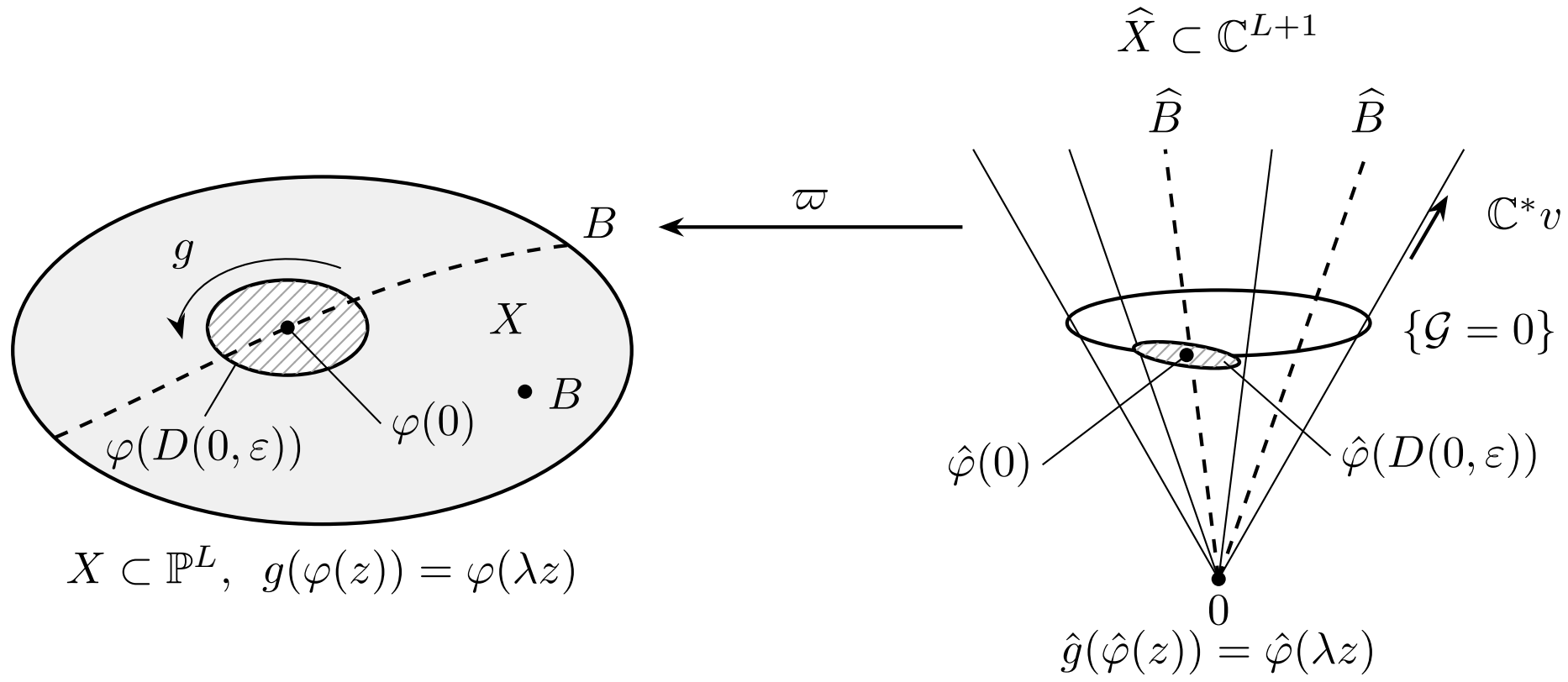
\begin{figure}[ht]
\centering
\begin{tikzpicture}[x=0.84cm,y=0.72cm,>={Stealth[length=2mm]},font=\small]
  \draw[thick,fill=gray!12] (0,0) ellipse [x radius=2.9,y radius=1.9];
  \node at (1.75,0.35) {$X$};
  \draw[thick,fill=white] (-0.325,0.25) ellipse [x radius=0.75,y radius=0.52];
  \draw[pattern=north east lines,pattern color=gray!70] (-0.325,0.25) ellipse [x radius=0.75,y radius=0.52];
  \draw[dashed,thick] (-2.5,-0.95) .. controls (-1.2,-0.3) and (0.4,0.9) .. (2.3,1.15);
  \node at (2.6,1.4) {$B$};
  \fill (1.9,-0.45) circle [radius=1.6pt];
  \node[right] at (1.95,-0.47) {$B$};
  \fill (-0.325,0.25) circle [radius=1.6pt];
  \draw[thin] (-0.325,0.25) -- (0.55,-0.8) node[right,inner sep=2pt] {$\varphi(0)$};
  \draw[thin] (-0.81,-0.15) -- (-1.1,-0.72) node[below,inner sep=2pt] {$\varphi(D(0,\varepsilon))$};
  \draw[->] (0.175,0.90) arc [start angle=60,end angle=190,x radius=1.0,y radius=0.75];
  \node at (-1.3,1.05) {$g$};
  \node at (0,-2.4) {$X\subset\PP^L$, \ $g(\varphi(z))=\varphi(\lambda z)$};

  \coordinate (O) at (8.4,-2.5);
  \foreach \xt in {6.1,7.0,8.9,10.7} {\draw (O) -- (\xt,2.2);}
  \draw[dashed,thick] (O) -- (7.9,2.2);
  \draw[dashed,thick] (O) -- (9.8,2.2);
  \fill (O) circle [radius=1.6pt] node[below] {$0$};
  \node[above] at (7.9,2.2) {$\hB$};
  \node[above] at (9.8,2.2) {$\hB$};
  \draw[->,thick] (10.2,1.0) -- (10.55,1.7);
  \node[right] at (10.75,1.5) {$\CC^*v$};
  \draw[thick] (8.4,0.3) ellipse [x radius=1.42,y radius=0.36];
  \node[right] at (9.95,0.2) {$\{\cG=0\}$};
  \draw[thick,fill=white,rotate around={-8:(8.1,-0.05)}] (8.1,-0.05) ellipse [x radius=0.5,y radius=0.12];
  \draw[pattern=north east lines,pattern color=gray!70,rotate around={-8:(8.1,-0.05)}] (8.1,-0.05) ellipse [x radius=0.5,y radius=0.12];
  \fill (8.1,-0.05) circle [radius=1.6pt];
  \draw[thin] (8.1,-0.05) -- (6.75,-1.25) node[left,inner sep=2pt] {$\hvarphi(0)$};
  \draw[thin] (8.45,-0.12) -- (9.45,-1.15) node[right,inner sep=2pt] {$\hvarphi(D(0,\varepsilon))$};
  \node at (8.4,3.5) {$\hX\subset\CC^{L+1}$};
  \node at (8.4,-3.3) {$\hat g(\hvarphi(z))=\hvarphi(\lambda z)$};

  \draw[->,thick] (6.0,1.35) -- (3.15,1.35) node[midway,above] {$\varpi$};
\end{tikzpicture}
\caption{The equivariant lift of a rotation disc to the cone over $X$ (Lemma~\ref{lem:lift}). The centre of the disc may lie on the branch boundary $B$ (dashed).}
\label{fig:cone}
\end{figure}

\begin{proof}
The map $\varpi$ is a $\CC^*$-bundle over all of $X$, and holomorphic line bundles on a disc are trivial, so $\varphi$ has a holomorphic lift $\hvarphi_0\colon D(0,\varepsilon)\to\hX\setminus\{0\}$. The points $\hat g(\hvarphi_0(z))$ and $\hvarphi_0(\lambda z)$ lie in the same fibre of $\varpi$. Hence $\hat g(\hvarphi_0(z))=a(z)\,\hvarphi_0(\lambda z)$, where $a$ is nowhere zero and holomorphic, being locally a quotient of coordinates. Write $a=\exp A$ with $A(z)=\sum_{m\geq0}A_mz^m$. Since $|\lambda|=1$ and $M\geq2$, we have $|\lambda^m-M|\geq M-1\geq1$ for every $m\geq0$. The series
\[
\mathcal U(z):=\sum_{m\geq0}\frac{A_mz^m}{\lambda^m-M}
\]
therefore converges on $D(0,\varepsilon)$ and satisfies
\[
\mathcal U(\lambda z)-M\,\mathcal U(z)=A(z).
\]
Set $\hvarphi:=\exp(\mathcal U)\,\hvarphi_0$. As $\hat g$ is homogeneous of degree $M$,
\[
\hat g\bigl(\hvarphi(z)\bigr)=\exp\bigl(M\,\mathcal U(z)+A(z)\bigr)\hvarphi_0(\lambda z)=\hvarphi(\lambda z).
\]
Finally $\varpi\circ\hvarphi=\varphi$, so $\hvarphi$ is nonconstant and maps $D^*(0,\varepsilon)$ into $\varpi^{-1}(\Xo)=\hXo$ (Figure~\ref{fig:cone}).
\end{proof}

Write $\hvarphi(z)=\sum_{m\geq0}a_mz^m$ with $a_m\in\CC^{L+1}$, and let $n_0\geq1$ be the greatest common divisor of the integers $m\geq1$ with $a_m\neq0$. Then $\hvarphi(z)=\hat\psi(z^{n_0})$, where $\hat\psi(w):=\sum_ma_{n_0m}w^m$ converges on $D(0,\varepsilon^{n_0})$. The map $\hat\psi$ is nonconstant, takes values in $\hX\setminus\{0\}$, maps $D^*(0,\varepsilon^{n_0})$ into $\hXo$, and satisfies $\hat g(\hat\psi(w))=\hat\psi(\lambda^{n_0}w)$, where $\lambda^{n_0}\in S^1$ is not a root of unity. Replacing $(\varphi,\hvarphi,\lambda,\varepsilon)$ by $(\varpi\circ\hat\psi,\hat\psi,\lambda^{n_0},\varepsilon^{n_0})$, we assume from now on that $n_0=1$.

\begin{lemma}\label{lem:RI}
Let $\theta,\theta'\in S^1$. If $\hvarphi(\theta z)=\hvarphi(\theta'z)$ for all $z$ in a nonempty open subset of $D(0,\varepsilon)$, then $\theta=\theta'$.
\end{lemma}

\begin{proof}
By the identity theorem the equality holds on $D(0,\varepsilon)$. Comparing coefficients gives $a_m(\theta^m-\theta'^m)=0$, so $(\theta/\theta')^m=1$ for every $m\geq1$ with $a_m\neq0$. These integers have greatest common divisor $1$, hence $\theta=\theta'$.
\end{proof}

\subsection{The covering attached to the disc}\label{ssec:covering}

\begin{figure}[ht]
\centering
\begin{tikzpicture}[x=0.86cm,y=0.76cm,>={Stealth[length=2mm]},font=\small]
  \draw[thick,fill=gray!12] (1.6,0) circle [radius=1.35];
  \draw[thick,postaction={decorate,decoration={markings,mark=at position 0.3 with {\arrow{Stealth[length=2mm]}}}}]
        (2.3,0) arc [start angle=0,end angle=360,radius=0.7];
  \draw[fill=white] (1.6,0) circle [radius=2pt];
  \fill (2.3,0) circle [radius=1.5pt] node[right,inner sep=2pt] {$z_0$};
  \node at (2.28,0.72) {$\ell$};
  \draw[->] (2.0,-1.75) arc [start angle=-60,end angle=-120,radius=0.8];
  \node at (1.6,-2.2) {$z\mapsto\lambda^{-1}z$};
  \node at (0.35,1.45) {$D^*$};

  \draw[thick] (6.4,-1.55) rectangle (11.4,1.55);
  \draw[thick,fill=gray!12] (8.9,0) ellipse [x radius=1.7,y radius=1.0];
  \draw[line width=1.6pt,postaction={decorate,decoration={markings,mark=at position 0.3 with {\arrow{Stealth[length=2.4mm]}}}}]
        (9.8,0) arc [start angle=0,end angle=360,x radius=0.9,y radius=0.52];
  \fill[white] (8.82,-1.55) rectangle (8.98,1.55);
  \draw[pattern=north east lines,pattern color=gray] (8.82,-1.55) rectangle (8.98,1.55);
  \fill (9.8,0) circle [radius=1.5pt] node[right,inner sep=2pt] {$x_0$};
  \node at (10.08,0.48) {$\gamma_0$};
  \node at (7.75,-1.2) {$\hvarphi(D^*)$};
  \node[above] at (8.9,1.55) {$\hB$};
  \node at (10.95,-1.25) {$\hXo$};

  \foreach \yy in {3.3,5.0,6.7} {
    \draw[thick,fill=white] (5.9,\yy-0.55) -- (10.9,\yy-0.55) -- (11.6,\yy+0.55) -- (6.6,\yy+0.55) -- cycle;
    \fill[white] (8.77,\yy-0.57) -- (8.93,\yy-0.57) -- (9.63,\yy+0.57) -- (9.47,\yy+0.57) -- cycle;
    \draw[pattern=north east lines,pattern color=gray,thin] (8.77,\yy-0.55) -- (8.93,\yy-0.55) -- (9.63,\yy+0.55) -- (9.47,\yy+0.55) -- cycle;
  }
  \draw[decorate,decoration={coil,aspect=0.5,segment length=6.6mm,amplitude=1.4mm}] (5.95,2.75) -- (5.95,7.3);
  \draw[thick,fill=gray!12] (9.2,5.0) ellipse [x radius=1.25,y radius=0.36];
  \draw[line width=1.2pt] (9.2,5.0) ellipse [x radius=0.72,y radius=0.2];
  \draw[fill=white] (9.2,5.0) circle [radius=1.6pt];
  \fill (9.92,5.0) circle [radius=1.5pt];
  \draw[thin] (9.92,5.0) -- (10.55,5.0);
  \node[right,inner sep=1pt] at (10.55,5.0) {$\tilde x_0$};
  \node at (7.1,5.0) {$\kappa_0(D^*)$};
  \draw[->] (9.81,5.6) arc [start angle=50,end angle=130,x radius=0.95,y radius=0.3];
  \node at (10.2,5.82) {$G$};
  \fill (10.9,6.75) circle [radius=1.5pt];
  \fill (10.75,3.3) circle [radius=1.5pt];
  \draw[->] (10.9,6.75) .. controls (12.9,6.2) and (12.9,3.9) .. (10.83,3.32);
  \node at (12.75,5.0) {$G$};
  \node at (6.9,7.65) {$\tX$};

  \draw[->,thick] (3.15,0) -- (6.25,0) node[midway,above] {$\hvarphi$};
  \draw[->,thick] (10.3,2.6) -- (10.3,1.7) node[midway,right] {$\pi$};
  \draw[->,thick] (2.75,1.0) .. controls (4.6,1.6) and (5.0,4.4) .. (5.6,4.9);
  \node at (5.05,2.25) {$\kappa_0$};

  \draw (-0.3,2.75) rectangle (4.15,6.6);
  \draw[thick] (0.0,3.25) .. controls (1.2,4.05) and (2.6,4.05) .. (3.9,3.3);
  \node at (3.6,3.05) {$\hB$};
  \fill[gray!45] (1.9,3.84) .. controls (0.9,4.6) and (1.2,5.5) .. (2.5,6.1)
                 .. controls (2.6,5.2) and (2.9,4.6) .. (1.9,3.84);
  \draw[thick] (1.9,3.84) .. controls (0.9,4.6) and (1.2,5.5) .. (2.5,6.1);
  \draw[thick] (1.9,3.84) .. controls (2.9,4.6) and (2.6,5.2) .. (2.5,6.1);
  \fill (1.9,3.84) circle [radius=1.6pt] node[below right,inner sep=2pt] {$p$};
  \fill (2.5,6.1) circle [radius=1.6pt];
  \node at (0.9,5.1) {$\beta_0$};
  \node at (3.1,5.0) {$\beta_1$};
\end{tikzpicture}
\caption{The covering attached to the punctured disc and the lifted inverse branch $G$ (Lemma~\ref{lem:G}). Inset: two paths from $p=\hvarphi(0)\in\hB$, otherwise contained in $\hXo$, joined by a relative homotopy (Remark~\ref{rem:relative-homotopy}).}
\label{fig:covering}
\end{figure}
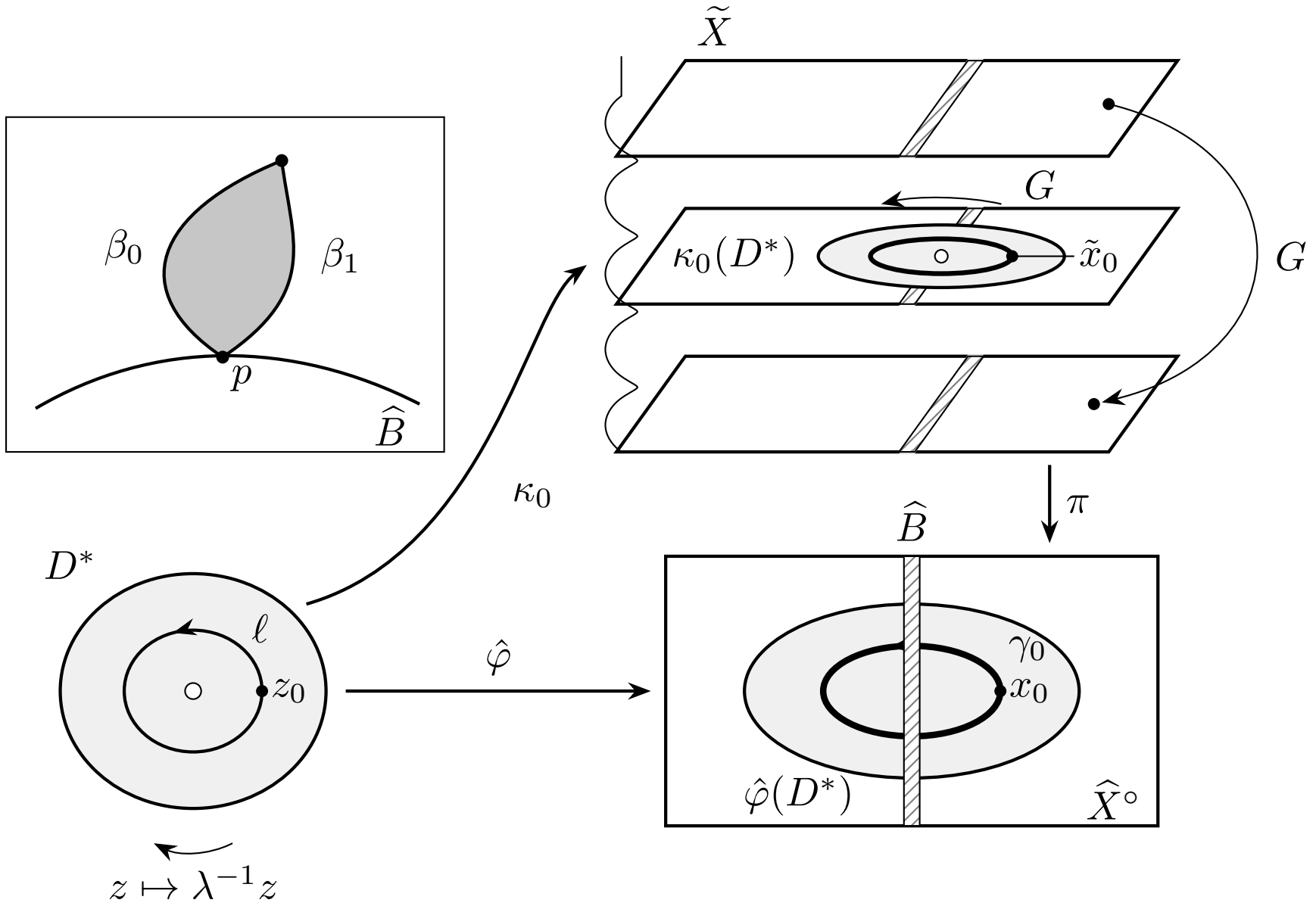

Fix $z_0\in D^*(0,\varepsilon)$ and set $x_0:=\hvarphi(z_0)\in\hXo$. The loop $\ell(t):=z_0\exp(2\pi it)$, $t\in[0,1]$, generates $\pi_1(D^*(0,\varepsilon),z_0)$. Let $\gamma_0:=[\hvarphi\circ\ell]\in\pi_1(\hXo,x_0)$, and let
\[
\pi\colon(\tX,\tilde x_0)\to(\hXo,x_0)
\]
be the connected covering with $\pi_*\pi_1(\tX,\tilde x_0)=\langle\gamma_0\rangle$ \cite[Chapter~III, Theorem~8.1]{Bre93} (Figure~\ref{fig:covering}). We endow $\tX$ with the complex structure for which $\pi$ is locally biholomorphic; it is a connected complex manifold, and it is second countable, as $\pi_1(\hXo)$ is countable. Since $\hvarphi_*\pi_1(D^*(0,\varepsilon),z_0)=\langle\gamma_0\rangle$, the lifting theorem \cite[Chapter~III, Theorem~4.1]{Bre93} gives a unique holomorphic map
\begin{equation}\label{eq:kappa0}
\kappa_0\colon D^*(0,\varepsilon)\to\tX\qquad\text{with}\qquad\pi\circ\kappa_0=\hvarphi|_{D^*(0,\varepsilon)}\quad\text{and}\quad\kappa_0(z_0)=\tilde x_0.
\end{equation}
It is nonconstant. As $\pi_*$ is injective and $\pi_*[\kappa_0\circ\ell]=\gamma_0$, the class $[\kappa_0\circ\ell]$ generates $\pi_1(\tX,\tilde x_0)$. Below, $D^*$ stands for $D^*(0,\varepsilon)$.

\begin{lemma}\label{lem:G}
There is a holomorphic map $G\colon\tX\to\tX$ such that
\begin{equation}\label{eq:G}
\hat g\circ\pi\circ G=\pi\qquad\text{and}\qquad G\bigl(\kappa_0(z)\bigr)=\kappa_0(\lambda^{-1}z)\quad(z\in D^*).
\end{equation}
Consequently $\hat g^m\circ\pi\circ G^m=\pi$ and $G^m\circ\kappa_0(z)=\kappa_0(\lambda^{-m}z)$ for every $m\geq0$.
\end{lemma}

\begin{proof}
By \eqref{eq:equivariant}, $\hat g(\hvarphi(\lambda^{-1}z))=\hvarphi(z)\in\hXo$ for $z\in D^*$. Let $Y_0$ be the connected component of $\hat g^{-1}(\hXo)$ which contains the connected set $\hvarphi(\lambda^{-1}D^*)$, and set $y_0:=\hvarphi(\lambda^{-1}z_0)$, so that $\hat g(y_0)=x_0$. By Lemma~\ref{lem:covering}(1), $Y_0\subset\hXo$; by Lemma~\ref{lem:covering}(2), the restriction $\eta:=\hat g|_{Y_0}\colon Y_0\to\hXo$ of a finite covering to a connected component is a covering.

We first lift $\pi$ through $\eta$. The loop $\hvarphi\circ(\lambda^{-1}\ell)$ lies in $Y_0$ and is based at $y_0$, and its image under $\eta$ is $\hvarphi\circ\ell$. Hence $\gamma_0\in\eta_*\pi_1(Y_0,y_0)$, and there is a unique holomorphic map $\sigma\colon\tX\to Y_0$ with $\eta\circ\sigma=\pi$ and $\sigma(\tilde x_0)=y_0$. The maps $\sigma\circ\kappa_0$ and $\hvarphi(\lambda^{-1}\,\cdot)$ from $D^*$ to $Y_0$ are lifts through $\eta$ of $\hvarphi|_{D^*}$ which agree at $z_0$. Therefore
\begin{equation}\label{eq:hkappa}
\sigma\circ\kappa_0(z)=\hvarphi(\lambda^{-1}z)\qquad(z\in D^*).
\end{equation}

We now lift $\sigma$ through $\pi$. Set $\tilde x_1:=\kappa_0(\lambda^{-1}z_0)$, so that $\pi(\tilde x_1)=y_0=\sigma(\tilde x_0)$. By \eqref{eq:hkappa},
\[
\sigma_*[\kappa_0\circ\ell]=[\hvarphi\circ(\lambda^{-1}\ell)]=\pi_*[\kappa_0\circ(\lambda^{-1}\ell)],
\]
where $\kappa_0\circ(\lambda^{-1}\ell)$ is a loop based at $\tilde x_1$. As $[\kappa_0\circ\ell]$ generates $\pi_1(\tX,\tilde x_0)$, we get $\sigma_*\pi_1(\tX,\tilde x_0)\subset\pi_*\pi_1(\tX,\tilde x_1)$, and there is a unique holomorphic map $G\colon\tX\to\tX$ with $\pi\circ G=\sigma$ and $G(\tilde x_0)=\tilde x_1$. Then $\hat g\circ\pi\circ G=\eta\circ\sigma=\pi$. Moreover $G\circ\kappa_0$ and $\kappa_0(\lambda^{-1}\,\cdot)$ are lifts through $\pi$ of $\hvarphi(\lambda^{-1}\,\cdot)$ which agree at $z_0$, so they coincide. The last assertion follows by induction.
\end{proof}

\begin{remark}\label{rem:relative-homotopy}
If $\varphi(0)\notin B$, the loop $\hvarphi\circ\ell$ is contractible in $\hvarphi(D(0,\varepsilon))\subset\hXo$, so $\pi$ is the universal covering of $\hXo$, as in \cite[Section~3]{Le22}. If $\hX$ is smooth at $\hvarphi(0)$ and $\hB$ is, near this point, a smooth hypersurface transverse to $\hvarphi$, then $\pi$ is the covering induced by \emph{relative homotopy} of \cite[Section~4.3]{Le22}, whose points are the classes of paths issued from $\hvarphi(0)$ and otherwise contained in $\hXo$, up to homotopy through such paths with fixed endpoint. 
\end{remark}

\subsection{Limits of the inverse branches}\label{ssec:limits}

For $m\geq0$ we set $k_m:=\pi\circ G^m\colon\tX\to\hXo\subset\CC^{L+1}$. By Lemma~\ref{lem:G},
\begin{equation}\label{eq:km}
\hat g^m\circ k_m=\pi,\qquad k_{a+a'}=k_a\circ G^{a'}\qquad\text{and}\qquad k_m\circ\kappa_0(z)=\hvarphi(\lambda^{-m}z)\quad(z\in D^*).
\end{equation}
Denote by 
\[
S = \overline{ \{k_m \in \mathcal C(\tX,\CC^{L+1}), m \geq 0  \}} 
\] where the closure is taken with respect to the compact-open topology. An element $h$ of $S$ may meet $\hB$, as $B$ may have components of any codimension; we lift it to $\tX$ on the set
\[
\Omega_h:=\tX\setminus h^{-1}(\hB).
\]

\begin{lemma}\label{lem:S}
\begin{enumerate}[label=\textup{(\arabic*)}]
\item The space $S$ is compact and metrizable. Its elements are holomorphic maps $\tX\to\hX\setminus\{0\}$, and every $h\in S$ is the limit of a sequence $(k_{a_j})$ with integers $a_j\geq0$.
\item Let $h\in S$ and $k_{a_j}\to h$. Then $\lambda^{-a_j}$ converges to some $\theta_h\in S^1$, and
\begin{equation}\label{eq:theta}
h\circ\kappa_0(z)=\hvarphi(\theta_hz)\in\hXo\qquad(z\in D^*).
\end{equation}
The number $\theta_h$ depends only on $h$, and the map $h\mapsto\theta_h$ is continuous on $S$.
\item For every $h\in S$, the set $\Omega_h$ is open, dense and connected in $\tX$, and contains $\kappa_0(D^*)$.
\end{enumerate}
\end{lemma}

\begin{proof}
(1) As in \cite[Proposition~3.9 and Lemma~4.13]{Le22}, we use the Green function. By \eqref{eq:km} and Lemma~\ref{lem:green}, $\cG\circ\pi=\cG\circ\hat g^m\circ k_m=M^m\,\cG\circ k_m$, so that
\[
\bigl|\log|k_m(x)|\bigr|\leq C_0+\bigl|\cG(\pi(x))\bigr|\qquad(x\in\tX,\ m\geq0).
\]
Thus $(k_m)$ is locally uniformly bounded, and locally uniformly bounded away from $0$. Since $\tX$ is second countable, $\mathcal C(\tX,\CC^{L+1})$ is metrizable, and $S$ is compact by Montel's theorem. Its elements are limits of sequences $(k_{a_j})$; they are holomorphic, take values in the closed set $\hX$ and satisfy the same bound, so they do not vanish.

(2) By \eqref{eq:km}, $h(\kappa_0(z))=\lim_j\hvarphi(\lambda^{-a_j}z)$ for $z\in D^*$. If $\theta^*\in S^1$ is a limit point of $(\lambda^{-a_j})$, then $h\circ\kappa_0=\hvarphi(\theta^*\,\cdot)$ on $D^*$. By Lemma~\ref{lem:RI}, two limit points coincide, so $\lambda^{-a_j}$ converges, and its limit $\theta_h$ satisfies \eqref{eq:theta}. The same argument shows that $\theta_h$ does not depend on the sequence $(a_j)$. If $h_i\to h$ in $S$ and $\theta_{h_i}\to\theta^*$ along a subsequence, then
\[
\hvarphi(\theta^*z)=\lim_i h_i(\kappa_0(z))=\hvarphi(\theta_hz).
\]
Lemma~\ref{lem:RI} gives $\theta^*=\theta_h$, so $h\mapsto\theta_h$ is continuous.

(3) As $h$ takes values in $\hX$ and the $Q_i$ are polynomials, $h^{-1}(\hB)=\{Q_1\circ h=\dots=Q_l\circ h=0\}$ is an analytic subset of $\tX$. By \eqref{eq:theta} it does not meet $\kappa_0(D^*)$.
\end{proof}

\begin{lemma}\label{lem:L}
Let $h\in S$, and let $(a_j)$ be a sequence of integers $a_j\geq0$ with $k_{a_j}\to h$. Then $G^{a_j}$ converges locally uniformly on $\Omega_h$ to a holomorphic map $\tilde h\colon\Omega_h\to\tX$ with
\begin{equation}\label{eq:lift}
\pi\circ\tilde h=h|_{\Omega_h}\qquad\text{and}\qquad\tilde h\circ\kappa_0(z)=\kappa_0(\theta_hz)\quad(z\in D^*).
\end{equation}
Moreover $\tilde h$ is the unique lift of $h|_{\Omega_h}$ through $\pi$ with $\tilde h(\tilde x_0)=\kappa_0(\theta_hz_0)$, so it depends only on $h$.
\end{lemma}

\begin{proof}
Let $\Omega'$ be the set of points $x\in\Omega_h$ such that $(G^{a_j}(x))_j$ converges in $\tX$. It contains $\kappa_0(D^*)$, since $G^{a_j}(\kappa_0(z))=\kappa_0(\lambda^{-a_j}z)\to\kappa_0(\theta_hz)$ by Lemmas~\ref{lem:G} and~\ref{lem:S}(2).

Let $x\in\Omega_h$. Choose an open set $V\subset\hXo$ containing $h(x)$ which is evenly covered by $\pi$. Choose a connected open neighbourhood $\Theta$ of $x$ whose closure is a compact subset of $\Omega_h$ with $h(\overline\Theta)\subset V$. As $k_{a_j}\to h$ uniformly on $\overline\Theta$, there is a compact set $V_1\subset V$ with $k_{a_j}(\overline\Theta)\subset V_1$ for $j$ large. Then, for $j$ large, the connected set $G^{a_j}(\Theta)$ lies in a single sheet $\mathcal S_j$ of $\pi^{-1}(V)$.

Suppose that $\Theta$ contains a point $x'\in\Omega'$, and let $w:=\lim_jG^{a_j}(x')$. Then $\pi(w)=h(x')\in V$. Let $\mathcal S$ be the sheet which contains $w$. As $\mathcal S$ is open, $\mathcal S_j=\mathcal S$ for $j$ large, and
\[
G^{a_j}|_\Theta=(\pi|_{\mathcal S})^{-1}\circ k_{a_j}|_\Theta\longrightarrow(\pi|_{\mathcal S})^{-1}\circ h|_\Theta
\]
uniformly on $\Theta$, since $(\pi|_{\mathcal S})^{-1}$ is uniformly continuous on $V_1$. Thus $\Theta\subset\Omega'$ as soon as $\Theta$ meets $\Omega'$. It follows that $\Omega'$ is open and closed in the connected set $\Omega_h$, so that $\Omega'=\Omega_h$, and that the convergence is locally uniform. The limit $\tilde h$ is locally of the form $(\pi|_{\mathcal S})^{-1}\circ h$, hence holomorphic with $\pi\circ\tilde h=h$, and $\tilde h(\kappa_0(z))=\kappa_0(\theta_hz)$. The last assertion follows from the uniqueness of lifts on the connected set $\Omega_h$.
\end{proof}

\subsection{The semigroup of limits and an idempotent}\label{ssec:semigroup}

\begin{lemma}\label{lem:W}
Let $h,h'\in S$, and let $(a_j)$, $(a'_j)$ be sequences of nonnegative integers with $k_{a_j}\to h$ and $k_{a'_j}\to h'$. Then $k_{a_j+a'_j}$ converges in $S$ to an element $h\star h'$ which depends only on $(h,h')$, and
\begin{equation}\label{eq:star}
h\star h'=h\circ\tilde h'\ \text{ on }\Omega_{h'}\qquad\text{and}\qquad\theta_{h\star h'}=\theta_h\theta_{h'}.
\end{equation}
The product $\star$ on $S$ is commutative and associative.
\end{lemma}

\begin{proof}
By \eqref{eq:km}, $k_{a_j+a'_j}=k_{a_j}\circ G^{a'_j}$. As $k_{a_j}\to h$ locally uniformly on $\tX$ and $G^{a'_j}\to\tilde h'$ locally uniformly on $\Omega_{h'}$ (Lemma~\ref{lem:L}), we get $k_{a_j+a'_j}\to h\circ\tilde h'$ locally uniformly on $\Omega_{h'}$. Every limit point of $(k_{a_j+a'_j})$ in the compact metrizable space $S$ is continuous on $\tX$ and equals $h\circ\tilde h'$ on the dense set $\Omega_{h'}$ (Lemma~\ref{lem:S}(3)). So there is exactly one limit point, and the sequence converges. Its limit depends only on $h$ and $\tilde h'$, hence only on $(h,h')$ by Lemma~\ref{lem:L}. Finally $\lambda^{-(a_j+a'_j)}\to\theta_h\theta_{h'}$, which gives the second identity in \eqref{eq:star}.

Commutativity follows from $a_j+a'_j=a'_j+a_j$. If $h''\in S$ and $k_{a''_j}\to h''$, the first part, applied to the pairs of sequences $(a_j+a'_j,a''_j)$ and $(a_j,a'_j+a''_j)$, shows that $(h\star h')\star h''$ and $h\star(h'\star h'')$ are both the limit of $(k_{a_j+a'_j+a''_j})$.
\end{proof}

\begin{lemma}\label{lem:idempotent}
There are $\rho\in S$ and a sequence of integers $n_j\to\infty$ such that
\begin{equation}\label{eq:idempotent}
k_{n_j}\to\rho,\qquad\rho\star\rho=\rho,\qquad\theta_\rho=1,\qquad\rho\circ\kappa_0=\hvarphi|_{D^*}\qquad\text{and}\qquad\cG\circ\rho=0.
\end{equation}
\end{lemma}

\begin{proof}
We follow the extraction in the proof of \cite[Theorem~2.1.29]{Aba89}. By compactness of $S$, there are indices $m_1<m_2<\cdots$ with $m_{j+1}\geq2m_j+j$ such that $k_{m_j}$ converges to some $h\in S$. Set $s_j:=m_{j+1}-m_j$ and $t_j:=s_j-m_j$, so that $s_j\geq t_j\geq j$. Choose an infinite set $J$ of indices along which $k_{s_j}\to\rho$ and $k_{t_j}\to h'$ for some $\rho,h'\in S$. Lemma~\ref{lem:W}, applied along $J$ to the pairs of sequences $(m_j,s_j)$ and $(t_j,m_j)$, gives $h\star\rho=\lim k_{m_{j+1}}=h$ and $h'\star h=\lim k_{s_j}=\rho$. Hence
\[
\rho\star\rho=(h'\star h)\star\rho=h'\star(h\star\rho)=h'\star h=\rho.
\]
By \eqref{eq:star}, $\theta_\rho^2=\theta_\rho$, so $\theta_\rho=1$, and \eqref{eq:theta} gives $\rho\circ\kappa_0=\hvarphi$ on $D^*$. Let $(n_j)$ be the sequence $(s_j)_{j\in J}$. Finally, for $x\in\tX$ we have $\cG(k_{n_j}(x))=M^{-n_j}\cG(\pi(x))\to0$ and $k_{n_j}(x)\to\rho(x)\neq0$, so $\cG(\rho(x))=0$ by continuity of $\cG$ on $\CC^{L+1}\setminus\{0\}$.
\end{proof}

\subsection{A closed retract on which the inverse branch is an automorphism}\label{ssec:retract}

Set $\Omega:=\Omega_\rho=\tX\setminus\rho^{-1}(\hB)$, and let $\tilde\rho\colon\Omega\to\tX$ be the map given by Lemma~\ref{lem:L} for $h=\rho$ and the sequence $(n_j)$. By \eqref{eq:lift}, \eqref{eq:star} and \eqref{eq:idempotent},
\begin{equation}\label{eq:retraction}
G^{n_j}\to\tilde\rho\ \text{ locally uniformly on }\Omega,\qquad\pi\circ\tilde\rho=\rho|_\Omega,\qquad\tilde\rho\circ\kappa_0=\kappa_0,\qquad\rho\circ\tilde\rho=\rho\ \text{ on }\Omega.
\end{equation}

\begin{proposition}\label{prop:retract}
The map $\tilde\rho$ is a holomorphic retraction of $\Omega$. Its image $\tR:=\tilde\rho(\Omega)$ is a closed connected complex submanifold of $\tX$ of positive dimension, which contains $\kappa_0(D^*)$. Moreover $\pi=\rho$ on $\tR$, and $\pi(\tR)$ is contained in the compact subset $\hX\cap\{\cG=0\}$ of $\CC^{L+1}\setminus\{0\}$.
\end{proposition}

\begin{proof}
Let $x\in\Omega$. By \eqref{eq:retraction}, $\rho(\tilde\rho(x))=\rho(x)\notin\hB$, so $\tilde\rho(x)\in\Omega$. The maps $\tilde\rho\circ\tilde\rho$ and $\tilde\rho$ are lifts of $\rho|_\Omega$ through $\pi$ on the connected set $\Omega$, since $\pi\circ\tilde\rho\circ\tilde\rho=\rho\circ\tilde\rho=\rho$, and they agree on $\kappa_0(D^*)$. Hence $\tilde\rho\circ\tilde\rho=\tilde\rho$. By Lemma~\ref{lem:rossi}, the set $\tR$ is a closed complex submanifold of $\Omega$. It is connected, and it contains the image of the nonconstant map $\kappa_0$ by \eqref{eq:retraction}, so $\dim\tR\geq1$.

We now show that $\tR$ is closed in $\tX$. The set $\cZ:=\{x\in\tX:\rho(x)=\pi(x)\}$ is closed in $\tX$, and $\cZ\subset\Omega$ because $\pi$ takes values in $\hXo$. For $x\in\tR$ we have $\rho(x)=\pi(\tilde\rho(x))=\pi(x)$, so $\tR\subset\cZ$. Hence the closure of $\tR$ in $\tX$ lies in $\Omega$, where $\tR$ is closed.

Finally, for $x\in\tR$ we have $\pi(x)=\rho(x)$ and $\cG(\rho(x))=0$ by \eqref{eq:idempotent}. We conclude by Lemma~\ref{lem:green}, as for the centre manifold in \cite[Section~3]{Le22}.
\end{proof}

In particular, $K\cap\tR$ is compact for every compact set $K\subset\tX$, and $\tR$ is second countable.

\begin{proposition}\label{prop:automorphism}
We have $G(\tR)\subset\tR$, and $G_{\mathcal R}:=G|_{\tR}$ is a holomorphic automorphism of $\tR$. Moreover
\begin{equation}\label{eq:GM}
\pi\circ G_{\mathcal R}^{-1}=\hat g\circ\pi\ \text{ on }\tR,\qquad G_{\mathcal R}^m\circ\kappa_0(z)=\kappa_0(\lambda^{-m}z)\quad(m\in\ZZ,\ z\in D^*),\qquad G_{\mathcal R}^{n_j}\to\Id_{\tR}
\end{equation}
locally uniformly on $\tR$.
\end{proposition}

\begin{proof}
By Lemma~\ref{lem:G}, $\hat g\circ k_m\circ G=\pi\circ G^m=k_m$, and letting $m=n_j\to\infty$ gives $\hat g\circ\rho\circ G=\rho$. If $x\in\Omega$, then $\hat g(\rho(G(x)))=\rho(x)\notin\hB$, so $\rho(G(x))\notin\hB$ by Lemma~\ref{lem:covering}(1). Thus $G(\Omega)\subset\Omega$, and by \eqref{eq:retraction} and the continuity of $G$,
\[
\tilde\rho\bigl(G(x)\bigr)=\lim_jG^{n_j}\bigl(G(x)\bigr)=G\bigl(\tilde\rho(x)\bigr)\qquad(x\in\Omega).
\]
For $x\in\tR$ it follows that $G(x)=\tilde\rho(G(x))\in\tR$. The map $G_{\mathcal R}$ is holomorphic, and $G_{\mathcal R}^{n_j}\to\tilde\rho|_{\tR}=\Id_{\tR}$ locally uniformly by \eqref{eq:retraction}.

The map $G_{\mathcal R}$ is injective: if $G(x)=G(y)$ with $x,y\in\tR$, then $G^{n_j}(x)=G^{n_j}(y)$ for all $j$, and $x=y$ in the limit. Its differential is injective at every point $x$, because
\[
D_x(G_{\mathcal R}^{n_j})
 =D_{G(x)}(G_{\mathcal R}^{n_j-1})\circ D_xG_{\mathcal R}
\]
tends to the identity.

To prove surjectivity, let $y\in\tR$, and let $u\colon U\to\CC^{\dim\tR}$ be a chart of $\tR$ with $u(y)=0$ and $u(U)\supset\overline b(0,2r)$, a closed ball. For $j$ large, $G^{n_j}$ maps $u^{-1}(\overline b(0,r))$ into $u^{-1}(b(0,2r))$, and $\hat G_j:=u\circ G^{n_j}\circ u^{-1}$ satisfies $|\hat G_j(w)-w|<r$ on $\overline b(0,r)$. The map $w\mapsto w-\hat G_j(w)$ sends $\overline b(0,r)$ continuously into itself, so it has a fixed point $w_0$ by Brouwer's theorem \cite[Chapter~II, Corollary~11.12]{Bre93}. Then $\hat G_j(w_0)=0$, and $y=G^{n_j}(u^{-1}(w_0))\in G(\tR)$.

So $G_{\mathcal R}$ is a bijective local biholomorphism of $\tR$, hence an automorphism. The first identity of \eqref{eq:GM} follows from \eqref{eq:G}, and the second from Lemma~\ref{lem:G}, applied at the point $\lambda z$ when $m<0$.
\end{proof}

\subsection{A one-parameter group rotating the disc}\label{ssec:flow}

We use the argument of Ueda for rotation domains \cite[Appendix, Lemmas~A.1, A.2 and Theorem~A.3]{Ued98}, suggested by Bedford: the automorphism group is a Lie group by Kobayashi's theorem, and the closure of the group generated by the map contains a one-parameter group. We add the character $\theta$, which shows that this group rotates the disc.

\begin{proposition}\label{prop:flow}
There is a continuous one-parameter subgroup $(h_t)_{t\in\RR}$ of $\Aut(\tR)$ such that $(t,x)\mapsto h_t(x)$ is continuous on $\RR\times\tR$ and
\begin{equation}\label{eq:rotation}
h_t\bigl(\kappa_0(z)\bigr)=\kappa_0\bigl(\exp(it)z\bigr)\qquad(t\in\RR,\ z\in D^*).
\end{equation}
\end{proposition}

\begin{proof}
By Proposition~\ref{prop:retract}, the holomorphic map $\pi|_{\tR}\colon\tR\to\CC^{L+1}$ is bounded, and it is locally injective because $\pi$ is. Hence $\tR$ is hyperbolic (Lemma~\ref{lem:hyperbolic}), and $\Aut(\tR)$, with the compact-open topology, is a Lie group (Theorem~\ref{thm:kobayashi}). Let $\Gamma$ be the closure of $\{G_{\mathcal R}^m:m\in\ZZ\}$ in $\Aut(\tR)$. It is a closed abelian subgroup, hence a Lie group, and its identity component $\Gamma^\circ$ is open in $\Gamma$ and generated by $\exp(\Lie\Gamma)$ (Section~\ref{ssec:tools}). By \eqref{eq:GM}, $G_{\mathcal R}^{n_j}\to\Id$ in $\Gamma$, so there is $j_1$ with $G_{\mathcal R}^{n_j}\in\Gamma^\circ$ for $j\geq j_1$.

Let $\psi\in\Gamma$. As $\Gamma$ is metrizable, $\psi=\lim_iG_{\mathcal R}^{m_i}$ for some integers $m_i$. For $z\in D^*$, \eqref{eq:GM} gives $\psi(\kappa_0(z))=\lim_i\kappa_0(\lambda^{-m_i}z)$, so $\pi(\psi(\kappa_0(z)))=\lim_i\hvarphi(\lambda^{-m_i}z)$. As in the proof of Lemma~\ref{lem:S}(2), Lemma~\ref{lem:RI} shows that $\lambda^{-m_i}$ converges to a number $\theta(\psi)\in S^1$ which depends only on $\psi$, that $\theta\colon\Gamma\to S^1$ is continuous, and that
\begin{equation}\label{eq:character}
\psi\circ\kappa_0(z)=\kappa_0\bigl(\theta(\psi)z\bigr)\qquad(\psi\in\Gamma,\ z\in D^*).
\end{equation}
By \eqref{eq:character}, $(\psi\psi')(\kappa_0(z))=\kappa_0(\theta(\psi)\theta(\psi')z)$; applying $\pi$ and Lemma~\ref{lem:RI}, we see that $\theta$ is a homomorphism. Moreover $\theta(G_{\mathcal R})=\lambda^{-1}$.

The group $\theta(\Gamma^\circ)$ is a connected subgroup of $S^1$. It contains $\lambda^{-n_j}$ for $j\geq j_1$, and $\lambda^{-n_j}\neq1$ because $\lambda$ is not a root of unity. Hence $\theta(\Gamma^\circ)=S^1$. The continuous homomorphism $\theta$ is smooth on $\Gamma^\circ$, and its differential $d\theta\colon\Lie(\Gamma)\to\Lie(S^1)=i\RR$ is nonzero: otherwise $\theta(\exp(tv))=\exp(t\,d\theta(v))=1$ for all $v\in\Lie(\Gamma)$ \cite[Proposition~20.8(g)]{Lee13}, and $\theta$ would be trivial on $\Gamma^\circ$. Choose $v\in\Lie(\Gamma)$ with $d\theta(v)=i$, and set $h_t:=\exp(tv)$. Then $(h_t)$ is a continuous one-parameter subgroup of $\Gamma\subset\Aut(\tR)$, the map $(t,x)\mapsto h_t(x)$ is continuous, and $\theta(h_t)=\exp(it)$. We conclude by \eqref{eq:character}.
\end{proof}

The next lemma is classical (Bochner and Montgomery); for a hyperbolic manifold it is contained in \cite[Theorem~5.4.2(3)]{Kob98}. We give a short proof which uses only the continuity of $(t,x)\mapsto h_t(x)$.

\begin{lemma}\label{lem:generator}
There is a holomorphic vector field $\cV$ on $\tR$ such that, for every $x\in\tR$, the map $t\mapsto h_t(x)$ is of class $C^1$ and
\[
\frac{d}{dt}h_t(x)=\cV\bigl(h_t(x)\bigr)\qquad\text{and}\qquad D_xh_t\bigl(\cV(x)\bigr)=\cV\bigl(h_t(x)\bigr)\qquad(t\in\RR).
\]
\end{lemma}

\begin{proof}
Fix $x_1\in\tR$, a chart $u\colon U\to\CC^{\dim\tR}$ at $x_1$, and balls $b_0\Subset b_1\Subset b_2\Subset u(U)$ centred at $u(x_1)$. By continuity of $(s,x)\mapsto h_s(x)$ and compactness of $\overline b_1$, there is $\delta_0>0$ with $h_s(u^{-1}(\overline b_1))\subset u^{-1}(b_2)$ for $|s|\leq\delta_0$. The maps $\hat h_s:=u\circ h_s\circ u^{-1}\colon\overline b_1\to b_2$ are holomorphic in $w$ and jointly continuous in $(s,w)$. For $0<\delta\leq\delta_0/2$ set
\[
\mathcal I_\delta(w):=\frac1\delta\int_0^\delta\hat h_s(w)\,ds,
\]
a holomorphic map on $b_1$. As $\delta\to0$, $\mathcal I_\delta\to\Id$ uniformly on $\overline b_1$, by uniform continuity of $(s,w)\mapsto\hat h_s(w)$ and $\hat h_0=\Id$; by the Cauchy estimates, $D\mathcal I_\delta\to\mathrm I$ uniformly on $b_0$. Fix $\delta$ with $\|D\mathcal I_\delta-\mathrm I\|\leq\frac12$ on the convex set $b_0$. Then $|\mathcal I_\delta(w)-\mathcal I_\delta(w')-(w-w')|\leq\frac12|w-w'|$ on $b_0$, so $\mathcal I_\delta$ is injective on $b_0$ with invertible differential, hence a biholomorphism of $b_0$ onto an open set.

For $w$ near $u(x_1)$ and $|t|\leq\delta_0/2$ small enough that $\hat h_t(w)\in b_0$, the group law $\hat h_s(\hat h_t(w))=\hat h_{s+t}(w)$, $0\leq s\leq\delta$, gives
\[
\mathcal I_\delta\bigl(\hat h_t(w)\bigr)=\frac1\delta\int_t^{t+\delta}\hat h_s(w)\,ds,
\]
which is $C^1$ in $t$ with derivative $(\hat h_{t+\delta}(w)-\hat h_t(w))/\delta$. Hence $t\mapsto\hat h_t(w)$ is $C^1$ near $0$, and its derivative at $t=0$ is
\[
\widehat\cV(w)
 =\bigl(D\mathcal I_\delta(w)\bigr)^{-1}
   \frac{\hat h_\delta(w)-w}{\delta},
\]
which is holomorphic in $w$.

The vector $\frac{d}{dt}\big|_{t=0}h_t(x)$ does not depend on the chart, so these local expressions define a holomorphic vector field $\cV$ on $\tR$. The group law $h_{t+s}=h_s\circ h_t$ gives $\frac{d}{dt}h_t(x)=\cV(h_t(x))$ for every $t$, and differentiating $h_t\circ h_s=h_s\circ h_t$ with respect to $s$ at $s=0$ gives $D_xh_t(\cV(x))=\cV(h_t(x))$.
\end{proof}

\subsection{The maximal punctured disc}\label{ssec:maximal-disc}

\begin{figure}[ht]
\centering
\begin{tikzpicture}[x=0.9cm,y=0.78cm,>={Stealth[length=2mm]},font=\small]
  \draw[thick,rounded corners=10pt] (0,0) rectangle (11.4,6.2);
  \node at (0.55,5.75) {$\tX$};
  \fill[gray!10] (11.4,5.5) .. controls (9.0,5.45) and (7.6,5.5) .. (6.0,5.6)
        .. controls (4.0,5.75) and (0.9,5.6) .. (0.9,3.1)
        .. controls (0.9,0.7) and (3.2,0.45) .. (5.6,0.45)
        .. controls (7.6,0.45) and (9.0,0.6) .. (11.4,0.7) -- cycle;
  \draw[dashed] (11.4,5.5) .. controls (9.0,5.45) and (7.6,5.5) .. (6.0,5.6)
        .. controls (4.0,5.75) and (0.9,5.6) .. (0.9,3.1)
        .. controls (0.9,0.7) and (3.2,0.45) .. (5.6,0.45)
        .. controls (7.6,0.45) and (9.0,0.6) .. (11.4,0.7);
  \node at (1.32,3.1) {$\Omega$};
  \fill[gray!28] (1.75,3.0) .. controls (1.75,4.7) and (3.6,5.2) .. (5.0,5.1)
        .. controls (6.6,5.0) and (7.4,4.75) .. (11.4,4.6) -- (11.4,3.2)
        .. controls (8.2,3.0) and (7.4,1.3) .. (5.2,1.0)
        .. controls (3.4,0.8) and (1.75,1.3) .. (1.75,3.0);
  \draw[decorate,decoration={snake,amplitude=0.5mm,segment length=4mm}]
        (11.4,4.6) .. controls (7.4,4.75) and (6.6,5.0) .. (5.0,5.1)
        .. controls (3.6,5.2) and (1.75,4.7) .. (1.75,3.0)
        .. controls (1.75,1.3) and (3.4,0.8) .. (5.2,1.0)
        .. controls (7.4,1.3) and (8.2,3.0) .. (11.4,3.2);
  \node at (2.75,1.75) {$\tR$};
  \fill[gray!60] (4.4,3.0) ellipse [x radius=0.75,y radius=0.55];
  \foreach \rx/\ry in {0.4/0.3,0.75/0.55,1.1/0.8,1.45/1.05,1.8/1.3} {
    \draw[postaction={decorate,decoration={markings,mark=at position 0.13 with {\arrow{Stealth[length=1.6mm]}}}}]
          (4.4,3.0) ++(\rx,0) arc [start angle=0,end angle=360,x radius=\rx,y radius=\ry];
  }
  \foreach \a in {60,105,150,195,240,285,330} {
    \draw[thin] (4.4,3.0) ++({0.1*cos(\a)},{0.08*sin(\a)}) -- ++({1.7*cos(\a)},{1.22*sin(\a)});
  }
  \draw[line width=1.3pt,->] (4.52,3.02) .. controls (6.0,3.3) and (7.6,3.75) .. (11.1,3.92);
  \draw[fill=white] (4.4,3.0) circle [radius=2pt];
  \fill (5.06,3.13) circle [radius=1.5pt];
  \draw[thin] (5.06,3.13) -- (6.3,1.85) node[below right,inner sep=1pt] {$\tilde x_0$};
  \draw[thin] (3.95,3.3) -- (3.5,4.3);
  \node at (3.4,4.5) {$\kappa_0(D^*)$};
  \node at (4.95,4.6) {$h_t$};
  \node at (8.6,4.25) {$\cW=-i\cV$};
  \node[align=center] at (9.15,1.4) {leaves every\\ compact set of $\tR$};
  \draw[thin] (9.15,1.95) -- (9.25,3.8);

  \draw[->,thick] (4.4,-0.15) -- (4.4,-1.05) node[midway,right] {$\pi$};
  \draw[thin,fill=white] (4.4,-2.75) ellipse [x radius=2.2,y radius=1.45];
  \draw[pattern=north east lines,pattern color=gray!60,draw=none] (4.4,-2.75) ellipse [x radius=2.2,y radius=1.45];
  \fill[gray!60] (4.4,-2.75) ellipse [x radius=0.9,y radius=0.6];
  \fill (4.4,-2.75) circle [radius=1.6pt];
  \draw[line width=1.3pt] (4.4,-2.75) -- (6.47,-2.26);
  \foreach \a in {19,75,128,205,262,311} {
    \fill (4.4,-2.75) ++({2.2*cos(\a)},{1.45*sin(\a)}) circle [radius=1.6pt];
  }
  \draw[thin] (4.4,-2.75) -- (2.0,-2.75);
  \node[left,inner sep=2pt] at (2.0,-2.75) {$\tau(0)$};
  \node[right] at (6.6,-2.15) {$\tau^*(\zeta)\in\hB$};
  \node[right] at (7.1,-3.3) {$\tau(D(0,R))\subset\hX\cap\{\cG=0\}$};
  \draw[thin] (7.05,-3.3) -- (5.9,-3.35);
\end{tikzpicture}
\caption{The retract $\tR$, the one-parameter group $(h_t)$ and the radial field $\cW$ (Proposition~\ref{prop:flow}, Lemma~\ref{lem:maximal-disc}). Below: the maximal disc $\tau$, whose radial limits lie in $\hB$ (Lemma~\ref{lem:radial-limits}).}
\label{fig:retract}
\end{figure}
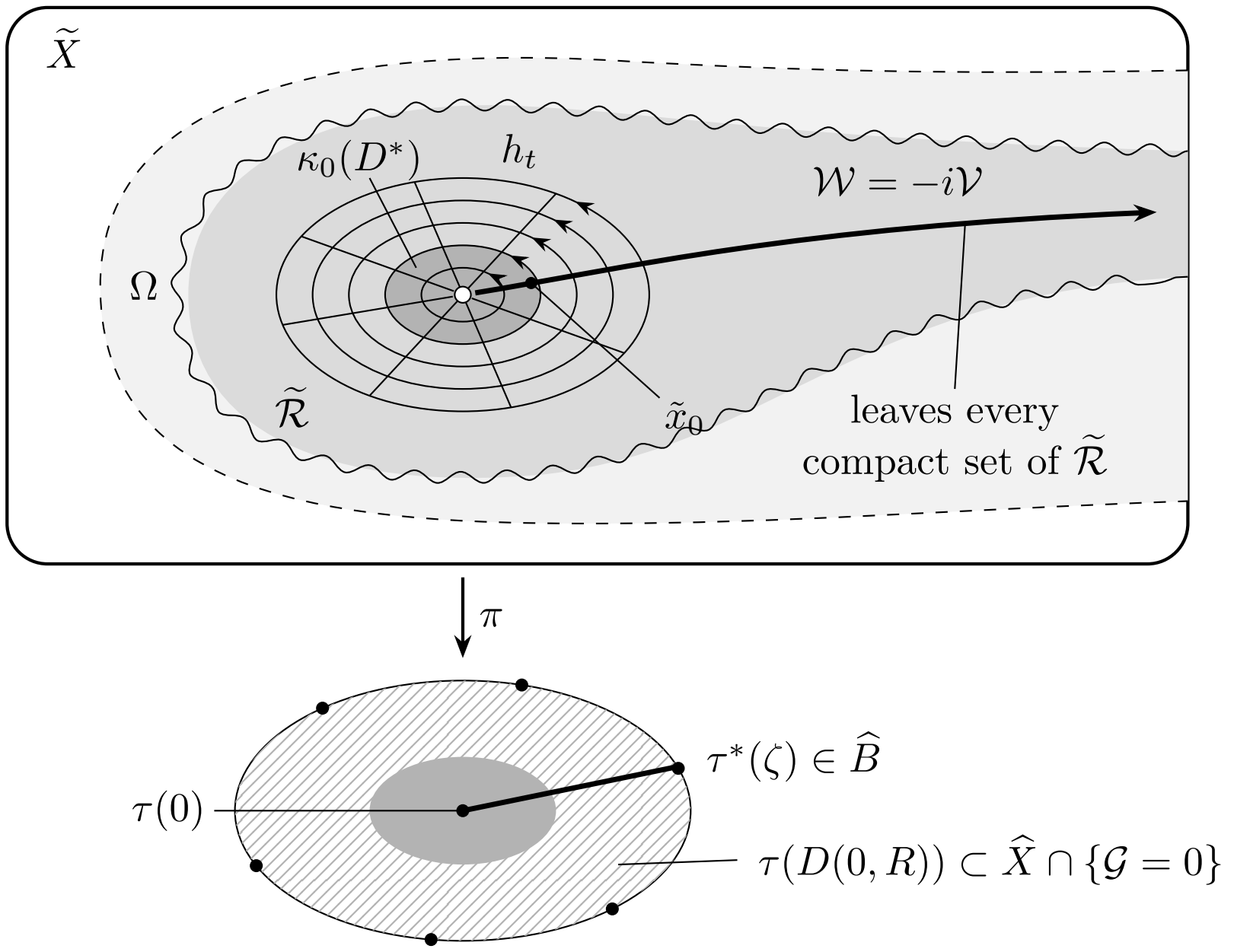

Let $\cW:=-i\cV$, a holomorphic vector field on $\tR$, and let $\xi\colon(s_-,s_+)\to\tR$ be the maximal integral curve of $\cW$, viewed as a smooth real vector field, with $\xi(0)=\tilde x_0=\kappa_0(z_0)$ (Figure~\ref{fig:retract}).

\begin{lemma}\label{lem:maximal-disc}
We have $s_-=-\infty$ and $s_+<\infty$. Set $R:=|z_0|\exp(s_+)$. Then $\varepsilon\leq R<\infty$, and the following hold.
\begin{enumerate}[label=\textup{(\arabic*)}]
\item There is a holomorphic map $\kappa\colon D^*(0,R)\to\tR$ with $\kappa=\kappa_0$ on $D^*$ and $\kappa(z_0\exp(s+it))=h_t(\xi(s))$ for $s<s_+$ and $t\in\RR$.
\item The map $\tau:=\pi\circ\kappa$ extends to a bounded holomorphic map $\tau\colon D(0,R)\to\CC^{L+1}$, with $\tau=\hvarphi$ on $D(0,\varepsilon)$.
\item For every $t\in\RR$ and every compact set $K\subset\tR$, there is $s_K<s_+$ such that $h_t(\xi(s))\notin K$ for $s_K<s<s_+$.
\end{enumerate}
\end{lemma}

\begin{proof}
As $D_xh_t$ is complex linear, Lemma~\ref{lem:generator} gives $D_xh_t(\cW(x))=\cW(h_t(x))$. Hence $h_t\circ\xi$ is an integral curve of $\cW$ through $h_t(\tilde x_0)$. It is maximal, with the same interval $(s_-,s_+)$: if it extended, applying $h_{-t}$ would extend $\xi$. Set $A(s,t):=h_t(\xi(s))$ on $(s_-,s_+)\times\RR$. Then
\[
\partial_sA=\cW(A)\qquad\text{and}\qquad\partial_tA=\cV(A)=i\,\partial_sA.
\]
Both derivatives are continuous, since $A$ is continuous by Proposition~\ref{prop:flow}. So $A$ is of class $C^1$ and holomorphic in $w=s+it$.

Set $s_\varepsilon:=\log(\varepsilon/|z_0|)$ and $\alpha(w):=\kappa_0(z_0\exp w)$ for $\operatorname{Re}w<s_\varepsilon$. By \eqref{eq:rotation}, $\alpha(s+it)=h_t(\alpha(s))$, so $\partial_t\alpha=\cV(\alpha)$. As $\alpha$ is holomorphic, $\partial_s\alpha=-i\,\partial_t\alpha=\cW(\alpha)$. Thus $s\mapsto\alpha(s)$, $s<s_\varepsilon$, is an integral curve of $\cW$ with $\alpha(0)=\tilde x_0$. By uniqueness and maximality, $s_-=-\infty$, $s_+\geq s_\varepsilon$, and $A=\alpha$ on $\{\operatorname{Re}w<s_\varepsilon\}$.

By \eqref{eq:rotation}, $h_{2\pi}(\tilde x_0)=\tilde x_0$, so $h_{2\pi}\circ\xi=\xi$ and $A(w+2\pi i)=A(w)$. With $R:=|z_0|\exp(s_+)\in[\varepsilon,\infty]$, the map $w\mapsto z_0\exp w$ is a covering of $D^*(0,R)$ by $\{\operatorname{Re}w<s_+\}$, with deck group generated by $w\mapsto w+2\pi i$. So $\kappa(z_0\exp w):=A(w)$ defines a holomorphic map $\kappa\colon D^*(0,R)\to\tR$, and $\kappa=\kappa_0$ on $D^*$. This proves~(1).

The map $\tau=\pi\circ\kappa$ is holomorphic and bounded on $D^*(0,R)$ by Proposition~\ref{prop:retract}, and $\tau=\hvarphi$ on $D^*$. Hence it extends holomorphically to $D(0,R)$, with $\tau(0)=\hvarphi(0)$. If $R=\infty$, then $\tau$ is constant by Liouville's theorem, whereas $\hvarphi$ is not. Therefore $R<\infty$, that is, $s_+<\infty$, and (2) holds.

For (3), there is $\delta>0$ such that every integral curve of $\cW$ starting at a point of the compact set $K$ exists on $(-\delta,\delta)$. If $h_t(\xi(s_1))\in K$ for some $s_1\in(s_+-\delta,s_+)$, then the maximal integral curve $h_t\circ\xi$ extends beyond $s_+$, which is impossible. So $s_K:=s_+-\delta$ works.
\end{proof}

\begin{lemma}\label{lem:radial-limits}
Let $\kappa$, $\tau$ and $R$ be as in Lemma~\ref{lem:maximal-disc}. For $\zeta\in S^1$, if the radial limit $\tau^*(\zeta):=\lim_{r\to R}\tau(r\zeta)$ exists, then $\tau^*(\zeta)\in\hB$.
\end{lemma}

\begin{proof}
We argue as in \cite[Section~3]{Le22}. Let $\zeta\in S^1$ be such that $\tau^*(\zeta)$ exists. It lies in the closed set $\hX$, because $\tau(D^*(0,R))\subset\pi(\tR)\subset\hX$. Suppose that $\tau^*(\zeta)\in\hXo$, and let $r_0:=|z_0|$. The path $\beta\colon[r_0,R]\to\hXo$ defined by $\beta(r):=\tau(r\zeta)$ for $r<R$ and $\beta(R):=\tau^*(\zeta)$ is continuous. Let $\tilde\beta$ be its lift through $\pi$ with $\tilde\beta(r_0)=\kappa(r_0\zeta)$. By uniqueness of lifts, $\tilde\beta(r)=\kappa(r\zeta)$ for $r<R$. As $\tR$ is closed in $\tX$ (Proposition~\ref{prop:retract}), the set $K:=\tilde\beta([r_0,R])\cap\tR$ is a compact subset of $\tR$, and it contains $\kappa(r\zeta)$ for $r_0\leq r<R$. Write $\zeta=\exp(it)z_0/|z_0|$ with $t\in\RR$. Then
\[
\kappa(r\zeta)=h_t\bigl(\xi(\log(r/r_0))\bigr)\in K\qquad(r_0\leq r<R),
\]
with $\log(r/r_0)\to s_+$ as $r\to R$. This contradicts Lemma~\ref{lem:maximal-disc}(3).
\end{proof}

\begin{proof}[End of the proof of Proposition~\ref{prop:key}]
By Fatou's theorem \cite[Theorem~11.32]{Rud87}, applied to the components of the bounded map $\tau$ after rescaling, the radial limit $\tau^*(\zeta)$ exists for almost every $\zeta\in S^1$, and it lies in $\hB$ by Lemma~\ref{lem:radial-limits}. Each $Q_i\circ\tau$ is a bounded holomorphic function on $D(0,R)$, whose radial limit $Q_i(\tau^*(\zeta))$ vanishes for almost every $\zeta$. By Lemma~\ref{lem:riesz}, $Q_i\circ\tau\equiv0$ for every $i$. As $\tau$ takes values in the closed set $\hX$, it follows that $\tau(D(0,R))\subset\hX\cap\{Q_1=\dots=Q_l=0\}=\hB$. But $\tau=\hvarphi$ maps $D^*$ into $\hXo$. This contradiction completes the proof of Proposition~\ref{prop:key}.
\end{proof}

\subsection{Exclusion of rotation discs}\label{ssec:proof-cor}

\begin{proof}[Proof of Corollary~\ref{cor:no-rotation-disc}]
Let $g$, $X$ and $\lambda$ be as in Corollary~\ref{cor:no-rotation-disc}, and suppose that there is a nonconstant holomorphic map $\varphi\colon D(0,\varepsilon)\to X$ with $g\circ\varphi(z)=\varphi(\lambda z)$. If $\lambda^k=1$ for some $k\geq1$, then $g^k\circ\varphi=\varphi$, so the set $\Fix(g^k)$ contains the infinite set $\varphi(D(0,\varepsilon))$, which contradicts Lemma~\ref{lem:finite-fix}.

Assume now that $\lambda$ is not a root of unity. By Lemma~\ref{lem:Z}, the Zariski closure $Z$ of $\varphi(D(0,\varepsilon))$ is irreducible, of dimension at least one, with $g(Z)=Z$. Let $\nu\colon X_Z\to Z$ be its normalization and $g_Z$ the map induced by $g$. By Proposition~\ref{prop:normalization-boundary}, $g_Z$ is a polarized endomorphism of the normal projective variety $X_Z$, and it admits a branch boundary $B_Z$. By Lemma~\ref{lem:N}, $\varphi$ lifts to a nonconstant holomorphic map $\varphi_Z\colon D(0,\varepsilon)\to X_Z$ with
\[
g_Z\circ\varphi_Z(z)=\varphi_Z(\lambda z)\qquad\text{and}\qquad\varphi_Z\bigl(D^*(0,\varepsilon)\bigr)\subset X_Z\setminus B_Z.
\]
This contradicts Proposition~\ref{prop:key}, applied to $g_Z$, $B_Z$ and $\varphi_Z$.
\end{proof}

\subsection{Closed invariant sets}\label{ssec:closed-E}

We now allow the invariant set $E$ to be closed rather than algebraic. Throughout, $g$ is a polarized endomorphism of a normal projective variety $X$ of positive dimension, with $g^*H\simeq H^{\otimes M}$, and $E\subset X$ is a closed set such that
\begin{equation}\label{eq:hypE}
g(E)\subset E\qquad\text{and}\qquad\Sing X\cup\bigcup_{n\geq1}\Br(g^n)\subset E.
\end{equation}
We use the notation of Section~\ref{ssec:cone}, and we set $\mathcal K:=\hX\cap\{\cG=0\}$, which is a compact subset of $\CC^{L+1}\setminus\{0\}$ by Lemma~\ref{lem:green}. In this subsection $\Xo$ and $\hXo$ are defined with $E$ in place of $B$:
\[
\widehat E:=\varpi^{-1}(E)\cup\{0\},\qquad\Xo:=X\setminus E,\qquad\hXo:=\hX\setminus\widehat E=\varpi^{-1}(\Xo).
\]
The set $\widehat E$ is closed in $\hX$ and invariant under $\CC^*$.

\begin{lemma}\label{lem:pc-admissible}
The set $\PC(g)$ is closed and satisfies \eqref{eq:hypE}, and it is contained in every closed set $E$ which satisfies \eqref{eq:hypE}. Every branch boundary satisfies \eqref{eq:hypE}.
\end{lemma}

\begin{proof}
By Lemma~\ref{lem:boundary}, and since the normal variety $X$ is locally irreducible, the set $\Sigma:=\Sing X\cup\bigcup_{n\geq1}\Br(g^n)$ satisfies $g(\Sigma)\subset\Sigma$. Hence $g(\overline\Sigma)\subset\overline\Sigma$ by continuity, and $\PC(g)=\overline\Sigma$. The other assertions follow from the definitions.
\end{proof}

\begin{lemma}\label{lem:lift-extension}
Let $\lambda\in S^1$ and $0<r\leq\infty$, and let $\varphi\colon D(0,r)\to X$ be a nonconstant holomorphic map with $g(\varphi(z))=\varphi(\lambda z)$.
\begin{enumerate}[label=\textup{(\arabic*)}]
\item There is a holomorphic map $\hvarphi\colon D(0,r)\to\hX\setminus\{0\}$ with $\varpi\circ\hvarphi=\varphi$ and $\hat g(\hvarphi(z))=\hvarphi(\lambda z)$.
\item Every map $\hvarphi$ as in \textup{(1)} satisfies $\cG\circ\hvarphi=0$, that is, takes values in $\mathcal K$.
\item We have $r<\infty$.
\item Every holomorphic map $\varphi_1\colon D(0,r')\to X$ which extends $\varphi$, with $r'\geq r$, satisfies $g(\varphi_1(z))=\varphi_1(\lambda z)$.
\end{enumerate}
\end{lemma}

\begin{proof}
(1) The proof of Lemma~\ref{lem:lift} uses neither $B$ nor the assumption that $\lambda$ is not a root of unity, and it applies to a disc of any radius $r\leq\infty$: holomorphic line bundles on $D(0,r)$ are trivial, and the series $\mathcal U$ converges where $A$ does.

(2) The function $u:=\cG\circ\hvarphi$ is continuous on $D(0,r)$, and $u(\lambda z)=\cG(\hat g(\hvarphi(z)))=M\,u(z)$ by Lemma~\ref{lem:green}. Let $s\in[0,r)$. The rotation $z\mapsto\lambda z$ maps the circle $\{|z|=s\}$ onto itself, so $\max_{|z|=s}u=M\max_{|z|=s}u$, and this maximum vanishes because $M\geq2$. So does the minimum. Hence $u\equiv0$.

(3) If $r=\infty$, then by (1) and (2) the map $\varphi$ lifts to an entire map with values in the compact set $\mathcal K$, which is constant by Liouville's theorem; but $\varphi$ is not constant.

(4) The holomorphic maps $g\circ\varphi_1$ and $\varphi_1(\lambda\,\cdot)$ from $D(0,r')$ to $\PP^L$ agree on $D(0,r)$.
\end{proof}

A map $\hvarphi$ as in Lemma~\ref{lem:lift-extension}(1) is called an \emph{equivariant lift} of $\varphi$.

\begin{lemma}\label{lem:coprime}
Let $\lambda\in S^1$ be a number which is not a root of unity, let $\varphi\colon D(0,r)\to X$ be a rotation disc for $\lambda$, and let $\hvarphi(z)=\sum_{m\geq0}a_mz^m$ be an equivariant lift of $\varphi$. Let $n_0\geq1$ be the greatest common divisor of the integers $m\geq1$ with $a_m\neq0$. Then the series $\hat\psi(w):=\sum_{m\geq0}a_{n_0m}w^m$ converges on $D(0,r^{n_0})$, and $\hvarphi(z)=\hat\psi(z^{n_0})$. The map $\psi:=\varpi\circ\hat\psi$ is a rotation disc for $\lambda^{n_0}$ with equivariant lift $\hat\psi$, we have $\varphi(z)=\psi(z^{n_0})$, and the exponents $m\geq1$ which occur in $\hat\psi$ have greatest common divisor $1$.
\end{lemma}

\begin{proof}
This is the reduction which precedes Lemma~\ref{lem:RI}, in which the branch boundary plays no role.
\end{proof}

From now on, $\lambda\in S^1$ is not a root of unity, $\varphi\colon D(0,\varepsilon)\to X$ is a rotation disc for $\lambda$ with $\varphi(D^*(0,\varepsilon))\cap E=\emptyset$, and $\hvarphi$ is an equivariant lift of $\varphi$. Until the proof of Theorem~A we assume that $n_0=1$ in Lemma~\ref{lem:coprime}, so that Lemma~\ref{lem:RI} holds, and we write $D^*:=D^*(0,\varepsilon)$.

As $\Sing X\subset E$, the set $\Xo$ is a complex manifold; it need not be connected. Let $U_0$ be the connected component of $\Xo$ which contains $\varphi(D^*)$. Then $\widehat U_0:=\varpi^{-1}(U_0)$ is the connected component of $\hXo$ which contains $\hvarphi(D^*)$, and it is invariant under $\CC^*$. We carry out the construction on this connected component, which is a second-countable complex manifold.

\begin{proposition}\label{prop:construction-E}
Fix $z_0\in D^*$, and let $x_0:=\hvarphi(z_0)$ and the loop $\ell$ be as in Section~\ref{ssec:covering}. Let $\pi\colon(\tX,\tilde x_0)\to(\widehat U_0,x_0)$ be the connected covering with $\pi_*\pi_1(\tX,\tilde x_0)=\langle\gamma_0\rangle$, where $\gamma_0\in\pi_1(\widehat U_0,x_0)$ is the class of the loop $\hvarphi\circ\ell$. With the notation of Sections~\ref{ssec:covering}--\ref{ssec:maximal-disc}, set $\Omega_h:=\tX\setminus h^{-1}(\widehat E)$ for $h\in S$, and let $\Omega_h^0$ be the connected component of $\Omega_h$ which contains $\kappa_0(D^*)$. Set $\Omega:=\Omega_\rho^0$. Then the statements of Lemmas~\ref{lem:G}, \ref{lem:S}\textup{(1)--(2)}, \ref{lem:L}, \ref{lem:W}, \ref{lem:idempotent}, of Propositions~\ref{prop:retract}, \ref{prop:automorphism}, \ref{prop:flow} and of Lemmas~\ref{lem:generator}, \ref{lem:maximal-disc} hold with $\widehat U_0$, $\widehat E$ and $\Omega_h^0$ in place of $\hXo$, $\hB$ and $\Omega_h$. Moreover, for every $h\in S$ the set $\Omega_h$ is open and contains $\kappa_0(D^*)$, and $h(\Omega_h^0)\subset\widehat U_0$.
\end{proposition}

In particular $\tR$ is a closed connected submanifold of $\tX$ with $\pi(\tR)\subset\widehat U_0\cap\mathcal K$. We write $R_1$ for the radius which is denoted by $R$ in Lemma~\ref{lem:maximal-disc}; thus $\varepsilon\leq R_1<\infty$, and $\tau=\pi\circ\kappa$ on $D^*(0,R_1)$.

\begin{proof}
The modifications concern the covering, the domains of the lifted limits, and their connected components.

\emph{Lemma~\ref{lem:covering}.} By \eqref{eq:hypE}, $\hat g(\widehat E)\subset\widehat E$ and $\hat g^{-1}(\hXo)\subset\hXo$. The proof of Lemma~\ref{lem:covering}(2) uses $B$ only through the inclusion $\Br(g)\subset B$ and, for surjectivity, through the connectedness of $\Xo$. It shows that $\hat g\colon\hat g^{-1}(\hXo)\to\hXo$ is a proper local homeomorphism. This map is surjective, because $g$ is surjective and $\hat g(tv)=t^M\hat g(v)$. Hence it is a finite covering, and so is its restriction to a connected component $Y_0$ of $\hat g^{-1}(\hXo)$, onto the connected component of $\hXo$ which contains $\hat g(Y_0)$.

\emph{Lemma~\ref{lem:G}.} The component $Y_0$ of $\hat g^{-1}(\hXo)$ which contains $y_0$ is a connected subset of $\hXo$ which meets $\widehat U_0$ at $y_0\in\hvarphi(D^*)$, so $Y_0\subset\widehat U_0$; and $\hat g(Y_0)$ contains $x_0\in\widehat U_0$. So $\eta=\hat g|_{Y_0}\colon Y_0\to\widehat U_0$ is a covering. The map $\sigma$ takes values in $Y_0\subset\widehat U_0$, the base of $\pi$, so it can be lifted through $\pi$. The rest of the proof is unchanged.

\emph{Lemma~\ref{lem:S}.} Parts (1) and (2) do not involve $B$. Let $h\in S$. The set $\Omega_h$ is open because $\widehat E$ is closed, and it contains $\kappa_0(D^*)$ by \eqref{eq:theta}. The set $h(\Omega_h^0)$ is a connected subset of $\hXo$ which meets $\widehat U_0$, so it lies in $\widehat U_0$.

\emph{Lemma~\ref{lem:L}.} The proof uses only that $\Omega_h$ is connected and that $h$ maps it into the base of $\pi$. It applies to $\Omega_h^0$.

\emph{Lemma~\ref{lem:W}.} The proof shows that $k_{a_j+a'_j}\to h\circ\tilde h'$ locally uniformly on $\Omega_{h'}^0$. Every limit point of $(k_{a_j+a'_j})$ in $S$ is a holomorphic map $\tX\to\CC^{L+1}$ which equals $h\circ\tilde h'$ on the nonempty open set $\Omega_{h'}^0$. By the identity theorem on the connected manifold $\tX$, there is exactly one limit point. This replaces the density of $\Omega_{h'}$. The rest of the proof, and the proof of Lemma~\ref{lem:idempotent}, are unchanged.

\emph{Proposition~\ref{prop:retract}.} Let $x\in\Omega$. As $\rho(\tilde\rho(x))=\rho(x)\notin\widehat E$, the set $\tilde\rho(\Omega)$ is a connected subset of $\Omega_\rho$. It contains $\tilde\rho(\kappa_0(D^*))=\kappa_0(D^*)$, so $\tilde\rho(\Omega)\subset\Omega$. For the closedness of $\tR$, we have $\cZ\subset\Omega_\rho$ because $\pi$ takes values in $\widehat U_0$. Let $x$ be a point of the closure of $\tR$ in $\tX$. Then $x\in\cZ\subset\Omega_\rho$, and the connected component of $\Omega_\rho$ which contains $x$ is open, so it meets $\tR\subset\Omega$ and is equal to $\Omega$. As $\tR$ is closed in $\Omega$, we get $x\in\tR$.

\emph{Proposition~\ref{prop:automorphism}.} For $x\in\Omega$ the proof gives $\rho(G(x))\notin\widehat E$, that is, $G(x)\in\Omega_\rho$. The connected set $G(\Omega)$ contains $G(\kappa_0(D^*))=\kappa_0(D^*)$, hence $G(\Omega)\subset\Omega$.

The proofs of Proposition~\ref{prop:flow} and of Lemmas~\ref{lem:generator} and~\ref{lem:maximal-disc} involve only the manifold $\tR$, the bounded locally injective map $\pi|_{\tR}$, the automorphism $G_{\mathcal R}$, the map $\kappa_0$ and Lemma~\ref{lem:RI}.
\end{proof}

\begin{proposition}\label{prop:radial-E}\label{prop:R-maximal}
\begin{enumerate}[label=\textup{(\arabic*)}]
\item The map $\tau$ takes values in $\hX\setminus\{0\}$, and $\tau(D^*(0,R_1))\subset\widehat U_0\cap\mathcal K$.
\item Let $\zeta\in S^1$. If the radial limit $\tau^*(\zeta):=\lim_{r\to R_1}\tau(r\zeta)$ exists, then $\tau^*(\zeta)\in\widehat E$.
\item If $\varphi$ has no holomorphic extension $\varphi_1\colon D(0,r)\to X$ with $r>\varepsilon$ and $\varphi_1(D^*(0,r))\cap E=\emptyset$, then $R_1=\varepsilon$ and $\tau=\hvarphi$.
\end{enumerate}
\end{proposition}

\begin{proof}
(1) On $D^*(0,R_1)$ we have $\tau=\pi\circ\kappa$ with $\kappa$ valued in $\tR$, and $\tau(0)=\hvarphi(0)\neq0$.

(2) By (1), $\tau^*(\zeta)\in\mathcal K$. Suppose that $\tau^*(\zeta)\notin\widehat E$. The connected component of $\hXo$ which contains $\tau^*(\zeta)$ is open in $\hX$, so it contains $\tau(r\zeta)\in\widehat U_0$ for $r$ close to $R_1$. Hence $\tau^*(\zeta)\in\widehat U_0$, and the path $\beta$ of the proof of Lemma~\ref{lem:radial-limits} lies in $\widehat U_0$, the base of $\pi$. That proof now gives a contradiction.

(3) By (1), the map $\varpi\circ\tau\colon D(0,R_1)\to X$ is holomorphic and maps $D^*(0,R_1)$ into $U_0\subset X\setminus E$. It extends $\varphi$, because $\tau=\hvarphi$ on $D(0,\varepsilon)$. Hence $R_1\leq\varepsilon$.
\end{proof}

\begin{proof}[Proof of Theorem~A]
Let $\lambda$, $\varphi$ and $R$ be as in Theorem~A. By the definition of $R$ and the identity theorem, $\varphi$ has a unique holomorphic extension $\varphi\colon D(0,R)\to X$, and $\varphi(D^*(0,R))\cap E=\emptyset$. By Lemma~\ref{lem:lift-extension}(4) and~(3), this map is a rotation disc for $\lambda$, and $R<\infty$. This proves (1). From now on we assume that $\varepsilon=R$.

Let $\hvarphi\colon D(0,R)\to\mathcal K$ be an equivariant lift of $\varphi$ (Lemma~\ref{lem:lift-extension}), and let $n_0$, $\hat\psi$ and $\psi$ be as in Lemma~\ref{lem:coprime}. Then $\psi$ is a rotation disc for $\lambda^{n_0}$ on $D(0,R^{n_0})$, and $\psi(D^*(0,R^{n_0}))=\varphi(D^*(0,R))$ is disjoint from $E$. If $\psi_1$ is an extension of $\psi$ to a larger disc whose punctured image is disjoint from $E$, then $z\mapsto\psi_1(z^{n_0})$ is such an extension of $\varphi$; so there is none. Propositions~\ref{prop:construction-E} and~\ref{prop:radial-E} apply to $\psi$, $\hat\psi$ and $\lambda^{n_0}$, with $\varepsilon$ replaced by $R^{n_0}$. By Proposition~\ref{prop:radial-E}(3), the map of Lemma~\ref{lem:maximal-disc} for $\psi$ is $\hat\psi$, so by Proposition~\ref{prop:radial-E}(2) every radial limit of $\hat\psi$ which exists lies in $\widehat E$. As $\hvarphi(r\zeta)=\hat\psi(r^{n_0}\zeta^{n_0})$, the same holds for $\hvarphi$.

The map $\hvarphi$ is bounded. By Fatou's theorem \cite[Theorem~11.32]{Rud87}, $\hvarphi^*(\zeta)$ exists for almost every $\zeta\in S^1$. For such $\zeta$ we have $\hvarphi^*(\zeta)\in\widehat E\cap\mathcal K$, so that $\hvarphi^*(\zeta)\neq0$ and $\varphi(r\zeta)=\varpi(\hvarphi(r\zeta))$ converges to $\varpi(\hvarphi^*(\zeta))\in E$ as $r\to R$. If $\varphi^*(\zeta)$ exists, then $\varphi(r\lambda\zeta)=g(\varphi(r\zeta))$ converges to $g(\varphi^*(\zeta))$. This proves (2).

Let $A\subset X$ be a closed algebraic subset, and let $Q_1,\dots,Q_l$ be homogeneous polynomials such that
\[
\varpi^{-1}(A)=(\hX\setminus\{0\})\cap\{Q_1=\dots=Q_l=0\}.
\]
Let $N_A$ be the set of $\zeta\in S^1$ such that $\varphi^*(\zeta)$ exists and lies in $A$. If $\zeta\in N_A$ and $\hvarphi^*(\zeta)$ exists, then $\varpi(\hvarphi^*(\zeta))=\varphi^*(\zeta)\in A$, so the bounded holomorphic function $Q_i\circ\hvarphi$ has radial limit $0$ at $\zeta$. The set of points at which $Q_i\circ\hvarphi$ has radial limit $0$ is a Borel set. After removing the null set where $\hvarphi^*$ does not exist, we see that if $N_A$ is not of measure zero, this Borel set has positive measure for every $i$, and Lemma~\ref{lem:riesz} gives $Q_i\circ\hvarphi\equiv0$. Then $\hvarphi$ takes values in $\varpi^{-1}(A)$, and $A$ contains the image of $\varphi$. This proves (3).
\end{proof}

\begin{corollary}\label{cor:countable-algebraic}
Let $\lambda\in S^1$ be a number which is not a root of unity, and let $\varphi\colon D(0,\varepsilon)\to X$ be a rotation disc for $\lambda$ with $\varphi(D^*(0,\varepsilon))\cap E=\emptyset$. Then $E$ is not contained in a countable union of closed algebraic subsets of $X$ none of which contains $\varphi(D(0,\varepsilon))$. In particular, $E$ is not a countable union of closed algebraic subsets of $X$.
\end{corollary}

\begin{proof}
Suppose that $E\subset\bigcup_nA_n$, where the $A_n$ are closed algebraic subsets which do not contain $\varphi(D(0,\varepsilon))$. By item~(2) of Theorem~A, $\varphi^*(\zeta)$ exists and lies in $\bigcup_nA_n$ for almost every $\zeta\in S^1$. By item~(3), for each $n$ the set of $\zeta$ with $\varphi^*(\zeta)\in A_n$ has measure zero. This is a contradiction. A subset of $E$ does not contain $\varphi(D(0,\varepsilon))$, which gives the last assertion.
\end{proof}

When $E=B$ is a branch boundary, Corollary~\ref{cor:countable-algebraic}, with $A_1=B$, is Proposition~\ref{prop:key}. The assumption that $\lambda$ is not a root of unity excludes no rotation disc, by Lemma~\ref{lem:finite-fix}, as in the proof of Corollary~\ref{cor:no-rotation-disc}.

\section{Bounded inverse branches}\label{sec:bib}

Throughout this section, $Y$ is a normal complex space, $h\colon Y\to Y$ is a holomorphic map, and $q\in Y$ is a fixed point of $h$ at which $h$ is a local biholomorphism. We write $h_q^*$ for the automorphism of $V_q$ induced by $h$ (Definition~\ref{def:multiplier}); its transpose is $D_qh$. A family of holomorphic maps into $Y$ is \emph{normal} if each of its sequences has a subsequence converging locally uniformly to a holomorphic map into $Y$.

\begin{definition}\label{def:bib}
The fixed point $q$ has \emph{bounded inverse branches} if there are a connected open neighbourhood $U$ of $q$ and holomorphic maps $b_n\colon U\to Y$, $n\geq0$, with $h^n\circ b_n=\Id_U$ and $b_n(q)=q$, such that the following hold.
\begin{enumerate}[label=\textup{(\arabic*)}]
\item For every $k\geq0$ and every connected open neighbourhood $U_1\subset U$ of $q$ with $b_k(U_1)\subset U$,
\begin{equation}\label{eq:composition}
b_n\circ b_k=b_{n+k}\quad\text{on }U_1,\ \text{for all }n\geq0.
\end{equation}
\item The family $(b_n)_{n\geq0}$ is normal.
\item For every neighbourhood $W$ of $q$, there is a neighbourhood $U_0\subset U$ of $q$ with $b_n(U_0)\subset W$ for all $n$.
\end{enumerate}
\end{definition}

Since $h^n$ is a local biholomorphism at $q$, the germ of $b_n$ at $q$ is the inverse of the germ of $h^n$, so $b_n$ induces $(h_q^*)^{-n}$ on $V_q$.

Assume that $q$ has bounded inverse branches. Fix an open neighbourhood of $q$ in $U$ identified with a closed analytic subset of an open subset of $\CC^\varkappa$, with coordinates $z_1,\dots,z_\varkappa$ and $q=0$. For $r>0$ small, let $\beta(r)$ be the connected component of $Y\cap\BB(r)$ containing $q$, where $\BB(r)\subset\CC^\varkappa$ is the ball of radius $r$ about $0$. It is an open neighbourhood of $q$ in $Y$ and a closed analytic subset of $\BB(r)$, and its closure in $Y$ is compact and contained in this neighbourhood of $q$. These sets form a fundamental system of neighbourhoods of $q$, and $\overline{\beta(r')}\subset\beta(r)$ for $r'<r$, since $\overline{\beta(r')}$ is connected, contains $q$ and lies in $Y\cap\BB(r)$. By Definition~\ref{def:bib}(3), there are radii $r_0>r_1>r_2>0$ such that, for every $n\geq0$,
\begin{equation}\label{eq:radii}
b_n\bigl(\beta(r_1)\bigr)\subset\beta(r_0)\qquad\text{and}\qquad b_n\bigl(\beta(r_2)\bigr)\subset\beta(r_1/2).
\end{equation}
As $\beta(r_1)$ is connected and open and $b_k(\beta(r_1))\subset U$, \eqref{eq:composition} gives $b_n\circ b_k=b_{n+k}$ on $\beta(r_1)$ for all $n,k\geq0$.

At a regular point, part~(1) of the next proposition is the Cauchy estimate in the proof of \cite[Proposition~3.3.1]{SZ26b}; part~(2) uses the extraction in the proof of \cite[Theorem~2.1.29]{Aba89}, as does Lemma~\ref{lem:idempotent}.

\begin{proposition}\label{prop:local-retract}
Assume that $q$ has bounded inverse branches, and let $r_0>r_1>r_2$ be radii as in \eqref{eq:radii}. Let $E^c$ \textup(resp.\ $E^u$\textup) be the sum of the generalized eigenspaces of $D_qh$ on $T_qY$ for the eigenvalues of modulus one \textup(resp.\ greater than one\textup), and set $c:=\dim E^c$. Denote by $(E^c)^*$ and $(E^u)^*$ the annihilators of $E^u$ and $E^c$ in $V_q$.
\begin{enumerate}[label=\textup{(\arabic*)}]
\item The powers $(h_q^*)^{-n}$, $n\geq0$, are bounded. Hence every eigenvalue of $D_qh$ has modulus at least one, $T_qY=E^c\oplus E^u$ and $V_q=(E^c)^*\oplus(E^u)^*$, and $D_qh$ and $h_q^*$ are diagonalizable on $E^c$ and on $(E^c)^*$, respectively.
\item There are integers $n_j\to\infty$ and a holomorphic map $\rho\colon U\to Y$ such that $b_{n_j}\to\rho$ locally uniformly on $U$, $\rho(q)=q$,
\[
\rho\circ\rho=\rho\quad\text{on }\beta(r_2)\qquad\text{and}\qquad b_k\circ\rho=\rho\circ b_k\quad\text{on }\beta(r_1),\ \text{for every }k\geq0,
\]
and $D_q\rho$ is the projection of $T_qY$ onto $E^c$ along $E^u$. Equivalently, $\rho^*$ is the projection of $V_q$ onto $(E^c)^*$ along $(E^u)^*$.
\end{enumerate}
\end{proposition}

\begin{proof}
(1) The functions $z_k\circ b_n$ are defined on $\beta(r_1)$, vanish at $q$ and are bounded by $r_0$ there, by \eqref{eq:radii}, so their classes $(h_q^*)^{-n}[z_k]$ in $V_q$ are bounded, by Lemma~\ref{lem:jet} applied on $\beta(r_1)$. As the classes $[z_k]$ span $V_q$, the powers $(h_q^*)^{-n}$ are bounded. A linear map with bounded powers has its eigenvalues in the closed unit disc, and its eigenvalues of modulus one are semisimple. Applied to $(h_q^*)^{-1}$, whose transpose is $(D_qh)^{-1}$, this shows that every eigenvalue of $D_qh$ has modulus at least one, so that $T_qY=E^c\oplus E^u$, and that $D_qh$ is diagonalizable on $E^c$. The annihilators of complementary subspaces are complementary, so $V_q=(E^c)^*\oplus(E^u)^*$. Finally, $h_q^*$ preserves $(E^c)^*$, as $D_qh$ preserves $E^u$, and restriction to $E^c$ identifies $(E^c)^*$ with the dual of $E^c$ and $h_q^*|_{(E^c)^*}$ with the transpose of $D_qh|_{E^c}$, which is diagonalizable.

(2) By Definition~\ref{def:bib}(2), there are integers $m_1<m_2<\cdots$ with $m_{j+1}\geq2m_j+j$ such that $b_{m_j}$ converges locally uniformly on $U$ to a holomorphic map $u$. Set $s_j:=m_{j+1}-m_j$ and $t_j:=s_j-m_j$, so that $s_j\geq t_j\geq j$. Choose an infinite set $J$ of indices along which $b_{s_j}\to\rho$ and $b_{t_j}\to u'$ locally uniformly on $U$, for holomorphic maps $\rho,u'\colon U\to Y$. Then $\rho(q)=q$.

Let $(a_j)$, $(c_j)$ be sequences of nonnegative integers such that $b_{a_j}\to v$ and $b_{c_j}\to w$ locally uniformly on $U$. On $\beta(r_1)$ we have $b_{a_j+c_j}=b_{a_j}\circ b_{c_j}$, and this converges locally uniformly on $\beta(r_1)$ to $v\circ w$, since $b_{c_j}\to w$ on $\beta(r_1)$ with values in $\overline{\beta(r_0)}$, $b_{a_j}\to v$ uniformly on this compact subset of $U$, and $v$ is continuous.

We apply this along $J$ to the pairs of sequences $(m_j,s_j)$ and $(t_j,m_j)$, and, for a fixed $k\geq0$, to $(k,s_j)$ and $(s_j,k)$. As $m_j+s_j=m_{j+1}$ and $t_j+m_j=s_j$, we obtain, on $\beta(r_1)$,
\[
u=u\circ\rho,\qquad \rho=u'\circ u,\qquad
b_k\circ\rho=\lim_{j\in J}b_{k+s_j}=\rho\circ b_k.
\]
For $x\in\beta(r_2)$, the point $\rho(x)$ lies in $\overline{\beta(r_1/2)}\subset\beta(r_1)$, by \eqref{eq:radii}, so that
\[
\rho(\rho(x))=u'(u(\rho(x)))=u'(u(x))=\rho(x).
\]
Let $(n_j)$ be the sequence $(s_j)_{j\in J}$.

The functions $z_k\circ b_{n_j}$ vanish at $q$ and converge locally uniformly on $\beta(r_1)$ to $z_k\circ\rho$. As the classes $[z_k]$ span $V_q$, Lemma~\ref{lem:jet} gives $\rho^*=\lim_j(h_q^*)^{-n_j}$; transposing, $D_q\rho=\lim_j(D_qh)^{-n_j}$. These maps preserve $E^c$ and $E^u$ and tend to $0$ on $E^u$; on $E^c$ they are diagonalizable with eigenvalues of modulus one, so their limit is invertible there. Being idempotent, $D_q\rho$ is the projection onto $E^c$ along $E^u$, and its transpose $\rho^*$ is the projection of $V_q$ onto $(E^c)^*$ along $(E^u)^*$.
\end{proof}

\section{First route: the Main Theorem}\label{sec:route-one}

Throughout this section, $f$ is a PCF endomorphism of $\PP^N$ of degree $d\geq2$, and $p$ is a periodic point of $f$ of period $a$. Let $E^\times\subset T_p\PP^N$ be the sum of the generalized eigenspaces of $D_pf^a$ for its nonzero eigenvalues. For $n\geq1$, the iterate $f^n$ is PCF of degree $d^n$, as $\PC(f^n)=\PC(f)$: on $\PP^N$ the sets $\Br(f^j)=\bigcup_{k\leq j}f^k(C(f))$ increase with $j$ (Section~\ref{ssec:branch-boundary}).

\begin{theorem}[Sun--Zhao {\cite[Theorem~3.2.1 and Lemma~3.2.2]{SZ26b}}]\label{thm:core}
Assume that $\dim E^\times>0$. There are an integer $b\geq1$ and an irreducible algebraic subvariety $W\subset\PP^N$ of dimension $\dim E^\times$ with $f^{ab}(W)=W$, such that the germ $(W,p)$ has a smooth irreducible component $\mathcal B_W$ with $T_p\mathcal B_W=E^\times$, and $f^{ab}$ induces an automorphism of the germ $(\mathcal B_W,p)$. On the normalization of $W$, the point corresponding to $\mathcal B_W$ is a smooth point, fixed by the map induced by $f^{ab}$, at which the normalization map has injective differential.
\end{theorem}

Set $F:=f^{ab}$; the variety $W$ is called the \emph{algebraic core} of $p$. Let $\nu\colon X\to W$ be the normalization, $g\colon X\to X$ the map induced by $F$ (Section~\ref{ssec:invariant}), and $q\in\Reg X$ the fixed point of $g$ corresponding to $\mathcal B_W$, so that $\nu(q)=p$; here the points of $\nu^{-1}(p)$ correspond to the irreducible components of the germ $(W,p)$, the germ of $X$ at each point being mapped onto the corresponding component \cite[Chapter~8, \S\S4.1--4.2]{GrRe84}. We call $(X,g,q)$ the \emph{normalized core} of $p$.

Since $X$ is smooth at $q$, $\dim X=\dim\mathcal B_W$ and $D_q\nu$ is injective, $\nu$ maps the germ $(X,q)$ biholomorphically onto its image, which is the component $\mathcal B_W$ corresponding to $q$; so $D_q\nu$ is an isomorphism of $T_qX$ onto $E^\times$. Differentiating $\nu\circ g=F\circ\nu$ at $q$ gives $D_pF\circ D_q\nu=D_q\nu\circ D_qg$. Since $D_pF=(D_pf^a)^b$ preserves $E^\times$, we obtain
\begin{equation}\label{eq:core-spectrum}
D_qg\ \text{is conjugate to}\ \bigl(D_pf^a|_{E^\times}\bigr)^b.
\end{equation}
This map is invertible, so $g$ is a local biholomorphism at $q$.

\begin{proposition}[Sun--Zhao {\cite[Lemmas~2.3.2, 2.3.3 and~A.0.2]{SZ26b}}]\label{prop:sz-inverse}
Let $j\geq1$, let $Z\subset\PP^N$ be an irreducible algebraic subvariety of positive dimension with $f^j(Z)=Z$, let $\nu\colon X\to Z$ be its normalization and $g\colon X\to X$ the map induced by $f^j$. Then $H:=\nu^*\cO(1)$ is ample, $g^*H\simeq H^{\otimes d^j}$, and $g$ admits a branch boundary. Let $q\in\Reg X$ be a fixed point of $g$ at which $g$ is a local biholomorphism. There are a connected neighbourhood $U\subset\Reg X$ of $q$ and holomorphic maps $b_n\colon U\to X$, $n\geq0$, such that $g^n\circ b_n=\Id_U$, $b_n(q)=q$ and $g\circ b_{n+1}=b_n$, the family $(b_n)_{n\geq0}$ is normal, and for every neighbourhood $W'$ of $q$ there is a neighbourhood $U_0\subset U$ of $q$ with $b_n(U_0)\subset W'$ for all $n$.
\end{proposition}

In particular, $g$ is a polarized endomorphism of the normal projective variety $X$, and $q$ has bounded inverse branches (Definition~\ref{def:bib}). Indeed, replacing $U$ by the connected component of its interior containing $q$, which preserves the other properties, we may assume that $U$ is open. If $U_1\subset U$ is a connected open neighbourhood of $q$ with $b_k(U_1)\subset U$, then $b_n\circ b_k$ and $b_{n+k}$ are holomorphic right inverses of $g^{n+k}$ on $U_1$, as $g^{n+k}\circ b_n\circ b_k=g^k\circ b_k=\Id_{U_1}$, and both fix $q$. They agree near $q$, where $g^{n+k}$ is injective, hence on the connected manifold $U_1$ by the identity theorem. This gives \eqref{eq:composition}.

The germ $\cR$ in the following proposition was constructed by Astorg \cite[Lemma~8]{Ast20}, by an idempotent extraction as in the proof of \cite[Theorem~2.1.29]{Aba89} and a linearization lemma of Forn\ae ss and Sibony \cite[Lemma~6.8]{FS94}. We linearize by Ces\`aro averages, as in \cite[Section~3]{Le22}.

\begin{proposition}\label{prop:centre-germ}
In the setting of Section~\textup{\ref{sec:bib}}, assume that $q\in\Reg Y$ has bounded inverse branches, choose the embedding of Section~\textup{\ref{sec:bib}} to be a chart which identifies a neighbourhood of $q$ in $U$ with a neighbourhood of $0$ in $T_qY$, with differential the identity at $q$, and let $E^c$, $c$ and $\rho$ be as in Proposition~\textup{\ref{prop:local-retract}}. Near $q$, the fixed-point set of $\rho$ is a complex submanifold $\cR$ of dimension $c$ with $T_q\cR=E^c$. Moreover $h(\cR)=\cR$ and $b_1(\cR)\subset\cR$ as germs at $q$, and the germ of $h|_{\cR}$ at $q$ is holomorphically conjugate to $D_qh|_{E^c}$.
\end{proposition}

\begin{proof}
Set $A:=D_qh$, and give $T_qY$ the Euclidean norm of the chart, so that $\beta(r)$ is the ball of radius $r$ about $q$. Let $P$ be the projection of $T_qY$ onto $E^c$ along $E^u$. It commutes with $A$, and as $A|_{E^c}$ is diagonalizable with eigenvalues of modulus one (Proposition~\ref{prop:local-retract}(1)), there is $C>0$ with $\|A^kP\|\leq C$ for all $k\in\ZZ$.

Let $\Theta_1$ be the connected component containing $q$ of $\{x\in\beta(r_2):\rho(x)\in\beta(r_2)\}$. The set $\rho(\Theta_1)$ is connected, contains $q$, and lies in $\{x\in\beta(r_2):\rho(x)\in\beta(r_2)\}$ because $\rho\circ\rho=\rho$ on $\beta(r_2)$. So $\rho(\Theta_1)\subset\Theta_1$, and $\rho|_{\Theta_1}$ is a holomorphic retraction of $\Theta_1$. By Lemma~\ref{lem:rossi}, its fixed-point set $\cR:=\rho(\Theta_1)$ is a closed connected submanifold of $\Theta_1$. As $\rho$ maps $\Theta_1$ into $\cR$ and is the identity on $\cR$, we get $T_q\cR=\operatorname{Im}D_q\rho=E^c$ by Proposition~\ref{prop:local-retract}(2), and $\dim\cR=c$.

By Proposition~\ref{prop:local-retract}(2), $b_k\circ\rho=\rho\circ b_k$ on $\beta(r_1)$ for every $k\geq0$. By Definition~\ref{def:bib}(3), there is $r_3\in(0,r_2]$ with $b_k(\beta(r_3))\subset\Theta_1$ for all $k$. For $x\in\cR\cap\beta(r_3)$, the point $b_k(x)\in\Theta_1$ satisfies $\rho(b_k(x))=b_k(\rho(x))=b_k(x)$, so $b_k(x)\in\cR$. Since $b_1$ is a local inverse of $h$ at $q$, we also get $h\circ\rho=\rho\circ h$ near $q$, hence $h(\cR)\subset\cR$ near $q$. As $D_q(h|_{\cR})=A|_{E^c}$ is invertible, $h$ maps a neighbourhood of $q$ in $\cR$ onto a neighbourhood of $q$ in $\cR$, so that $h(\cR)=\cR$ as germs.

It remains to linearize $h|_{\cR}$. On $\beta(r_1)$ set
\[
\mathcal L_n:=\frac1n\sum_{k=0}^{n-1}A^kP\circ b_k
\qquad(n\geq1).
\]
These holomorphic maps take values in $E^c$ and are bounded by $Cr_0$, by \eqref{eq:radii}. By Montel's theorem, a subsequence converges locally uniformly on $\beta(r_1)$ to a holomorphic map $\mathcal L\colon\beta(r_1)\to E^c$.

For $x\in\beta(r_2)$ we have $b_1(x)\in\beta(r_1)$, and $b_k\circ b_1=b_{k+1}$ gives
\[
\mathcal L_n(b_1(x))
 =A^{-1}\left(\mathcal L_n(x)
   +\frac{A^nP\,b_n(x)-Px}{n}\right).
\]
The second term in the parenthesis is $O(1/n)$, uniformly on $\beta(r_2)$, so $\mathcal L\circ b_1=A^{-1}\circ\mathcal L$ on $\beta(r_2)$. Moreover
\[
D_q\mathcal L_n=\frac1n\sum_{k=0}^{n-1}A^kPA^{-k}=P
\]
for every $n$, so $D_q\mathcal L=P$, and $\mathcal L|_{\cR}$ is a local biholomorphism of $(\cR,q)$ onto $(E^c,0)$. As $b_1(\cR)\subset\cR$ near $q$, the chart $\mathcal L|_{\cR}$ conjugates the germ of $b_1|_{\cR}$ to $A^{-1}|_{E^c}$, hence the germ of $h|_{\cR}$, which is its inverse, to $A|_{E^c}$.
\end{proof}

\begin{proof}[Proof of the Main Theorem]
If $f^{a'}(p)=p$, then $a'=ka$ for some $k\geq1$, and the eigenvalues of $D_pf^{a'}=(D_pf^a)^k$ are the $k$-th powers of those of $D_pf^a$. As $\lambda^k=0$ if and only if $\lambda=0$, and $|\lambda^k|>1$ if and only if $|\lambda|>1$, it suffices to treat $D_pf^a$. If $E^\times=\{0\}$, all eigenvalues of $D_pf^a$ vanish. Otherwise, let $b$, $W$, $F=f^{ab}$ and the normalized core $(X,g,q)$ be as above. By Proposition~\ref{prop:sz-inverse} for $j=ab$ and $Z=W$, and the paragraph after it, $q$ has bounded inverse branches, as $g$ is a local biholomorphism at $q$. By Proposition~\ref{prop:local-retract}(1), applied to $Y:=X$ and $h:=g$, every eigenvalue of $D_qg$ has modulus at least one.

Suppose that $D_qg$ has an eigenvalue $\lambda$ of modulus one, so that $c\geq1$. By Propositions~\ref{prop:local-retract}(1) and~\ref{prop:centre-germ}, there are a germ $\cR\ni q$ of complex manifold with $g(\cR)=\cR$ and a chart $\Psi\colon(\cR,q)\to(\CC^c,0)$ with $\Psi\circ g=\Lambda\circ\Psi$, where $\Lambda=\operatorname{diag}(\lambda,\mu_2,\dots,\mu_c)$ and $|\mu_i|=1$. As $\Lambda$ preserves polydiscs centred at $0$, for $\varepsilon$ small the map $\varphi_0(z):=\Psi^{-1}(z,0,\dots,0)$ is defined on $D(0,\varepsilon)$. It is injective and satisfies $g\circ\varphi_0(z)=\varphi_0(\lambda z)$. The map $\psi:=\nu\circ\varphi_0$ is nonconstant, because $\nu$ has finite fibres, and $F\circ\psi(z)=\psi(\lambda z)$ for $z\in D(0,\varepsilon)$. The map $F$ is PCF, so it is polarized by $\cO(1)$ and $\PC(F)$ is a branch boundary for $F$ (Section~\ref{ssec:branch-boundary}). Thus $\psi$ contradicts Corollary~\ref{cor:no-rotation-disc} applied to $F$ on $\PP^N$. Hence every eigenvalue of $D_qg$ has modulus greater than one. By \eqref{eq:core-spectrum}, if $\lambda\neq0$ is an eigenvalue of $D_pf^a$, then $\lambda^b$ is an eigenvalue of $D_qg$, so $|\lambda|>1$.
\end{proof}

\begin{remark}\label{rem:routes}
For the projective-space case, both arguments use preservation of branch boundaries only on normalizations of invariant subvarieties of projective space, where it is due to Sun and Zhao \cite[Lemmas~2.3.2 and~A.0.2]{SZ26b}. They therefore do not require Appendix~\ref{app:whitney}. Indeed, any rotation disc obtained on such a normalization descends through its finite map to a nonconstant rotation disc for an iterate of the original PCF endomorphism.
\end{remark}

\section{Second route: proof of Theorem~C}\label{sec:route-two}

Sections~\ref{ssec:fitting}--\ref{ssec:local-aut} reduce Theorem~C to a fixed point at which the map is a local biholomorphism, which is treated in Section~\ref{ssec:claimS}.

\subsection{The invertible part of a local endomorphism}\label{ssec:fitting}

Let $(R,\fm)$ be a complete Noetherian local $\CC$-algebra with residue field $\CC$, and $\sigma$ a local endomorphism of $\CC$-algebras of $R$. A \emph{formal arc} of $R$ is a local homomorphism of $\CC$-algebras $R\to\CC[[t]]$. The map $\sigma$ induces an endomorphism of each finite dimensional space $R/\fm^k$, with Fitting decomposition $R/\fm^k=R_k^{\rm nil}\oplus R_k^{\rm inv}$ into $\sigma$-stable subspaces, $\sigma$ being nilpotent on $R_k^{\rm nil}$ and bijective on $R_k^{\rm inv}$. Let $\fn=\fn(R,\sigma)$ be the set of $h\in R$ with $h\bmod\fm^k\in R_k^{\rm nil}$ for every $k\geq1$. If $p$ is a periodic point, of period $a$, of an endomorphism $f$ of $\PP^N$, with $D_pf^a$ not nilpotent, then for $R=\hcO_{\PP^N,p}$ and $\sigma=(f^a)^*$ the ideal $\fn$ is that of the formal invariant germ of \cite[Section~3.1]{SZ26b}, by \cite[Lemma~3.1.1]{SZ26b} and parts~(1) and~(2) of the following lemma; part~(4) is the argument of \cite[proof of Lemma~3.2.4]{SZ26b}.

\begin{lemma}\label{lem:fitting}
\begin{enumerate}[label=\textup{(\arabic*)}]
\item The set $\fn$ is a $\sigma$-stable ideal of $R$. An element $h\in R$ lies in $\fn$ if and only if, for every $k$, there is $j$ with $\sigma^j(h)\in\fm^k$. For every $k\geq1$ there is $j_k$ such that $\sigma^j(\fn)\subset\fm^k$ for $j\geq j_k$. Moreover $\fn(R,\sigma^\kappa)=\fn(R,\sigma)$ for every $\kappa\geq1$.
\item The endomorphism of $R/\fn$ induced by $\sigma$ is an automorphism, and $\fn=0$ if and only if $\sigma$ is an automorphism of $R$.
\item The image of $\fn$ in $\fm/\fm^2$ is the generalized kernel of $\sigma$ on $\fm/\fm^2$. Hence the eigenvalues of $\sigma$ on the cotangent space $\fm/(\fn+\fm^2)$ of $R/\fn$ are the nonzero eigenvalues of $\sigma$ on $\fm/\fm^2$, with multiplicities.
\item Let $\beta$ be a formal arc of $R$ such that, for every $n\geq1$, there is a formal arc $\beta_n$ of $R$ with $\beta=\beta_n\circ\sigma^n$. Then $\beta(\fn)=0$.
\end{enumerate}
\end{lemma}

\begin{proof}
(1) The space $R_k^{\rm nil}$ is the generalized kernel of $\sigma$ on $R/\fm^k$, which is also that of $\sigma^\kappa$, and it is the kernel of $\sigma^j$ for $j\geq j_k:=\dim R/\fm^k$. As $\sigma(\fm^k)\subset\fm^k$, an element $h$ lies in $\fn$ if and only if, for every $k$, some $j$ satisfies $\sigma^j(h)\in\fm^k$, and then $\sigma^j(h)\in\fm^k$ for all $j\geq j_k$. This gives the second, third and fourth assertions. As $\sigma^j(hh')=\sigma^j(h)\sigma^j(h')$, the set $\fn$ is an ideal, and it is $\sigma$-stable.

(2) The equivariant projections $R/\fm^{k+1}\to R/\fm^k$ are onto, so they map $R^{\rm nil}_{k+1}$ onto $R^{\rm nil}_k$ and $R^{\rm inv}_{k+1}$ onto $R^{\rm inv}_k$. As $R=\varprojlim R/\fm^k$, we get $\fn=\varprojlim R^{\rm nil}_k$, the projection $\fn\to R_k^{\rm nil}$ is onto for every $k$, and $R/\fn\simeq\varprojlim R^{\rm inv}_k$ equivariantly, so $\sigma$ is bijective on $R/\fn$. If $\fn=0$, then $R_k^{\rm nil}=0$ for all $k$, so $\sigma$ is bijective on every $R/\fm^k$, hence on $R$. Conversely, an automorphism satisfies $\sigma(\fm^k)=\fm^k$, so it is bijective on $R/\fm^k$ and $\fn=0$.

(3) We have $R/\fm^2=\CC\cdot1\oplus\fm/\fm^2$ with $\sigma(1)=1$, so $R_2^{\rm nil}$ is the generalized kernel of $\sigma$ on $\fm/\fm^2$, and $\fn$ maps onto it by the proof of~(2).

(4) Let $h\in\fn$ and $k\geq1$. For $n\geq j_k$ we have $\beta(h)=\beta_n(\sigma^n(h))\in\beta_n(\fm^k)\subset t^k\CC[[t]]$, since $\beta_n$ is local. Hence $\beta(h)=0$.
\end{proof}

\subsection{Reduction to a normalized invariant subvariety}\label{ssec:reduction}

Throughout this subsection, $g_0$ is a polarized endomorphism of a normal projective variety $X_0$, $B_0$ is a branch boundary for $g_0$, and $q_0\in X_0$ is a fixed point of $g_0$. We set $R:=\hcO_{X_0,q_0}$, with maximal ideal $\fm$, $\sigma:=g_0^*$, $\fn:=\fn(R,\sigma)$ and $R_{\mathrm{inv}}:=R/\fn$. By Zariski's theorem on analytic normality, $R$ is a normal domain \cite[Theorem~32.2(i) and the Corollary of Theorem~32.6]{Mat89}.

We use twice that a domain $D$ which is a finite module over a complete Noetherian local ring $(A,\fm_A)$ is a complete local ring whose residue field is finite over $A/\fm_A$ \cite[\S9, Lemma~2, and Theorems~8.7 and~8.15]{Mat89}, since a domain has no idempotent other than $0$ and $1$.

\begin{lemma}\label{lem:arc}
Let $\widetilde R$ be a complete Noetherian local domain containing $\CC$, with residue field $\CC$ and $\dim\widetilde R\geq1$. For every nonzero $c\in\widetilde R$, there is a local homomorphism of $\CC$-algebras $\beta\colon\widetilde R\to\CC[[\tau]]$ with $\beta(c)\neq0$.
\end{lemma}

\begin{proof}
By Cohen's structure theorem \cite[Theorem~29.4(iii) and its proof]{Mat89}, applied with the coefficient field $\CC$, $\widetilde R$ is finite over a subring $P=\CC[[u_1,\dots,u_\varrho]]$, with $\varrho=\dim\widetilde R$. Let $T^l+p_1T^{l-1}+\dots+p_l$ be the minimal polynomial of $c$ over the fraction field of $P$. Its coefficients lie in $P$ \cite[Theorem~9.2]{Mat89}, as $P$ is normal and $c$ is integral over $P$, and $p_l\neq0$, as the polynomial is irreducible and $c\neq0$. Moreover $p_l\in c\widetilde R$. Choose $v\in\CC^\varrho\setminus\{0\}$ with $p_l(\tau v)\neq0$ in $\CC[[\tau]]$. The kernel $\mathfrak q$ of $P\to\CC[[\tau]]$, $u\mapsto\tau v$, is a prime ideal with $P/\mathfrak q\simeq\CC[[\tau]]$. By the lying over theorem \cite[Theorem~9.3(i)]{Mat89}, there is a prime $\mathfrak A$ of $\widetilde R$ with $\mathfrak A\cap P=\mathfrak q$. The ring $\widetilde R/\mathfrak A$ is a complete local domain of dimension one, finite over $P/\mathfrak q$. Its integral closure $D$ is the integral closure of $P/\mathfrak q\simeq\CC[[\tau]]$ in a finite extension of $\CC((\tau))$, which is separable, so it is finite over $P/\mathfrak q$ \cite[\S33, Lemma~1]{Mat89}. By the fact above, $D$ is a complete local domain with residue field $\CC$; it is normal of dimension one, hence a discrete valuation ring. If $t$ is a uniformizer, $D$ is finite over its subring $\CC[[t]]$ \cite[Theorem~29.4(iii) and its proof]{Mat89}, and $D=\CC[[t]]$ by Nakayama's lemma, as $D/tD=\CC$. Let $\beta$ be the composite $\widetilde R\to\widetilde R/\mathfrak A\hookrightarrow D=\CC[[t]]$, written in the variable $\tau$. It is local, since the maximal ideal of $D$ lies over that of $\widetilde R/\mathfrak A$ \cite[\S9, Lemma~2]{Mat89}. As $P/\mathfrak q\to\widetilde R/\mathfrak A$ is injective and $p_l\notin\mathfrak q$, we get $\beta(p_l)\neq0$, hence $\beta(c)\neq0$.
\end{proof}

\begin{lemma}\label{lem:reduced}
The ring $R_{\mathrm{inv}}$ is reduced.
\end{lemma}

\begin{proof}
For a formal arc $\beta$ of $R$ and $h\in R$, let $\ord\beta(h)\in\{0,1,\dots,\infty\}$ be the order of the power series $\beta(h)$. For $K\geq0$, let $\mathfrak a_K$ be the set of $h\in R$ with $\ord\beta(h)\geq K$ for every formal arc $\beta$ of $R$. The ideals $\mathfrak a_K$ form a decreasing sequence, and $\fm^K\subset\mathfrak a_K$ for every $K$.

We claim that $\bigcap_K\mathfrak a_K=0$. Indeed, if $\dim R=0$, then $R=\CC$ and $\mathfrak a_1=0$. If $\dim R\geq1$ and $h\neq0$, Lemma~\ref{lem:arc} gives a formal arc $\beta$ with $\beta(h)\neq0$, so $h\notin\mathfrak a_K$ for $K>\ord\beta(h)$.

For every $k$, there is $K$ with $\mathfrak a_K\subset\fm^k$. This is a lemma of Chevalley \cite[Exercise~8.7]{Mat89}, applied to the decreasing sequence $(\mathfrak a_K)$ of ideals of the complete local ring $R$, whose intersection is $0$.

Let $a\in R$ with $a^i\in\fn$ for some $i\geq1$, and let $k\geq1$. Choose $K$ with $\mathfrak a_K\subset\fm^k$. By Lemma~\ref{lem:fitting}(1), there is $j$ with $\sigma^j(a^i)\in\fm^{iK}$. Set $b:=\sigma^j(a)$. For every formal arc $\beta$ of $R$, we have $i\ord\beta(b)=\ord\beta(b^i)\geq iK$, since $\CC[[t]]$ is a domain and $\beta$ is local. Thus $b\in\mathfrak a_K\subset\fm^k$. By Lemma~\ref{lem:fitting}(1) again, $a\in\fn$.
\end{proof}

\begin{proposition}\label{prop:reduction}
Let $\lambda\neq0$ be an eigenvalue of $g_0^*$ on $V_{q_0}$. There exist $\kappa\geq1$ and an irreducible closed subvariety $Z\subset X_0$ containing $q_0$, with $\dim Z\geq1$ and $g_0^\kappa(Z)=Z$. Let $\nu\colon X\to Z$ be its normalization and $g$ the map induced by $g_0^\kappa$ \textup(Section~\textup{\ref{ssec:invariant}}\textup). There is a fixed point $q\in\nu^{-1}(q_0)$ of $g$ with the following properties.
\begin{enumerate}[label=\textup{(\arabic*)}]
\item The map $g$ is a polarized endomorphism of $X$ and admits a branch boundary.
\item There are a complete Noetherian local $\CC$-algebra $S$, a local endomorphism $\tilde\sigma$ of $S$, and a surjective local homomorphism $\pr\colon\hcO_{X,q}\to S$ such that $\pr\circ g^*=\tilde\sigma\circ\pr$. There is a nonzero element $a'\in\fm_S$ such that, for every $J\geq1$, some $k\geq1$ satisfies $(\tilde\sigma-\lambda^\kappa)^ka'\in\fm_S^J$.
\item For every nonzero $h\in\cOa_{X,q}$, there is a formal arc $\hat\alpha_0$ of $\hcO_{X,q}$ with $\hat\alpha_0(h)\neq0$ such that, for every $n\geq1$, $\hat\alpha_0=\hat\alpha_n\circ(g^*)^n$ for some formal arc $\hat\alpha_n$ of $\hcO_{X,q}$.
\end{enumerate}
\end{proposition}

\begin{proof}
\emph{Step 1: a branch of $\mathrm{Spf}\,R_{\mathrm{inv}}$.} By Lemma~\ref{lem:fitting}, $\fn$ does not change if $g_0$ is replaced by an iterate, and $\sigma$ induces an automorphism of $R_{\mathrm{inv}}$, which permutes the finitely many minimal primes of $R_{\mathrm{inv}}$. Choose $\kappa\geq1$ such that $\sigma^\kappa$ fixes each of them. We replace $g_0$ by $g_0^\kappa$, which is polarized with the branch boundary $B_0$, and $\lambda$ by $\lambda^\kappa$. So we assume that $\kappa=1$, $g_0$ and $\lambda$ standing for $g_0^\kappa$ and $\lambda^\kappa$.

Let $\mathcal E_J$ be the generalized $\lambda$-eigenspace of $\sigma$ on $\fm_{R_{\mathrm{inv}}}/\fm_{R_{\mathrm{inv}}}^J$. Equivariant surjections between finite dimensional spaces map generalized eigenspaces onto generalized eigenspaces, so the projections $\mathcal E_{J+1}\to\mathcal E_J$ are onto, and $\mathcal E_2\neq0$ by Lemma~\ref{lem:fitting}(3). Choosing successive compatible lifts of a nonzero element of $\mathcal E_2$ and using completeness, we obtain $\tilde a\in\fm_{R_{\mathrm{inv}}}\setminus\fm_{R_{\mathrm{inv}}}^2$ whose class modulo $\fm_{R_{\mathrm{inv}}}^J$ lies in $\mathcal E_J$ for every $J$. By Lemma~\ref{lem:reduced}, the minimal primes of $R_{\mathrm{inv}}$ have zero intersection, so one of them, $\fP$, does not contain $\tilde a$. Set $R':=R_{\mathrm{inv}}/\fP$, a complete local domain with the automorphism $\sigma'$ induced by $\sigma$, and let $a'\neq0$ be the image of $\tilde a$. Then $a'\in\fm_{R'}$, so $\dim R'\geq1$, and for every $J$ there is $k$ with $(\sigma'-\lambda)^ka'\in\fm_{R'}^J$.

\emph{Step 2: the subvariety $Z$.} Let $\fp$ be the kernel of $\cOa_{X_0,q_0}\to R'$, and let $Z\subset X_0$ be the irreducible closed subvariety whose ideal at $q_0$ is $\fp$. The relation $g_0^*\fp\subset\fp$ implies that $g_0$ maps a neighbourhood of $q_0$ in $Z$ into $Z$. Thus $g_0(Z)\subset Z$. Since $g_0$ is finite, $\dim g_0(Z)=\dim Z$, and $g_0(Z)=Z$. The completion $\hcO_{Z,q_0}=R/\fp R$ maps onto $R'$, so $\dim Z\geq\dim R'\geq1$. Let $\nu\colon X\to Z$ and $g$ be as in the statement. Assertion~(1) follows from Proposition~\ref{prop:normalization-boundary}.

\emph{Step 3: the point $q$.} Let $\widetilde R$ be the integral closure of $R'$ in its fraction field. It is a finite $R'$-module \cite[\S33, Lemma~1 and p.~263]{Mat89}, hence a complete local domain with residue field $\CC$, by the fact recalled before Lemma~\ref{lem:arc}. The automorphism $\sigma'$ extends to an automorphism $\tilde\sigma$ of $\widetilde R$. The injection $\cOa_{Z,q_0}=\cOa_{X_0,q_0}/\fp\hookrightarrow R'$ extends to an embedding $\omega$ of $\CC(X)=\CC(Z)$ into the fraction field of $\widetilde R$. Let $\widetilde\cO$ be the integral closure of $\cOa_{Z,q_0}$ in $\CC(Z)$. It is a semi-local ring, whose maximal ideals are the prime ideals over the maximal ideal of $\cOa_{Z,q_0}$ \cite[\S9, Lemma~2]{Mat89}. They correspond to the points of $\nu^{-1}(q_0)$, and $\cOa_{X,x}$ is the localization of $\widetilde\cO$ at the maximal ideal of $x$ \cite[Chapter~III, \S8, Theorem~3 and its Corollary]{Mum99}, \cite[Proposition~5.12]{AM69}.

As $\widetilde R$ is integrally closed, $\omega(\widetilde\cO)\subset\widetilde R$. The preimage of $\fm_{\widetilde R}$ in $\widetilde\cO$ is a prime ideal with residue field contained in $\CC$, hence a maximal ideal; let $q\in\nu^{-1}(q_0)$ be the corresponding point. The elements of $\widetilde\cO$ outside this ideal are mapped to units, so $\omega$ maps $\cOa_{X,q}$ into $\widetilde R$ by a local homomorphism, which is injective as the restriction of a field embedding.

Since $g_0(q_0)=q_0$ and $\nu\circ g=g_0\circ\nu$, the map $g$ preserves $\nu^{-1}(q_0)$, and $g^*$ preserves $\widetilde\cO$. The identity $\omega\circ g^*=\tilde\sigma\circ\omega$ holds on $\cOa_{Z,q_0}$, hence on $\CC(Z)$. For $h\in\widetilde\cO$, we have $h(g(q))=0$ if and only if $\omega(g^*h)=\tilde\sigma(\omega(h))$ lies in $\fm_{\widetilde R}$, that is, if and only if $h(q)=0$. So $g(q)$ and $q$ have the same maximal ideal in $\widetilde\cO$, and $g(q)=q$.

\emph{Step 4: the complete quotient $S$.} The local homomorphism $\omega$ extends by continuity to $\hat\omega\colon\hcO_{X,q}\to\widetilde R$, with $\hat\omega\circ g^*=\tilde\sigma\circ\hat\omega$. Set $S:=\hcO_{X,q}/\ker\hat\omega$ and let $\pr$ be the quotient map. This is a complete Noetherian local ring, identified with the image of $\hat\omega$ in $\widetilde R$. It is preserved by $\tilde\sigma$, and $\tilde\sigma(\fm_S)=\pr(g^*\widehat\fm_q)\subset\pr(\widehat\fm_q)=\fm_S$, so the restriction of $\tilde\sigma$ to $S$ is local. The composite $\hcO_{Z,q_0}\to\hcO_{X,q}\to\widetilde R$ agrees with $\hcO_{Z,q_0}\to R'\hookrightarrow\widetilde R$: both are continuous and agree on $\cOa_{Z,q_0}$. It follows that $R'\subset S$ and $\fm_{R'}\subset\fm_S$. Therefore $a'\in\fm_S$ and $(\tilde\sigma-\lambda)^ka'\in\fm_{R'}^J\subset\fm_S^J$, as required in~(2).

\emph{Step 5: formal arcs with lifts of every order.} Let $h\in\cOa_{X,q}$ be nonzero. Then $\omega(h)\neq0$. As $\widetilde R$ is finite over $R'$, $\dim\widetilde R=\dim R'\geq1$, and Lemma~\ref{lem:arc} gives a local homomorphism $\beta\colon\widetilde R\to\CC[[\tau]]$ with $\beta(\omega(h))\neq0$. The maps $\hat\alpha_n:=\beta\circ\tilde\sigma^{-n}\circ\hat\omega$, $n\geq0$, are local homomorphisms, hence formal arcs. They satisfy $\hat\alpha_0(h)\neq0$ and $\hat\alpha_n\circ(g^*)^n=\beta\circ\tilde\sigma^{-n}\circ\tilde\sigma^n\circ\hat\omega=\hat\alpha_0$.
\end{proof}

\subsection{The induced map is a local biholomorphism}\label{ssec:local-aut}

In this subsection, $X$, $m=\dim X\geq1$, $g$, $B$ and the fixed point $q$, possibly singular, are as in Section~\ref{ssec:arcs}, and we use the notation of Appendix~\ref{app:arcs}: the objects $\Pi$, $\Delta$, $\delta$, $B_1$, $\Omega_0$, $U$, $V_n$, $\Pi_n$, $U^\circ$ and $V_n^\circ$ defined there. The proof of the next proposition follows the proof of \cite[Theorem~3.2.1]{SZ26b}: arcs which lift through every iterate annihilate the non-invertible part (Lemma~\ref{lem:fitting}(4); \cite[Lemma~3.2.4]{SZ26b}), and a family of such arcs of rank $m$ (Proposition~\ref{prop:family}; \cite[Lemma~4.2.1]{SZ26b}) forces this part to vanish. Step~3 is new.

\begin{proposition}\label{prop:local-aut}
Assume that there is a formal arc $\hat\alpha_0$ of $\hcO_{X,q}$ with $\hat\alpha_0(\delta\circ\Pi)\neq0$ such that, for every $n\geq1$, $\hat\alpha_0=\hat\alpha_n\circ(g^*)^n$ for some formal arc $\hat\alpha_n$ of $\hcO_{X,q}$. Then $\fn(\hcO_{X,q},g^*)=0$, and $g$ is a local biholomorphism at $q$.
\end{proposition}

\begin{proof}
\emph{Step 1: a family of arcs.} Let $\varepsilon$, $\BB$ and $\Xi\colon D(0,\varepsilon)\times\BB\to X$ be given by Proposition~\ref{prop:family}. Thus $\Xi(0,s)=q$, and $\Xi$ maps $D^*(0,\varepsilon)\times\BB$ into $X\setminus B_1\subset\Reg X$ with rank $m$. Moreover, for $n\geq1$, every arc $\Xi(\cdot,s)$ is equal to $g^n\circ\eta_{s,n}$ for a holomorphic arc $\eta_{s,n}$ with $\eta_{s,n}(0)=q$. Identify a neighbourhood $W$ of $q$ in $X$ with a closed analytic subset of an open set of $\CC^\varkappa$, with coordinates $z_1,\dots,z_\varkappa$ and $q=0$. Shrinking $\varepsilon$ and $\BB$, we may assume that $\Xi$ takes values in $W$.

The functions $z_k\circ\Xi$ vanish on $\{0\}\times\BB$, so their expansions in $t$ lie in $t\,\cO(\BB)[[t]]$. Substituting them for $z$ defines a ring homomorphism $\CC[[z]]\to\cO(\BB)[[t]]$, which kills the convergent equations of $W$. As $\hcO_{X,q}$ is the quotient of $\CC[[z]]$ by the ideal generated by these equations, we obtain a homomorphism $\Xi^*\colon\hcO_{X,q}\to\cO(\BB)[[t]]$. For $s\in\BB$, its composite with the evaluation at $s$ is the formal arc given by the Taylor expansion of $\Xi(\cdot,s)$ at $0$. It is equal to $\hat\eta_{s,n}\circ(g^*)^n$, where $\hat\eta_{s,n}$ is the formal arc given by $\eta_{s,n}$, for every $n\geq1$. By Lemma~\ref{lem:fitting}(4), it vanishes on $\fn_X:=\fn(\hcO_{X,q},g^*)$. A power series in $t$ whose coefficients vanish at every point of $\BB$ is zero, so $\Xi^*(\fn_X)=0$. The kernel $\mathfrak K$ of $\Xi^*$ is a prime ideal, as $\cO(\BB)[[t]]$ is a domain, and $\fn_X\subset\mathfrak K$.

\emph{Step 2: $\mathfrak K=0$.} As in Section~\ref{ssec:reduction}, $\hcO_{X,q}$ is a normal domain; its dimension is $m$. Assume that $\mathfrak K\neq0$. Then the complete local domain $\Lambda:=\hcO_{X,q}/\mathfrak K$ has dimension $\varrho'<m$, and by Cohen's structure theorem \cite[Theorem~29.4(iii) and its proof]{Mat89} it is finite over a subring $P'=\CC[[v_1,\dots,v_{\varrho'}]]$, where the $v_i$ are images of elements of $\widehat\fm_q$, which we also denote by $v_i$. For each $j$, let $Q_j(w;v)\in P'[w]$ be the minimal polynomial of the class of $z_j$ over the fraction field of $P'$. It is monic with coefficients in $P'$ \cite[Theorem~9.2]{Mat89}, and since the characteristic is zero, $\partial_wQ_j$ does not vanish at the class of $z_j$. Thus $Q_j(z_j;v)\in\mathfrak K$ and $\partial_wQ_j(z_j;v)\notin\mathfrak K$. Set $\mathsf Z_j:=\Xi^*(z_j)$ and $\mathsf V_i:=\Xi^*(v_i)$, elements of $t\,\cO(\BB)[[t]]$. Then $Q_j(\mathsf Z_j;\mathsf V)=0$, and $\partial_wQ_j(\mathsf Z_j;\mathsf V)\neq0$. Differentiating with respect to $\tau\in\{t,s_1,\dots,s_m\}$, the formal chain rule gives
\[
\partial_wQ_j(\mathsf Z_j;\mathsf V)\,\partial_\tau\mathsf Z_j=-\sum_{i=1}^{\varrho'}\partial_{v_i}Q_j(\mathsf Z_j;\mathsf V)\,\partial_\tau\mathsf V_i.
\]
Over the fraction field of $\cO(\BB)[[t]]$, the rows $(\partial_\tau\mathsf Z_j)_\tau$ lie in the span of the $\varrho'$ rows $(\partial_\tau\mathsf V_i)_\tau$, so that $\operatorname{rank}(\partial_\tau\mathsf Z_j)_{j,\tau}\leq\varrho'$. This matrix is the Taylor expansion in $t$ of the Jacobian matrix of $z\circ\Xi$. Every minor of order $\varrho'+1$ vanishes as a formal series, hence as a holomorphic function on $D(0,\varepsilon)\times\BB$. It follows that $z\circ\Xi$ has rank at most $\varrho'<m$.

On $D^*(0,\varepsilon)\times\BB$, however, $\Xi$ has rank $m$ and takes values in $\Reg X$, where $z$ is an immersion. This contradiction proves $\mathfrak K=0$, and therefore $\fn_X=0$.

\emph{Step 3: $g$ is a local biholomorphism.} By Lemma~\ref{lem:fitting}(2), $g^*$ is an automorphism of $\hcO_{X,q}$. Let $\cO:=\cO_{X,q}$ be the analytic local ring, with maximal ideal $\fm$, and $\sigma:=g^*\colon\cO\to\cO$. As $q$ is isolated in $g^{-1}(q)$, the ideal $\sigma(\fm)\cO$ has zero set $\{q\}$, so it contains a power $\fm^k$ by the Hilbert--R\"uckert Nullstellensatz \cite[Chapter~I, Theorem~1.72]{GLS07}, and $\cO$ is a finite module over itself via $\sigma$, by (N3). Let $x\in\fm$. As $g^*$ maps $\widehat\fm$ onto $\widehat\fm$, we can write $x=g^*(y)$ with $y\in\widehat\fm$. Choose $y'\in\fm$ with $y-y'\in\widehat\fm^k$. Then $x-\sigma(y')\in g^*(\widehat\fm^k)\cap\cO=\fm^k\subset\sigma(\fm)\cO$, as the automorphism $g^*$ preserves $\widehat\fm^k$ and $\widehat\fm^k\cap\cO=\fm^k$. Thus $\fm=\sigma(\fm)\cO$, so $\cO/\sigma(\fm)\cO=\CC$, and Nakayama's lemma shows that $\sigma$ is surjective. It is injective, since $g^*$ is injective on $\hcO_{X,q}\supset\cO$. Hence $\sigma$ is an automorphism of the analytic algebra $\cO_{X,q}$, and the germ of $g$ at $q$ is an isomorphism \cite[Chapter~I, Remark~1.47.1(2)]{GLS07}.
\end{proof}

\subsection{Fixed points at which the map is a local biholomorphism}\label{ssec:claimS}

\begin{proposition}\label{thm:claimS}
Let $g$ be a polarized endomorphism of a normal projective variety $X$ which admits a branch boundary $B$, and let $q$ be a fixed point of $g$ at which $g$ is a local biholomorphism. Then every eigenvalue of the automorphism $g_q^*$ of $V_q$ induced by $g$ has modulus greater than one.
\end{proposition}

Throughout this subsection, $X$, $g$, $B$ and $q$ are as in Proposition~\ref{thm:claimS}. We may assume that $m:=\dim X\geq1$, since otherwise $V_q=0$. We use the objects of Section~\ref{ssec:arcs} for $(X,g,B,q)$.

The next lemma extends \cite[Lemma~2.3.3]{SZ26b} (Proposition~\ref{prop:sz-inverse}) to a fixed point which may be singular; the proof is theirs, with the projection $\Pi$ of Section~\ref{ssec:arcs} in place of a hypersurface germ containing $B$, and the injectivity of $g^{n+k}$ on $V_{n+k}$ in place of the identity theorem.

\begin{lemma}\label{lem:branches}
For every $n\geq0$, the map $g^n$ maps $V_n$ biholomorphically onto $U$; let $b_n\colon U\to V_n$ be its inverse. Then $b_n(q)=q$, the map induced by $b_n$ on the cotangent space $V_q$ is $(g_q^*)^{-n}$, and, if $U_1\subset U$ is a connected neighbourhood of $q$ with $b_k(U_1)\subset U$, then $b_n\circ b_k=b_{n+k}$ on $U_1$ for every $n\geq0$. Moreover, every sequence in $(b_n)_{n\geq0}$ has a subsequence which converges locally uniformly on $U$ to a holomorphic map $U\to X$, and for every neighbourhood $W$ of $q$ there is a neighbourhood $U_0\subset U$ of $q$ with $b_n(U_0)\subset W$ for every $n$. In particular, $q$ has bounded inverse branches \textup(Definition~\textup{\ref{def:bib}}\textup).
\end{lemma}

\begin{proof}
The map $g^n\colon V_n\to U$ is proper with finite fibres, as $V_n$ is open and closed in $g^{-n}(U)$. By Lemma~\ref{lem:projection}, $V_n\cap g^{-n}(q)\subset\Pi_n^{-1}(0)=\{q\}$, and $g^n$ restricts to a finite covering $V_n^\circ\to U^\circ$ of connected manifolds. Here $V_n^\circ=V_n\cap g^{-n}(U^\circ)$, since $\Pi_n=\Pi\circ g^n$, $B_1=\Pi^{-1}(\Delta)$ and $g^n(V_n)\subset U$. As $g$ is a local biholomorphism at the fixed point $q$, there is a neighbourhood $W_1$ of $q$ on which $g^n$ is injective, and by properness an open set $W_0\ni q$ with $V_n\cap g^{-n}(W_0)\subset W_1$. The set $U^\circ$ is dense in $U$, as $B_1$ is a proper algebraic subset of the irreducible variety $X$. So a point of $W_0\cap U^\circ$ has exactly one preimage in $V_n$, and the covering $V_n^\circ\to U^\circ$ has one sheet.

We claim that $g^n$ is injective on $V_n$. The spaces $V_n$ and $U$ are open in $X$, and $U$ is normal, so $g^n|_{V_n}$ is open by (N2). If $g^n(a)=g^n(b)$ with $a\neq b$, choose disjoint open neighbourhoods $W_a$, $W_b$ of $a$, $b$ in $V_n$. The open set $g^n(W_a)\cap g^n(W_b)$ contains a point of $U^\circ$, which has two preimages in $V_n^\circ$, a contradiction. Moreover $g^n(V_n)=U$: this image is closed in $U$, as $g^n|_{V_n}$ is proper, and it contains the dense subset $U^\circ$, the base of the covering $V_n^\circ\to U^\circ$. Thus $g^n\colon V_n\to U$ is a continuous proper bijection, hence a homeomorphism. Its inverse $b_n$ is holomorphic on $U^\circ$, since $U^\circ\cap\Br(g^n)=\emptyset$. In local embeddings of $V_n$, its components are continuous and holomorphic outside the nowhere dense analytic set $B_1\cap U$, hence holomorphic by (N1). Near $q$, $b_n$ is the inverse of the germ of $g^n$, so it induces $(g_q^*)^{-n}$ on $V_q$.

The set $b_n\circ b_k(U_1)$ is connected, contains $q$ and lies in $g^{-n-k}(U)$, hence in $V_{n+k}$. As $g^{n+k}$ is injective on $V_{n+k}$ and $g^{n+k}\circ b_n\circ b_k=\Id_{U_1}$, we get $b_n\circ b_k=b_{n+k}$ on $U_1$.

By Proposition~\ref{prop:fakhruddin} and the paragraph following it, we may embed $X$ into $\PP^L$ so that $g$ extends to an endomorphism $\Phi$ of $\PP^L$ of algebraic degree $M\geq2$. Viewed as maps into $\PP^L$, the $b_n$ satisfy $\Phi^n\circ b_n=\iota_U$, the inclusion of $U$. By Theorem~\ref{thm:ueda-normality}, which holds for any complex space as source, applied to $U$ and $\psi_0=\iota_U$, the family $(b_n)$ is normal, and its limits take values in the closed set $X$. Since this family is normal and $b_n(q)=q$ for every $n$, it is equicontinuous at $q$, which gives the last assertion.
\end{proof}

For $\eta>0$, let $\mathcal D_\eta:=\{w\in\CC^c:|w_i|<\eta\text{ for all }i\}$.

\begin{lemma}\label{lem:TOR}
Let $A=\operatorname{diag}(\mu_1,\dots,\mu_c)$ with $|\mu_i|=1$, and let $\Sigma_\eta\subset\mathcal D_\eta$ be a closed analytic subset of positive dimension at $0\in\Sigma_\eta$, with $A(\Sigma_\eta)=\Sigma_\eta$. There is $\ell\geq1$ such that one of the following holds.
\begin{enumerate}[label=\textup{(\alph*)}]
\item $A^\ell$ fixes infinitely many points of $\Sigma_\eta$.
\item There are $\mu'\in S^1$, $r>1$ and a nonconstant holomorphic map $\psi\colon D(0,r)\to\Sigma_\eta$ with $\psi(0)=0$ and $A^\ell\psi(z)=\psi(\mu'z)$ for $|z|<r$.
\end{enumerate}
\end{lemma}

\begin{proof}
Let $K$ be the closure of $\{A^n:n\in\ZZ\}$ in the compact group of diagonal unitary matrices; $\Sigma_\eta$ is $K$-invariant, being closed in $\mathcal D_\eta$ and invariant under $A^{\pm1}$. The group $K$ is a compact abelian Lie group, its identity component $K^0$ is open of finite index $\ell$ and is a torus $(S^1)^s$, and $K^0$ acts by characters: there are $\chi_1,\dots,\chi_c\in\ZZ^s$ with $k\cdot w=(k^{\chi_i}w_i)_i$ for $k\in K^0\simeq(S^1)^s$, in multi-index notation \cite[Theorem~20.12, Propositions~7.14 and~7.15]{Lee13}, \cite[Chapter~VI, \S XIII, Proposition~1, Chapter~IV, \S XIII, Proposition~1, and p.~213]{Che46}. As $K/K^0$ is generated by the class of $A$, we have $A^\ell\in K^0$. For $u\in\ZZ^s$ and $z\in\CC^*$, let $\vartheta_u(z)$ act by $\vartheta_u(z)\cdot w:=(z^{\langle\chi_i,u\rangle}w_i)_i$; for $z\in S^1$ it is an element of $K^0$.

We claim that, for $x\in\Sigma_\eta$ and $u\in\ZZ^s$, the set $A_x:=\{z\in\CC^*:\vartheta_u(z)\cdot x\in\mathcal D_\eta\}$ is an annulus containing $S^1$, and that $\vartheta_u(z)\cdot x\in\Sigma_\eta$ for every $z\in A_x$. Indeed, each condition $|z|^{\langle\chi_i,u\rangle}|x_i|<\eta$ is void or bounds $|z|$ from above or from below, and it holds on $S^1$. Let $f\in\cO(\mathcal D_\eta)$ vanish on $\Sigma_\eta$. The function $z\mapsto f(\vartheta_u(z)\cdot x)$ is holomorphic on $A_x$ and vanishes on $S^1$, since $K^0$ preserves $\Sigma_\eta$, so it vanishes on $A_x$. By Cartan's Theorem~A on the polydisc $\mathcal D_\eta$, $\Sigma_\eta$ is the common zero set of such functions \cite[Theorem~7.2.11]{Hor90} (the polydisc is Stein by \cite[Corollary~2.5.8 and p.~116]{Hor90}), which proves the claim.

Let $F_0:=\{w:w_i=0\text{ whenever }\chi_i\neq0\}$, the fixed-point set of $K^0$. If $\Sigma_\eta\cap F_0$ is infinite near $0$, then (a) holds. Otherwise, as $\Sigma_\eta$ has positive dimension at $0$, there are points $x_n\in\Sigma_\eta\setminus F_0$ with $x_n\to0$. Choose $v\in\ZZ^s$ with $\langle\chi_i,v\rangle\neq0$ whenever $\chi_i\neq0$. Replacing $v$ by $-v$ and passing to a subsequence, we may assume that each $x_n$ has a nonzero coordinate $x_{n,i}$ with $\langle\chi_i,v\rangle>0$. For $n$ large, let $t_n\geq0$ be minimal with $\max\{e^{t_n\langle\chi_i,v\rangle}|x_{n,i}|:\langle\chi_i,v\rangle>0\}=\eta/2$. All coordinates of $\vartheta_v(e^{t_n})\cdot x_n$ have modulus at most $\max(\eta/2,|x_n|)<\eta$, so $y_n:=\vartheta_v(e^{t_n})\cdot x_n\in\Sigma_\eta$ by the claim. The coordinates $y_{n,i}$ with $\langle\chi_i,v\rangle\leq0$ tend to $0$, and the others have modulus at most $\eta/2$, one of them equal to $\eta/2$. A limit point $y$ of $(y_n)$ lies in $\mathcal D_\eta$, hence in $\Sigma_\eta$; it is nonzero, $|y_i|\leq\eta/2$ for every $i$, and $\langle\chi_i,v\rangle>0$ for every $i$ in the support $T(y):=\{i:y_i\neq0\}$.

Assume that the characters $\chi_i$, $i\in T(y)$, are not proportional. Then there is $u'\in\ZZ^s$ such that the ratio $\langle\chi_i,u'\rangle/\langle\chi_i,v\rangle$ is not constant on $T(y)$. Let $p/p'$, with $p'\geq1$, be its minimum, and $u:=p'u'-pv$. Then $\langle\chi_i,u\rangle\geq0$ on $T(y)$, with equality exactly on a proper nonempty subset $T'\subset T(y)$. So $\{0<|z|\leq1\}\subset A_y$ for this $u$, and $\vartheta_u(z)\cdot y\in\Sigma_\eta$ for $0<|z|\leq1$, by the claim. As $z\to0$, these points tend to the point $y'$ with $y'_i=y_i$ for $i\in T'$ and $y'_i=0$ otherwise. Thus $y'\in\Sigma_\eta$, and $y'$ has the properties of $y$, with a smaller support. After finitely many steps, we obtain such a point $y$ for which $\chi_i=n_i\chi_0$ for $i\in T(y)$, with $\chi_0\in\ZZ^s$ and integers $n_i\geq1$.

Set $\psi(z):=(z^{n_i}y_i)_{i\in T(y)}$, the other coordinates being $0$. It is defined and takes values in $\mathcal D_\eta$ on $D(0,r)$ for some $r>1$, as $|y_i|\leq\eta/2$. The character $k\mapsto k^{\chi_0}$ of $K^0$ is nontrivial, hence onto $S^1$. For $|z|=1$, choose $k\in K^0$ with $k^{\chi_0}=z$; then $\psi(z)=k\cdot y\in\Sigma_\eta$. As in the proof of the claim, $\psi(D(0,r))\subset\Sigma_\eta$. The map $\psi$ is nonconstant, $\psi(0)=0$, and with $\mu':=(A^\ell)^{\chi_0}\in S^1$, the value of this character at $A^\ell\in K^0$, we get $A^\ell\psi(z)=\psi(\mu'z)$, since $\mu_i^\ell=(A^\ell)^{\chi_i}=\mu'^{n_i}$ for $i\in T(y)$.
\end{proof}

\begin{proof}[Proof of Proposition~\ref{thm:claimS}]
By Lemma~\ref{lem:branches}, $q$ has bounded inverse branches $(b_n)$ on $U$, and we apply Section~\ref{sec:bib} to $Y:=X$ and $h:=g$, with the coordinates $z_1,\dots,z_\varkappa$, the sets $\beta(r)$ and the radii $r_0>r_1>r_2$ fixed there, and with $E^c$, $E^u$, $(E^c)^*$, $(E^u)^*$, $c$ and $\rho$ as in Proposition~\ref{prop:local-retract}. By Proposition~\ref{prop:local-retract}(1), every eigenvalue of $g_q^*$ has modulus at least one, and it suffices to show that $c=0$. We assume that $c\geq1$.

\emph{Step 1: the retract and its Zariski cotangent space.} Let $\Theta_1$ be the connected component of $\{x\in\beta(r_2):\rho(x)\in\beta(r_2)\}$ containing $q$, an open neighbourhood of $q$. The set $\rho(\Theta_1)$ is connected, contains $q$, and lies in $\{x\in\beta(r_2):\rho(x)\in\beta(r_2)\}$, since $\rho\circ\rho=\rho$ on $\beta(r_2)$; so $\rho(\Theta_1)\subset\Theta_1$. Let $\Sigma:=\{x\in\Theta_1:\rho(x)=x\}=\rho(\Theta_1)$, a closed analytic subset of $\Theta_1$, with its reduced structure. Let $\iota_\Sigma\colon\Sigma\to\Theta_1$ be the inclusion. As $\Theta_1$ is reduced, $\rho$ induces a holomorphic map $\rho_\Sigma\colon\Theta_1\to\Sigma$, and $\rho_\Sigma\circ\iota_\Sigma=\Id_\Sigma$, $\iota_\Sigma\circ\rho_\Sigma=\rho|_{\Theta_1}$. Let $V_{\Sigma,q}$ be the cotangent space of $\Sigma$ at $q$. On cotangent spaces, $\iota_\Sigma^*\circ\rho_\Sigma^*=\Id_{V_{\Sigma,q}}$ and $\rho_\Sigma^*\circ\iota_\Sigma^*=\rho^*$, which is the projection of $V_q$ onto $(E^c)^*$ along $(E^u)^*$, by Proposition~\ref{prop:local-retract}(2). Hence $\iota_\Sigma^*\colon V_q\to V_{\Sigma,q}$ is onto, with kernel $(E^u)^*$, and it restricts to an isomorphism $(E^c)^*\to V_{\Sigma,q}$. In particular $V_{\Sigma,q}\neq0$, so $q$ is not isolated in $\Sigma$, and $\dim_q\Sigma\geq1$.

\emph{Step 2: $b_k$ preserves $\Sigma$ near $q$.} By Proposition~\ref{prop:local-retract}(2), $b_k\circ\rho=\rho\circ b_k$ on $\beta(r_1)$ for every $k\geq0$. By Definition~\ref{def:bib}(3), there is $r_3\in(0,r_2]$ with $b_k(\beta(r_3))\subset\Theta_1$ for every $k$. For $x\in\Sigma\cap\beta(r_3)$, we get $\rho(b_k(x))=b_k(\rho(x))=b_k(x)$ with $b_k(x)\in\Theta_1$, so $b_k(x)\in\Sigma$.

\emph{Step 3: an equivariant embedding of $\Sigma$.} By Proposition~\ref{prop:local-retract}(1), there is a basis $v_1,\dots,v_c$ of $(E^c)^*$ with $g_q^*v_i=\mu_iv_i$, where $|\mu_i|=1$. Let $h_i$ be linear combinations of $z_1,\dots,z_\varkappa$ with $[h_i]=v_i$. For $n\geq1$, set
\[
\Psi_{n,i}:=\frac1n\sum_{k=0}^{n-1}\mu_i^k\,h_i\circ b_k
\quad\text{on }\beta(r_1).
\]
These functions are bounded by $\sup_{\beta(r_0)}|h_i|$, vanish at $q$, and satisfy
\[
[\Psi_{n,i}]
 =\frac1n\sum_{k=0}^{n-1}\mu_i^k(g_q^*)^{-k}v_i=v_i.
\]
For $x\in\beta(r_2)$, we have $b_1(x)\in\beta(r_1)$ and $b_k(b_1(x))=b_{k+1}(x)$, so that
\[
\Psi_{n,i}(b_1(x))
 =\mu_i^{-1}\left(\Psi_{n,i}(x)
   +\frac{\mu_i^nh_i(b_n(x))-h_i(x)}{n}\right).
\] By Lemma~\ref{lem:montel}, applied on $\beta(r_1)$, which has a countable basis, a subsequence of $(\Psi_{n,i})_n$, the same for all $i$, converges locally uniformly on $\beta(r_1)$ to a function $\Psi_i$. Then $\Psi_i\circ b_1=\mu_i^{-1}\Psi_i$ on $\beta(r_2)$, and $[\Psi_i]=v_i$ by Lemma~\ref{lem:jet}. Set $\Psi:=(\Psi_1,\dots,\Psi_c)$ and $A_c:=\operatorname{diag}(\mu_1,\dots,\mu_c)$, so that $\Psi\circ b_1=A_c^{-1}\circ\Psi$ on $\beta(r_2)$.

The germ of $\Psi|_\Sigma$ at $q$ is a closed embedding into $(\CC^c,0)$: its cotangent map sends $dw_i$ to $\iota_\Sigma^*v_i$, so it is onto $V_{\Sigma,q}$ by Step~1. Nakayama's lemma gives $\fm_{\Sigma,q}=(\Psi_1,\dots,\Psi_c)$. By (N3), $\cO_{\Sigma,q}$ is finite over $\cO_{\CC^c,0}$ via $\Psi^*$. A second application of Nakayama's lemma shows that this module is generated by $1$, so $\Psi^*$ is onto. Thus there are an open neighbourhood $\Sigma'\subset\Sigma\cap\beta(r_3)$ of $q$ in $\Sigma$ and an open set $W_0\ni0$ such that $\Psi$ maps $\Sigma'$ biholomorphically onto a closed analytic subset $\Sigma_0$ of $W_0$; let $\Upsilon\colon\Sigma_0\to\Sigma'$ be the inverse. By Step~2, there is an open neighbourhood $\Sigma''\subset\Sigma'$ of $q$ in $\Sigma$ with $b_1(\Sigma'')\subset\Sigma'$. Choose $\eta>0$ with $\mathcal D_\eta\subset W_0$ and $\Upsilon(\Sigma_0\cap\mathcal D_\eta)\subset\Sigma''$, and set $\Sigma_\eta:=\Sigma_0\cap\mathcal D_\eta$, a closed analytic subset of $\mathcal D_\eta$ of positive dimension at $0$.

Let $w\in\Sigma_\eta$ and $x:=\Upsilon(w)$. Then $b_1(x)\in\Sigma'$ and $\Psi(b_1(x))=A_c^{-1}w\in\mathcal D_\eta$, so $A_c^{-1}w\in\Sigma_\eta$ and $b_1(\Upsilon(w))=\Upsilon(A_c^{-1}w)$. Thus $A_c^{-1}(\Sigma_\eta)\subset\Sigma_\eta$. Choosing $l_j\to\infty$ with $A_c^{l_j}\to\Id$, we get $A_c=\lim A_c^{-(l_j-1)}$, so $A_c(\Sigma_\eta)\subset\Sigma_\eta$ as $\Sigma_\eta$ is closed in $\mathcal D_\eta$. Hence $A_c(\Sigma_\eta)=\Sigma_\eta$. Applying $g$ to the relation $b_1\circ\Upsilon=\Upsilon\circ A_c^{-1}$ at the point $A_cw$ gives $g\circ\Upsilon=\Upsilon\circ A_c$ on $\Sigma_\eta$, and therefore $g^k\circ\Upsilon=\Upsilon\circ A_c^k$ on $\Sigma_\eta$ for every $k\geq1$.

\emph{Step 4: conclusion.} Let $\ell$ be given by Lemma~\ref{lem:TOR}, applied to $A_c$ and $\Sigma_\eta$. In case~(a), $g^\ell$ fixes the points $\Upsilon(w)$, where $w\in\Sigma_\eta$ is fixed by $A_c^\ell$, which are infinitely many as $\Upsilon$ is injective; this contradicts Lemma~\ref{lem:finite-fix}. In case~(b), the map $\varphi:=\Upsilon\circ\psi$ is holomorphic and nonconstant on $D(0,r)$, and $g^\ell\circ\varphi(z)=\Upsilon(A_c^\ell\psi(z))=\Upsilon(\psi(\mu'z))=\varphi(\mu'z)$. Since $g^\ell$ is polarized and $B$ is a branch boundary for $g^\ell$, this contradicts Corollary~\ref{cor:no-rotation-disc} (when $g$ is induced by an iterate of a PCF endomorphism of $\PP^N$ on the normalization of an invariant subvariety, the corollary may be applied on $\PP^N$ instead, see Remark~\ref{rem:routes}). Hence $c=0$.
\end{proof}

\subsection{Conclusion of the proof}\label{ssec:proofC}

\begin{proof}[Proof of Theorem~C]
Let $B$ be a branch boundary for $g$, and let $x$ be a periodic point of $g$ of period $a$. The map $g^a$ is a polarized endomorphism of $X$, and $B$ is a branch boundary for it. Replacing $g$ by $g^a$, we may assume that $x$ is fixed. Let $\lambda\neq0$ be an eigenvalue of $D_xg$, that is, of $g^*$ on $V_x$ (Definition~\ref{def:multiplier}). We apply Proposition~\ref{prop:reduction} to $X_0:=X$, $g_0:=g$, $B_0:=B$, $q_0:=x$ and $\lambda$. Let $\kappa$, $Z$, $\nu\colon X'\to Z$, $q\in X'$, $S$, $\tilde\sigma$, $\pr$ and $a'$ be as there, and $g'$ the induced map. It suffices to show that $|\lambda^\kappa|>1$.

By Proposition~\ref{prop:reduction}(1), $g'$ is polarized and admits a branch boundary $B'$, and we may use the objects of Section~\ref{ssec:arcs} for $(X',g',B',q)$. As $\Pi$ is a finite surjective morphism, $\Delta\neq\PP^m$ and $X'$ is irreducible, the regular function $\delta\circ\Pi$ is a nonzero element of $\cOa_{X',q}$. By Proposition~\ref{prop:reduction}(3), applied to $h=\delta\circ\Pi$, the assumption of Proposition~\ref{prop:local-aut} holds, so $g'$ is a local biholomorphism at $q$. By Proposition~\ref{thm:claimS}, every eigenvalue of $g'^*_q$ on $V_q$ has modulus greater than one.

Since $S$ is complete Noetherian and $a'\neq0$, Krull's intersection theorem gives $j\geq1$ with $a'\in\fm_S^j\setminus\fm_S^{j+1}$. Proposition~\ref{prop:reduction}(2) shows that the nonzero class of $a'$ in $\fm_S^j/\fm_S^{j+1}$ is a generalized eigenvector with eigenvalue $\lambda^\kappa$. The quotient map $\pr$ gives an equivariant surjection $\widehat\fm_q^j/\widehat\fm_q^{j+1}\to\fm_S^j/\fm_S^{j+1}$, so $\lambda^\kappa$ is an eigenvalue of $\gr^jg'^*$ on $\widehat\fm_q^j/\widehat\fm_q^{j+1}$. By \eqref{eq:graded-eigenvalues}, $\lambda^\kappa$ is a product of $j$ eigenvalues of $g'^*_q$ on $V_q$, with repetitions allowed. Each factor has modulus greater than one. Therefore $|\lambda|^\kappa>1$ and $|\lambda|>1$.
\end{proof}

\begin{remark}\label{rem:reduction-dimension}
In the proof of Proposition~\ref{prop:reduction}, $\dim Z=\dim R'$, the map $\hat\omega$ is injective and $\hcO_{X,q}=\widetilde R$: the completed local ring of the normalization of $Z$ at $q$ is the normalization of the formal branch $R'$. Indeed, replace $g_0$ by $g_0^\kappa$, and let $\nu_q^*\colon R\to\hcO_{X,q}$ be induced by $\nu$. By Propositions~\ref{prop:reduction}(3), applied to $h=\delta\circ\Pi$, and~\ref{prop:local-aut}, $\fn(\hcO_{X,q},g^*)=0$, as in the proof of Theorem~C; so $\nu_q^*(\fn)=0$ by Lemma~\ref{lem:fitting}(1), and $\nu_q^*$ induces a homomorphism $R_{\mathrm{inv}}\to\hcO_{X,q}$, whose kernel $\fP''$ is prime, as $\hcO_{X,q}$ is a domain. By Step~4 of that proof, $\hat\omega\circ\nu_q^*$ is the map $R\to R'\subset\widetilde R$; hence $\fP''\subset\fP$, and $\fP''=\fP$ as $\fP$ is a minimal prime. Thus $R'\subset\hcO_{X,q}$, a finite extension since $\nu$ is finite, and $\dim R'=\dim\hcO_{X,q}=\dim Z$. As $S$ has the same dimension, $\ker\hat\omega=0$, and the normal domain $\hcO_{X,q}$, which lies between $R'$ and $\widetilde R$, is equal to $\widetilde R$. For $X_0=\PP^N$, Lemma~\ref{lem:fitting}(1)--(3) and Theorem~\ref{thm:core} show moreover that $R_{\mathrm{inv}}$ is the completed local ring of the smooth germ $\mathcal B_W$, and that $Z=W$.
\end{remark}

\section{Consequences}\label{sec:consequences}

\begin{corollary}\label{cor:conjugates}
Let $f$ be a PCF endomorphism of $\PP^N$ of degree $d\geq2$, and $\lambda$ a multiplier of $f$ at a periodic point. Then $\lambda$ is algebraic over $\QQ$. If $\lambda\neq0$, every conjugate of $\lambda$ over $\QQ$ has modulus greater than one.
\end{corollary}

\begin{proof}
Let $p$ be a periodic point of period $a$ such that $\lambda$ is an eigenvalue of $D_pf^a$, and $\sigma$ a field automorphism of $\CC$. Applying $\sigma$ to the coefficients of a homogeneous lift of $f$ gives an endomorphism $f^\sigma$ of degree $d$ with $\sigma\circ f=f^\sigma\circ\sigma$ on $\PP^N(\CC)$; its critical set is $C(f^\sigma)=\sigma(C(f))$. Since $f$ is PCF, its post-critical set is a finite union of forward images of $C(f)$. Applying $\sigma$ to this finite algebraic union gives $\PC(f^\sigma)=\sigma(\PC(f))$, so $f^\sigma$ is PCF. In an affine chart defined over $\QQ$, $\sigma(\lambda)$ is an eigenvalue of $D_{\sigma(p)}(f^\sigma)^a$.

Suppose that $\lambda$ is transcendental. As $\Qbar$ is countable, there is a transcendental number $t$ with $0<|t|<1$, and the isomorphism $\QQ(\lambda)\to\QQ(t)$ sending $\lambda$ to $t$ extends to an automorphism $\sigma$ of $\CC$. Then $t=\sigma(\lambda)$ is a multiplier of the PCF endomorphism $f^\sigma$ with $0<|t|<1$, which contradicts the Main Theorem. Thus $\lambda$ is algebraic. Let $\lambda\neq0$, and let $\lambda'$ be a conjugate of $\lambda$. The embedding $\QQ(\lambda)\to\CC$ sending $\lambda$ to $\lambda'$ extends to an automorphism $\sigma$ of $\CC$, and $\lambda'=\sigma(\lambda)$ is a nonzero multiplier of $f^\sigma$ at $\sigma(p)$; so $|\lambda'|>1$ by the Main Theorem.
\end{proof}

\begin{example}\label{ex:nonintegral}
For $d\geq2$, the polynomial
\[
p_d(z):=\frac{dz-z^d}{d-1}
\]
is critically fixed: its finite critical points satisfy $z^{d-1}=1$ and are fixed, and infinity is fixed. Thus $p_d$ is PCF, while $0$ is fixed with multiplier
\[
p_d'(0)=\frac{d}{d-1}.
\]
For $d\geq3$, this rational number is not an algebraic integer. Hence the algebraicity in Corollary~\ref{cor:conjugates} cannot be strengthened to integrality, even on $\PP^1$. These multipliers tend to $1$ as $d\to\infty$, so there is no uniform gap above $1$ independent of the degree.
\end{example}

Integrality does hold if $f$ is defined over a number field and has potentially good reduction \cite{Ing16} at every finite place. After a finite extension of each completion and a change of coordinates giving good reduction, a periodic point extends to a section over the valuation ring. The characteristic polynomial of the differential along that section has integral coefficients, so every multiplier is integral at every finite place.

Endomorphisms of $\PP^N$ of degree $d$ are the points of the complement $\mathcal E_d$ of the resultant hypersurface in the projective space of $(N+1)$-tuples of homogeneous polynomials of degree $d$ in $N+1$ variables; $\mathcal E_d$ is defined over $\QQ$. A \emph{holomorphic family} of endomorphisms of degree $d$ over a complex manifold $S$ is a holomorphic map $S\to\mathcal E_d$, $s\mapsto f_s$. Two endomorphisms $h$, $h'$ of $\PP^N$ are \emph{isospectral} if, for every $n\geq1$, the multisets $\{\det(T-D_zh^n):z\in\Fix(h^n)\}$ and $\{\det(T-D_zh'^n):z\in\Fix(h'^n)\}$ coincide. For multiplier maps and stability in algebraic families, see \cite[Sections~2.1.3--2.2]{GTV26}.

\begin{proposition}\label{prop:families}
Let $f$ be a PCF endomorphism of $\PP^N$ of degree $d\geq2$.
\begin{enumerate}[label=\textup{(\arabic*)}]
\item For every $n\geq1$, the scheme-theoretic intersection of the graph of $f^n$ with the diagonal of $\PP^N\times\PP^N$ is reduced. Hence $f^n$ has exactly $\sum_{i=0}^Nd^{ni}$ fixed points.
\item Let $(f_s)_{s\in S}$ be a holomorphic family of PCF endomorphisms of $\PP^N$ of degree $d$ over a connected complex manifold $S$. For every $n\geq1$, the set $\mathcal F_n:=\{(s,z)\in S\times\PP^N:f_s^n(z)=z\}$ is a complex submanifold of $S\times\PP^N$, and the projection $\mathcal F_n\to S$ is a covering map of degree $\sum_{i=0}^Nd^{ni}$. The maps $f_s$, $s\in S$, are pairwise isospectral.
\item The map $f$ is isospectral to a PCF endomorphism of $\PP^N$ of degree $d$ defined over $\Qbar$.
\end{enumerate}
\end{proposition}

\begin{proof}
(1) Let $z\in\Fix(f^n)$, of period $\ell$ dividing $n$. The eigenvalues of $D_zf^n=(D_zf^\ell)^{n/\ell}$ are powers of multipliers, hence $0$ or of modulus greater than one by the Main Theorem, so $D_zf^n-\Id$ is invertible. Thus the graph of $f^n$ and the diagonal meet transversally at $(z,z)$, and the number of fixed points is the intersection number of the graph with the diagonal, which is $\sum_{i=0}^Nd^{ni}$ by a computation in the cohomology ring of $\PP^N\times\PP^N$.

(2) The set $\mathcal F_n$ is closed in $S\times\PP^N$, so the projection $\mathcal F_n\to S$ is proper. By the proof of~(1) and the implicit function theorem, $\mathcal F_n$ is, near each of its points $(s_0,z_0)$, the graph of a holomorphic map $s\mapsto z(s)$. Thus $\mathcal F_n$ is a submanifold and the projection is a proper local homeomorphism, hence a covering map, of degree given by~(1). Along such a local section, the coefficients of $\det(T-D_{z(s)}f_s^n)$ are holomorphic in $s$, and algebraic by Corollary~\ref{cor:conjugates}, since the eigenvalues of $D_{z(s)}f_s^n$ are powers of multipliers of $f_s$. A continuous function on a connected set with values in the countable set $\Qbar$ is constant. Hence the multiset of characteristic polynomials at the points of $\Fix(f_s^n)$ is locally constant in $s$, hence constant on $S$.

(3) Let $\mathcal C\subset\mathcal E_d\times\PP^N$ be the critical divisor of the universal endomorphism $\mathcal U\colon(h,z)\mapsto(h,h(z))$, both defined over $\QQ$. The map $\mathcal U$ is proper with finite fibres, hence finite, so $\mathcal C_j:=\mathcal U^j(\mathcal C)$ is a closed algebraic subset defined over $\QQ$, whose fibre over $h$ is $h^j(C(h))$. For integers $u>v\geq1$, set
\[
\mathcal H_{u,v}:=\{h\in\mathcal E_d:h^u(C(h))\subset h^v(C(h))\}.
\]
Its complement is the projection of $\mathcal C_u\setminus\mathcal C_v$ to $\mathcal E_d$, so $\mathcal H_{u,v}$ is constructible and defined over $\QQ$, by Chevalley's theorem over $\QQ$ \cite[Chapter~I, \S8, Corollary~2]{Mum99}. If $h\in\mathcal H_{u,v}$, then $h^j(C(h))\subset\bigcup_{v\leq i<u}h^i(C(h))$ for all $j\geq v$, by induction on $j$, so $h$ is PCF. Conversely, if $h$ is PCF, each $h^j(C(h))$, $j\geq1$, is a union of irreducible components of $\PC(h)$, so $h^u(C(h))=h^v(C(h))$ for some $u>v\geq1$. Choose $u,v$ with $f\in\mathcal H_{u,v}$. The constructible set $\mathcal H_{u,v}$ is a finite disjoint union of smooth irreducible locally closed subvarieties defined over $\Qbar$; let $S$ be the one containing $f$. Then $S(\CC)$ is a complex manifold, connected by \cite[Chapter~7, \S2, Theorem~7.1]{Sha13II}, and $S$ has a point $f_0$ defined over $\Qbar$. All members of the holomorphic family $S(\CC)\hookrightarrow\mathcal E_d$ are PCF, so $f$ and $f_0$ are isospectral by~(2).
\end{proof}

\begin{corollary}\label{cor:tower}
Let $g$ be a polarized endomorphism of a normal projective variety $X$ which admits a branch boundary, for instance a PCF endomorphism of $\PP^N$. Let $Z\subset X$ be an irreducible closed subvariety with $g^k(Z)=Z$ for some $k\geq1$, $\nu\colon X_Z\to Z$ its normalization, and $g_Z$ the map induced by $g^k$. Then every periodic point of $g_Z$ is super-saddle. The same holds for the maps induced on the normalizations of the periodic subvarieties of $X_Z$, and so on recursively.
\end{corollary}

\begin{proof}
If $\dim Z=0$ there is nothing to prove. Otherwise $g^k$ is polarized, with the same branch boundary, so $g_Z$ is a polarized endomorphism of $X_Z$ which admits a branch boundary, by Proposition~\ref{prop:normalization-boundary} applied to $g^k$, and Theorem~C applies to $g_Z$. The last assertion follows by induction.
\end{proof}

For chains of periodic prime divisors, that the induced maps admit branch boundaries also follows from the theorem of Dupuy, Levy, Tucker and Wang \cite[Theorem~1.1 and Corollary~8.2]{DLTW26} and Proposition~\ref{prop:divisorial}(2); the corollary adds the super-saddle property and the subvarieties of higher codimension.

\begin{corollary}\label{cor:nonnormal}
Let $g_0$ be a polarized endomorphism of a projective variety $X_0$, $\nu\colon X\to X_0$ its normalization, and $g$ the induced endomorphism of $X$. If $g$ admits a branch boundary, then every periodic point of $g_0$ is super-saddle on its Zariski tangent space.
\end{corollary}

\begin{proof}
By Lemma~\ref{lem:pcf-transport}(d), applied to $\nu$, the map $g$ is polarized, so its periodic points are super-saddle by Theorem~C. Corollary~\ref{cor:descent}(a), applied to $\nu$, $g_0$ and $g$, gives the super-saddle property at every periodic point of $g_0$.
\end{proof}

\section{Examples}\label{sec:examples}

The examples test the hypotheses of the theorems and show singular post-critical sets and singular ambient points at periodic points.

\subsection{Rotation discs}\label{ssec:ex-discs}

\begin{example}\label{ex:sharp}
The following examples show that polarization, the algebraicity of $B$ and the condition $|\lambda|=1$ are needed in Corollary~\ref{cor:no-rotation-disc}. Let $\lambda=e^{2\pi i\xi}$ with $\xi\in\RR\setminus\QQ$.
\begin{enumerate}[label=\textup{(\arabic*)},leftmargin=*]
\item The rotation $g(z)=\lambda z$ of $\PP^1$ has $M=1$, the branch boundary $\{0,\infty\}$ and the rotation disc $\varphi(z)=z$.
\item For almost every $\xi$, $f(z)=z^2+\lambda z$ has a Siegel disc $U$ at $0$, linearized by $\varphi\colon\DD\to U$ \cite[Theorem~11.14 and Lemma~11.1]{Mil06}. The set $B:=\{\infty\}\cup\overline{\{f^n(\omega):n\geq1\}}$, $\omega:=-\lambda/2$, is closed and invariant and contains every $\Br(f^n)$. As $\partial U$ is infinite and lies in the closure of the orbit of $\omega$ \cite[Theorem~11.17]{Mil06}, this orbit avoids $U$ (an orbit entering $U$ stays on an invariant circle, so its closure would meet $\partial U$ in finitely many points); so $B\cap U=\emptyset$, and $B\supset\partial U$ is not algebraic.
\item For $g(z)=z^2$ and $B=\{0,\infty\}$, $\varphi=\exp$ satisfies $g\circ\varphi(z)=\varphi(2z)$ and $\varphi(\CC)\cap B=\emptyset$.
\item The map $g(z,w)=(z^2,\lambda w)$ of $\PP^1\times\PP^1$ is not polarized; the union $B$ of the lines $\{z=0\}$, $\{z=\infty\}$, $\{w=0\}$, $\{w=\infty\}$ is invariant and contains every $\Br(g^n)$, and $\varphi(z)=(1,z)$ is a rotation disc with $\varphi(\DD^*)\cap B=\emptyset$.
\end{enumerate}
\end{example}

\begin{example}\label{ex:discs}
Rotation discs occur for polarized maps without a branch boundary, with boundary in the post-critical set. With $\lambda$, $f$, $\omega$ as in Example~\ref{ex:sharp}(2), let $\psi\colon\DD\to\Delta$ be the linearization of $f$ on the Fatou component $\Delta$ of $0$; then $\partial\Delta\subset\PC(f)$ and $\Delta\cap\PC(f)=\emptyset$, as in Example~\ref{ex:sharp}(2).
\begin{enumerate}[label=\textup{(\alph*)},leftmargin=*]
\item If $\xi$ is a Bryuno number \cite[Theorem~11.10]{Mil06} and $c^2=\lambda$, the linearization $\psi_c$ of $f_c(z):=z(z-c)^2$ at its fixed point $0=f_c(c)$ is a rotation disc with $\psi_c(0)\in\PC(f_c)$ and $\psi_c(\DD^*)\cap\PC(f_c)=\emptyset$.
The finite critical points are $c$ and $c/3$, with $f_c(c)=0$. The boundary of the Siegel disc is contained in the closure of the critical orbits \cite[Theorem~11.17]{Mil06}. If the orbit of $c/3$ entered this disc, its closure would lie in finitely many points and a compact invariant circle in the disc. Together with the orbit $c\mapsto0\mapsto0$ and the fixed point at infinity, these orbits could not accumulate on its boundary. Thus the orbit of $c/3$ stays outside the disc, and its closure avoids it. The intersection of the post-critical set with the Siegel disc is therefore exactly $\{0\}$.

\item The endomorphism $F[s:p:t]:=[s^2-2pt+\lambda st:p^2+\lambda sp+\lambda^2pt:t^2]$ of $\PP^2$, with $F\circ\eta=\eta\circ(f,f)$ for $\eta(x,y):=(x+y,xy)$, has the rotation disc $\varphi:=(\psi,0)$, centred on the invariant conic $\{s^2=4p\}\subset\PC(F)$, with $\varphi(\DD^*)\cap\PC(F)=\emptyset$ and boundary in $\PC(F)$, inside a rotation domain \cite[Definition~4.1]{Ued98}.
\item For $g[z:w:t]:=[z^2+\lambda zt:w^2+bzt:t^2]$ with
\[
0<|b|\leq\frac{3}{16\sup_\DD|\psi|},
\]
a contraction argument in $H^\infty(\DD)$ gives a rotation disc $\varphi=(\psi,s)$ centred at the critical point $(0,0)$. Indeed, the operator
\[
T_b(v)(\zeta):=v(\lambda^{-1}\zeta)^2+b\psi(\lambda^{-1}\zeta)
\]
maps the closed ball of radius $1/4$ into itself and is $\frac12$-Lipschitz there. Its fixed point $s$ satisfies $s(0)=0$ and
\[
s(\lambda\zeta)=s(\zeta)^2+b\psi(\zeta),
\]
so $s$ is nonconstant. If $(z_n,w_n)$ is the orbit of the critical point $(\psi(\zeta),0)$ and $s_n:=s(\lambda^n\zeta)$, then
\[
z_n=\psi(\lambda^n\zeta),\qquad
w_{n+1}-s_{n+1}=(w_n-s_n)(w_n+s_n).
\]
Since $|w_n+s_n|\leq|w_n-s_n|+\frac12$, induction gives $|w_n-s_n|\leq\frac14(3/4)^n$. Recurrence of the irrational rotation now gives $\varphi(\DD)\subset\PC(g)$.

For all but countably many allowed values of $b$, this disc lies in no $\Br(g^n)$. To prove this, fix $n\geq1$ and $a\in f^{-n}(0)\setminus\{0\}$, where $f(z)=z^2+\lambda z$. Set $z_j=f^j(a)$, $w_0=0$ and $w_{j+1}=w_j^2+bz_j$. Then $w_n$ is a nonzero polynomial in $b$: if $r<n$ is the last index with $z_r\neq0$, then
\[
w_{r+1}=bz_r+O(b^2),\qquad
w_n=w_{r+1}^{\,2^{n-r-1}}.
\]
Exclude the zeros of these polynomials for all $n$ and $a$, a countable set. For every remaining $b\neq0$, the only point of $\{w=0\}$ mapped to $(0,0)$ by $g^n$ is $(0,0)$. In the coordinate $z=\psi(\zeta)$, the germ at $(0,0)$ of $g^n(\{w=0\})$ is the graph of $v_n:=T_b^n(0)$. Both $s$ and each $v_j$, $j\geq1$, have derivative $b\psi'(0)/\lambda\neq0$ at $0$. Writing $R(\zeta)=\lambda^{-1}\zeta$, we have
\[
s-v_n=((s-v_{n-1})\circ R)((s+v_{n-1})\circ R).
\]
For $n=1$ the second factor is $s\circ R\not\equiv0$; for $n\geq2$ its derivative at $0$ is nonzero. Induction gives $s-v_n\not\equiv0$ for every $n$. The other critical components are $\{z=-\lambda/2\}$ and the line at infinity; none of their images contains this nonvertical affine disc. Hence no $\Br(g^n)$ contains it. The punctured-avoidance hypothesis of Theorem~A is therefore a restriction even for a rotation disc contained in no individual iterate branch locus.

\end{enumerate}
\end{example}

\subsection{Quotients of tori and of products of the line}\label{ssec:ex-quotients}

An endomorphism $f$ of $\PP^n$ of degree $d\geq2$ is a \emph{Latt\`es map} if $\Theta\circ A=f\circ\Theta$ for a complex torus $T=\CC^n/\Lambda$, an affine map $A(z)=\alpha z+c$ of $T$, with $\alpha(\Lambda)\subset\Lambda$, and a finite surjective holomorphic map $\Theta\colon T\to\PP^n$ (Dupont \cite{Dup03} and Berteloot--Loeb \cite{BeLo01} assume that $\Theta$ is a quotient map by a finite group).

\begin{proposition}\label{prop:lattes}
Let $f$ be a Latt\`es map of $\PP^n$ of degree $d$. Every multiplier of $f$ at a periodic point of period $a$ has modulus $d^{ak/2}$ for some integer $k\geq1$.
\end{proposition}

\begin{proof}
The Hermitian form $\langle\cdot,\cdot\rangle$ of the ample line bundle $\mathcal L:=\Theta^*\cO(1)$ is positive definite, and that of $A^*\mathcal L\simeq\mathcal L^{\otimes d}$ is $\langle\alpha\cdot,\alpha\cdot\rangle=d\langle\cdot,\cdot\rangle$ \cite[Proposition~4.5.2, Lemmas~2.3.2 and~2.3.4]{BL04}; so the eigenvalues of $\alpha$ have modulus $d^{1/2}$. Let $\lambda$ be a multiplier of $f$ at a point $p$ of period $a$. At a periodic point of $A$ of period $b$ in $\Theta^{-1}(p)$, Corollary~\ref{cor:transfer-multipliers}(2) shows that $\lambda^{b/a}$ is a monomial of degree $k\geq1$ in the eigenvalues of $\alpha^b$.
\end{proof}

For $n\geq1$, the quotient map $\eta\colon(\PP^1)^n\to\PP^n$ by $\mathfrak S_n$ sends $([a_i:b_i])_i$ to the coefficients $[e_0:\dots:e_n]$ of $\prod_i(b_iX+a_iY)$, and a rational map $h$ of degree $d\geq2$ induces an endomorphism $F_h$ of $\PP^n$ of degree $d$ with $\eta\circ h^{\times n}=F_h\circ\eta$ (for $n=2$, Ueda \cite[Section~4.1]{Ued93}; see \cite{GHK23}).

\begin{proposition}\label{prop:symmetric}
$F_h$ is PCF if and only if $h$ is \cite{GHK23}. If $p$ is a periodic point of $F_h$ of period $a$ and $\lambda$ a multiplier at $p$, then for every periodic point $y=(z_1,\dots,z_n)\in\eta^{-1}(p)$ of $h^{\times n}$, of period $b$, one has $a\mid b$ and
\[
\lambda^{b/a}=\prod_{i=1}^n\bigl((h^b)'(z_i)\bigr)^{k_i}\qquad\text{with integers }k_i\geq0,\ 1\leq\textstyle\sum_ik_i\leq n .
\]
In particular, if $h$ is PCF, then $F_h$ is PCF and all periodic points of $F_h$ are super-saddle.
\end{proposition}

\begin{proof}
As $\PP^n$ is smooth, Corollary~\ref{cor:transfer-multipliers}(2) applies to every multiplier and every periodic point $y$ of the fibre, with a monomial of degree between $1$ and $N(\eta,p)$, and $D_y(h^{\times n})^b$ is diagonal with entries $(h^b)'(z_i)$. Moreover $N(\eta,p)\leq n$: near $p$, $\eta$ is a product of quotient maps $\CC^m\to\CC^m/\mathfrak S_m$ with $m\leq n$, whose coordinates are elementary symmetric polynomials of degrees at most $m$, algebraically independent. If $h$ is PCF, each $(h^b)'(z_i)$ is zero or has modulus greater than one by the one-dimensional result \cite[Corollary~14.5]{Mil06}; the displayed formula gives the super-saddle property for the periodic points of $F_h$.
\end{proof}

Let $\Phi$ be a reduced root system of rank $n$, with weight lattice $P$ and Weyl group $W$, and $\pi\colon T:=\operatorname{Hom}(P,\CC^*)\to T/W\simeq\CC^n$ the quotient map given by the fundamental characters \cite[Chapter~VI, \S3, no.~4, Theorem~1]{Bou68}; for $m\geq2$, $\pi(t^m)=P_m(\pi(t))$ for a polynomial map $P_m$. For the types $A_n$ and $C_n$, $P_m$ is, up to affine conjugacy \cite{FH91}, the Chebyshev map of Hoffman--Withers \cite{HW88} and Veselov \cite{Ves87,Ves91}, defined with the elementary symmetric polynomials $e=(e_1,\dots,e_n)$ and $T_m(s+s^{-1})=s^m+s^{-m}$ by
\[
P^{A_n}_m(e(t))=e(t_1^m,\dots,t_{n+1}^m),\qquad P^{C_n}_m(e(x))=e(T_m(x_1),\dots,T_m(x_n)),
\]
for $t\in\mathbb T:=\{t\in(\CC^*)^{n+1}:t_1\cdots t_{n+1}=1\}$ and $x\in\CC^n$.

\begin{proposition}\label{prop:chebyshev}
The maps $P^{A_n}_m$ and $P^{C_n}_m$ extend to endomorphisms of $\PP^n$ of algebraic degree $m$.
\end{proposition}

\begin{proof}
As $T_m^{-1}(\infty)=\{\infty\}$, $F_{T_m}$ maps $\{e_0\neq0\}$ into itself, where it is $P^{C_n}_m$. For $h(z)=z^m$, the endomorphism $S_m:=F_h=[e_0^m:\ast:\dots:\ast:e_{n+1}^m]$ of $\PP^{n+1}$ preserves the hyperplane $\Pi:=\{e_0=e_{n+1}\}\simeq\PP^n$, and $\Pi\cap\{e_0=1\}=e(\mathbb T)$, so $S_m|_\Pi$ extends $P^{A_n}_m$.
\end{proof}

For $n=2$, the description of type $C_2=B_2$ by symmetric products is due to Nakane \cite{Nak08} (see \cite[p.~995]{Uch09}). For $A_2$ and $A_3$, post-critical finiteness is stated in \cite[p.~999]{Uch09} and \cite{Uch17}, with references to Dinh and to Dinh and Sibony; Proposition~\ref{prop:foldings} below proves it for every $n$.

\begin{proposition}\label{prop:foldings}
Suppose that $P_m$ extends to an endomorphism $F$ of $\PP^n\supset\CC^n$. Then $F$ is PCF, its periodic cycles are super-saddle, and every multiplier of $F$ at a periodic point of period $a$ in $\CC^n$ has modulus $m^{ak}$ for some integer $k\geq1$.
\end{proposition}

\begin{proof}
The hyperplane $H_\infty:=\PP^n\setminus\CC^n$ satisfies $F^{-1}(H_\infty)=H_\infty$, and $F$ has degree $m$, as $P_m$ has the topological degree $m^n$ of $t\mapsto t^m$. Let $\lambda$ be a multiplier of $F$ at a point $p\in\CC^n$ of period $a$. At a periodic point of $t\mapsto t^m$ of period $b$ in $\pi^{-1}(p)$, the differential of $t\mapsto t^{m^b}$ is $m^b\Id$, and Corollary~\ref{cor:transfer-multipliers}(2) gives $|\lambda|^{b/a}=m^{bk}$, $k\geq1$. As $t\mapsto t^m$ is \'etale and commutes with $W$, Lemma~\ref{lem:pcf-transport}(c) and the paragraph following it show that $P_m$ is divisorially PCF. As $F^{-1}(H_\infty)=H_\infty$, the set $\PC(F)$ is contained in the union of $H_\infty$ and of the closure of $\PCdiv(P_m)$, a finite union of hypersurfaces; so $F$ is PCF. The super-saddle property is the Main Theorem.
\end{proof}

\begin{remark}\label{rem:foldings}
If $\Phi=\Phi_1\times\Phi_2$, then $P_m$ is the product of the maps of $\Phi_1$ and $\Phi_2$, and it extends to $\PP^n$ if both do. By Proposition~\ref{prop:chebyshev}, Proposition~\ref{prop:foldings} thus applies to $A_n$ and $C_n$ for every $n$, and to $A_1\times A_1\times A_1$, $A_1\times A_2$ and $A_1\times C_2$. For $m=2$, in the fundamental-character coordinates, $P_m$ has degree $3$ and topological degree $4$ for $G_2$, $8$ for $B_3$, so it extends neither to $\PP^2$ nor to $\PP^3$.
\end{remark}

\begin{example}\label{ex:cone-construction}
Let $q=[q_0:\dots:q_{n-1}]$ be a PCF endomorphism of $\PP^{n-1}$, $n\geq2$, of degree $d\geq2$. Then $\hat q:=[q_0:\dots:q_{n-1}:w^d]$ is a PCF endomorphism of $\PP^n$: $\PC(\hat q)$ is the union of $\{w=0\}$ and of the cone over $\PC(q)$ with vertex $o:=[0:\dots:0:1]$. The vertex is fixed with $D_o\hat q=0$. At a point $[z:1]\neq o$ of period $a$, the multipliers are $d^a$ and the eigenvalues of $D_{[z]}q^a$, which are powers of the multipliers of $q$ at $[z]$, by Euler's formula and \cite[Remark~3.3]{Le22}.
\end{example}

\subsection{Post-critical sets with singular points}\label{ssec:ex-singular-pc}

\begin{figure}[ht]
\centering
\begin{tikzpicture}[x=0.47cm,y=0.47cm,>={Stealth[length=2mm]},font=\small]
  \draw[dashed,gray] (0,0) circle [radius=4.2];
  \node[gray,anchor=north west,align=left,inner sep=1pt] at (2.75,-3.25) {circle at infinity\\ (schematic)};
  \draw[gray!60,->] (-3.75,0) -- (3.75,0) node[below left,gray,inner sep=3pt] {$\operatorname{Re}x$};
  \draw[gray!60,->] (0,-3.75) -- (0,3.75) node[below right,gray,inner sep=2pt] {$\operatorname{Im}x$};
  \draw[thick,fill=gray!18]
    plot[domain=0:360,samples=241,smooth] ({2*cos(\x)+cos(2*\x)},{2*sin(\x)-sin(2*\x)}) -- cycle;
  \node at (1.75,1.25) {$Q_4$};
  \fill (3,0) circle [radius=2.2pt];
  \draw[thin] (3.1,0.1) -- (4.6,1.2);
  \node[anchor=south west,align=left,inner sep=1pt] at (4.6,1.2) {fixed cusp $3$\\ multipliers $4$, $8$};
  \fill (-1.5,2.598) circle [radius=2.2pt];
  \fill (-1.5,-2.598) circle [radius=2.2pt];
  \node[anchor=south,inner sep=3pt] at (-1.5,2.6) {$3\omega$};
  \node[anchor=north,inner sep=3pt] at (-1.5,-2.6) {$3\omega^2$};
  \draw[<->,thick] (-1.8,2.3) .. controls (-3.4,1.2) and (-3.4,-1.2) .. (-1.8,-2.3);
  \draw[thin] (-3.0,0.25) -- (-4.6,1.2);
  \node[anchor=south east,align=right,inner sep=1pt] at (-4.6,1.2) {$2$-cycle\\ $Df^2$: $16$, $64$};
\end{tikzpicture}
\caption{The invariant real slice $\{y=\bar x\}$ of Example~\ref{ex:chebyshev-A2}, on which $f$ acts by $x\mapsto x^2-2\bar x$. It meets $Q_4$ in the deltoid $x=2e^{i\theta}+e^{-2i\theta}$, which bounds the image of the compact torus (grey), and its closure meets the line $\{w=0\}$ in a circle (dashed, drawn schematically), which does not contain the points $[1:0:0]$, $[0:1:0]$ of $\Sing\PC(f)$.}
\label{fig:chebyshev-A2}
\end{figure}
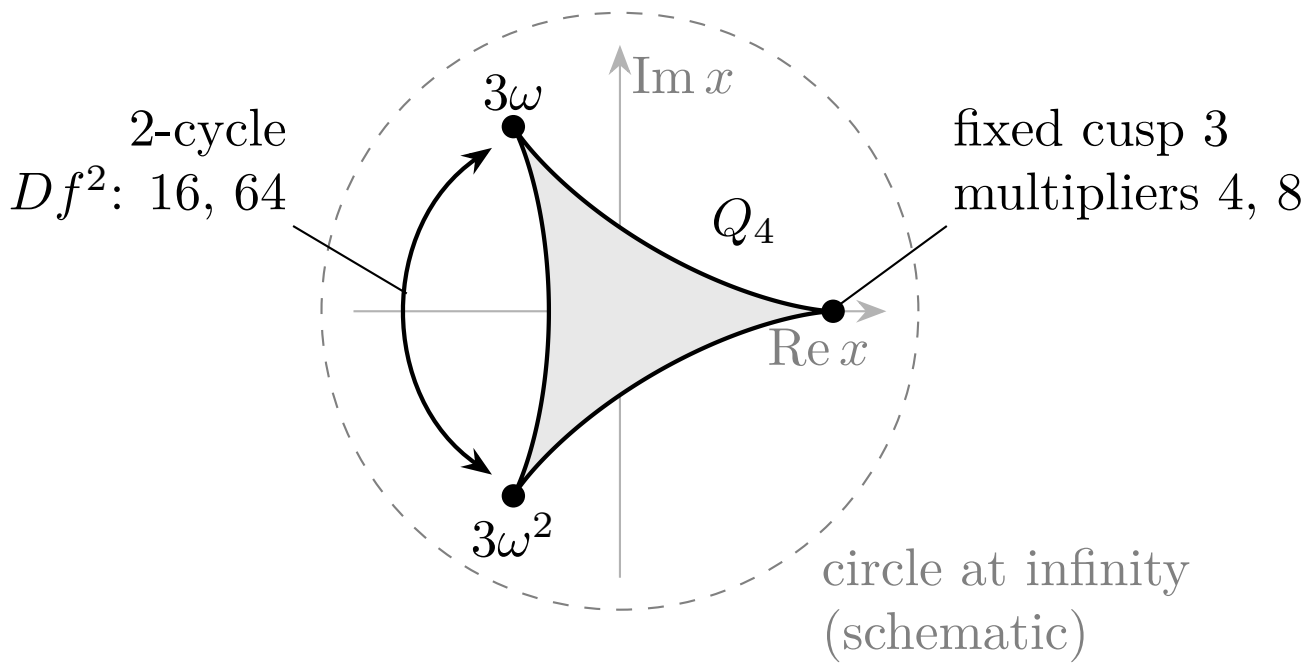

\begin{example}\label{ex:chebyshev-A2}
Let $f:=[x^2-2yw:y^2-2xw:w^2]$, the extension of $P^{A_2}_2$ in Proposition~\ref{prop:chebyshev}, with $(x,y)=e(t)$ on $\{w=1\}$. Its critical set $\{w=0\}\cup\{xy=w^2\}$ is mapped onto $\{w=0\}\cup Q_4$, where
\[
Q_4:=\{x^2y^2-4x^3w-4y^3w+18xyw^2-27w^4=0\}
\]
is the branch locus of $t\mapsto e(t)$, and both curves are invariant: $\PC(f)=\{w=0\}\cup Q_4$. The normalization of $Q_4$ is $\nu[s:u]:=[2s^3u+u^4:s^4+2su^3:s^2u^2]$, with $\nu(s)=e(s,s,s^{-2})$ for $u=1$, and $f\circ\nu=\nu\circ h$ for $h[s:u]:=[s^2:u^2]$. The set $\Sing\PC(f)$ consists of the cusps $\nu(1)=[3:3:1]$, $\nu(\omega)$, $\nu(\omega^2)$ of $Q_4$, where $\omega:=e^{2\pi i/3}$, and of the points $[1:0:0]$, $[0:1:0]$ of $Q_4\cap\{w=0\}$. The cusp $[3:3:1]$ is fixed, with multipliers $4$ and $8$; the cusps $\nu(\omega)$, $\nu(\omega^2)$ form a $2$-cycle, along which $Df^2$ has the eigenvalues $16$ and $64$; the points $[1:0:0]$ and $[0:1:0]$ are fixed, with $Df$ nilpotent. At a cusp $p$, $T_pQ_4=T_p\PP^2$, and these are the eigenvalues $\alpha^2$, $\alpha^3$ of Example~\ref{ex:transfer}(a), with $\alpha=h'(1)=2$, respectively $\alpha=(h^2)'(\omega)=4$ (Figure~\ref{fig:chebyshev-A2}).
\end{example}

\begin{table}[ht]
\centering
\footnotesize
\setlength{\tabcolsep}{3pt}
\renewcommand{\arraystretch}{1.15}
\begin{tabular}{@{}>{\raggedright\arraybackslash}p{0.45cm}>{\raggedright\arraybackslash}p{6.2cm}>{\raggedright\arraybackslash}p{4.3cm}>{\raggedright\arraybackslash}p{3.8cm}@{}}
\toprule
 & $f$ & $\PC(f)$; $\Sing\PC(f)$ & periodic points in $\Sing\PC(f)$\\
\midrule
(a) & $[x^2:y^2:z^2]$ & $3$ lines; $3$ points & all fixed: $\{0,0\}$\\
(b) & $[(-x+y+z)^2:(x-y+z)^2:(x+y-z)^2]$ \cite[Example~5.10]{Ued98} & $6$ lines; $7$ points & $[1:1:1]$ fixed: $\{-4,-4\}$\\
(c) & $[x^3-y^2z:-y^3:-z^3]$ \cite[Remark~2.12]{Rong08} & $2$ lines, a line and a cuspidal cubic exchanged; $3$ points & all fixed: $\{0,0\}$\\
\midrule
(d) & $F_h$, $h(z)=z^2-2$, $n=3$ (type $C_3$) & $3$ planes, a quartic $D$; $10$ curves, one the cuspidal edge of $D$ & a fixed point: $\{4,16,64\}$; a $2$-cycle: $\{16,-4,-64\}$\\
(e) & $[x_0^2:x_1^2-2x_0x_2:2x_0^2+x_2^2-2x_1x_3:x_3^2-2x_0x_2]$ (type $A_3$ \cite{Uch17}) & a plane, a sextic $D$; curves, one a cuspidal edge of $D$ & a fixed point: $\{4,8,16\}$; a $2$-cycle: $\{4,16,64\}$\\
(f) & $[(x_0-2x_1)^2:(x_0-2x_2)^2:(x_0-2x_3)^2:x_0^2]$ (Forn\ae ss--Sibony, see \cite[Section~2]{Jon98}) & $10$ planes, none invariant; $25$ lines & $[1:1:1:1]$ fixed: $\{-4,\pm4i\}$; a $2$-cycle: $\{-2,-4\pm4\sqrt3\}$\\
(g) & $[x_0^3-x_1^2x_2:-x_1^3:-x_2^3:x_3^3]$ & $3$ planes, a plane and a cone $V$ exchanged; $7$ curves & the apex, fixed: $\{0,0,0\}$; $2$-cycles: $\{9,0,0\}$\\
\bottomrule
\end{tabular}
\caption{PCF endomorphisms of $\PP^2$ and $\PP^3$ with singular post-critical sets, and the spectra of $Df$, respectively $Df^2$, at periodic points of period at most two in $\Sing\PC(f)$: all of them in rows (a)--(c), a selection in rows (d)--(g), with explicit representatives and characteristic polynomials in Example~\ref{ex:P2-isolated}.}
\label{tab:pc}
\end{table}

\begin{example}\label{ex:P2-isolated}\label{ex:P3-curves}\label{ex:cones}
Table~\ref{tab:pc} lists further examples. In (a)--(c) the set $\Sing\PC(f)$ is finite, in (d)--(g) it is a union of curves, and in (c)--(e) and (g) a periodic component of $\PC(f)$ is not normal. Row (g) is the cone over (c) (Example~\ref{ex:cone-construction}); its component $V=\{x_0^3=x_1^2x_2\}$ is not normal along the line $\{x_0=x_1=0\}$, and the normalization of $V$, the cone over the twisted cubic, is singular at the apex $[0:0:0:1]$, a fixed point at which all multipliers vanish.

The selected points can be specified as follows. In~(d), take $\eta(2,2,2)$ and $\eta(t,t,t)$, where $t^2+t-1=0$. The one-variable return multipliers are $4$ and $-4$, respectively, and the three symmetric-product multipliers are their first three powers. In~(e), take
\[
[1:4:6:4]\quad\text{and}\quad
[1:1+3\omega:-3:1+3\omega^2],\qquad \omega^3=1,\quad\omega\neq1.
\]
The second point is exchanged with its conjugate. The characteristic polynomials of the corresponding first-return differentials are $(T-4)(T-8)(T-16)$ and $(T-4)(T-16)(T-64)$.

In~(f), the line $[1:t:1:t]$ is invariant, and the induced map is $h(t)=(1-2t)^{-2}$. The fixed point $t=1$ has characteristic polynomial $(T+4)(T^2+16)$. The two roots of $4t^2-8t+1=0$ form a $2$-cycle, on which the characteristic polynomial of $Df^2$ is
\[
(T+2)(T^2+8T-32).
\]
The ten planes are $\{x_i=0\}$ and $\{x_i=x_j\}$, for $0\leq i<j\leq3$. Setting $x_4=0$ identifies their intersections with the type $A_4$ braid arrangement. Its codimension-two intersections correspond to partitions of five indices into three blocks: ten with one triple and fifteen with two pairs. This gives the $25$ lines in the table. In~(g), the displayed $2$-cycle is $[0:0:1:1]\leftrightarrow[0:0:1:-1]$; its spectrum follows from Example~\ref{ex:cone-construction}.
\end{example}

\subsection{Singular ambient varieties}\label{ssec:ex-singular-ambient}

\begin{example}\label{ex:kummer}
Let $A$ be an abelian surface, $\pi\colon A\to X:=A/\{\pm1\}$, and $g$ the endomorphism induced by $[k]$, $k\geq2$, polarized with $M=k^2$: $[k]^*L\simeq L^{\otimes k^2}$ for a symmetric ample $L$ \cite[Corollary~2.3.6]{BL04}, and the kernel of $\pi^*$ on $\operatorname{Pic}X$ is killed by $2$. The map $g$ is locally biholomorphic at $\pi(a)$ unless $a\notin A[2]$ and $ka\in A[2]$, so $\Br(g^n)\subset\Sing X$ for all $n$: $\PC(g)=\Sing X$, the sixteen $A_1$ points $\pi(A[2])$, is a branch boundary without divisorial part; an $A_1$ point is a singularity locally isomorphic to $\CC^2/\{\pm1\}$. For $k$ odd these points are fixed, and $g^*$ acts on the cotangent space $\langle u^2,uv,v^2\rangle$ of $\CC^2/\{\pm1\}$ by $k^2=M$.
\end{example}

\begin{example}\label{ex:P112}
Let $X:=\PP(1,1,2)$, with coordinates $[x:y:z]$ of weights $1,1,2$ and the $A_1$ point $v:=[0:0:1]$. Let $d=2k+1\geq3$, $a,b\in\CC$ with $((d-1)a)^{d-1}=((d-1)b)^{d-1}=d^{d-1}$, and $g[x:y:z]:=[x^d+axz^k:y^d+byz^k:z^d]$. Through $(x^2,xy,y^2,z)$, $X$ is the quadric cone $Q\subset\PP^3$, and $g^*\cO_Q(1)\simeq\cO_Q(d)$: $g$ is polarized with $M=d$. The set $B:=\{z=0\}\cup\{dx^{2k}+az^k=0\}\cup\{dy^{2k}+bz^k=0\}\cup\{v\}$ is a branch boundary: it is forward invariant, and its curves form the zero set of the Jacobian determinant of the lift of $g$ to $\CC^3$. Near $v$, $g$ lifts to the odd local biholomorphism $(x,y)\mapsto(x^d+ax,y^d+by)$; so $g$ is a local biholomorphism at the fixed point $v$, and $g^*$ has the eigenvalues $a^2,ab,b^2$ on the cotangent space at $v$, of modulus $(d/(d-1))^2\to1$. Thus Theorem~C admits no lower bound $|\lambda|\geq c>1$ independent of $M$. For $d=3$ and $a=b=3/2$, the spectrum at $v$ is $\{9/4\}$, and $\PC(g)=B$ consists of three invariant curves and the isolated point $v$.
\end{example}

\begin{example}\label{ex:mixed-dimension}
On $X:=(A\times\PP^1)/\langle\sigma_1,\sigma_2\rangle$, with $A$ an abelian surface, $\sigma_1(z,w):=(-z,1/w)$ and $\sigma_2(z,w):=(-z+a_0,w)$, $a_0\in A[2]\setminus\{0\}$, the map $g$ induced by $f:=([k],w^{k^2})$, $k\geq3$ odd, has a post-critical set with components of dimensions $2$, $1$ and $0$. Indeed, $A\times\{0,\infty\}$ is a branch boundary for $f$, which commutes with $\sigma_1$ and $\sigma_2$; the fixed points of $\sigma_2$ form the sixteen curves $\{2z=a_0\}\times\PP^1$, those of $\sigma_1$ are the points $(z,\pm1)$, $z\in A[2]$, which lie neither on these curves nor on $A\times\{0,\infty\}$, and $\sigma_1\sigma_2$ has none. By Lemma~\ref{lem:pcf-transport}(e) and the paragraph after it, the images of these three sets form a branch boundary for $g$; the images of the curves and of the points are $\Sing X$, eight curves and sixteen points, and the image of $A\times\{0,\infty\}$ is $\PCdiv(g)$ by Lemma~\ref{lem:pcf-transport}(a), as the quotient map is quasi-\'etale. So this branch boundary is $\PC(g)$ (Proposition~\ref{prop:divisorial}).
\end{example}

\begin{example}\label{ex:polarized}
Theorem~C fails for non-polarized maps, at a singular fixed point at which the map is a local biholomorphism. Let $g(z,w):=(z^2,\lambda w)$ on $\PP^1\times\PP^1$, $\lambda\in\CC^*$, $\sigma(z,w):=(1/z,-w)$, and $g_X$ the map induced on the normal surface $X:=(\PP^1\times\PP^1)/\sigma$. The image $B$ under $\pi$ of the lines $\{z=0\}$, $\{z=\infty\}$, $\{w=0\}$, $\{w=\infty\}$ is a branch boundary, and $x_0:=\pi(1,0)$ is a fixed $A_1$ point. In the coordinate $u:=(z-1)/(z+1)$, $\sigma(u,w)=(-u,-w)$ and $g(u,w)=(2u/(1+u^2),\lambda w)$, so $g_X^*$ has the eigenvalues $4$, $2\lambda$, $\lambda^2$ on the cotangent space $\langle u^2,uw,w^2\rangle$ at $x_0$; for $0<|\lambda|\leq1$, $\lambda^2$ is neither $0$ nor of modulus greater than one. The map $g_X$ is not polarized: if $g_X^*H\simeq H^{\otimes M}$ with $H$ ample, then $\pi^*H\simeq\cO(\alpha,\beta)$ with $\alpha,\beta>0$, and $g^*\cO(\alpha,\beta)=\cO(2\alpha,\beta)$ cannot be $\cO(M\alpha,M\beta)$.
\end{example}

\subsection{Non-normal varieties and the limits of Theorem~B}\label{ssec:ex-nonnormal}

\begin{example}\label{ex:transfer}
Let $\pi\colon Y\to X$ be finite and surjective, with $\pi\circ g=f\circ\pi$.
\begin{enumerate}[label=\textup{(\alph*)},leftmargin=*]
\item \emph{Cusp.} For $\pi(t)=(t^2,t^3)$ onto $X=\{v^2=u^3\}$, every germ $g(t)=\alpha t+\beta t^2+\cdots$ descends, $D_0\pi=0$, and $D_0f$ has the eigenvalues $\alpha^2$, $\alpha^3$. Here $N(\pi,0)=3$, and $\alpha^3$ is not a monomial of degree between $1$ and $2$ in $\alpha$ if $\alpha\notin\{0,\pm1\}$.
\item \emph{Node.} For $\pi(s)=(s^2-1,s(s^2-1))$ onto $X=\{v^2=u^2(u+1)\}$ and $g(s)=s-\frac14(s^2-1)+\frac34s(s^2-1)$, which fixes $\pm1$ with $g'(1)=2$ and $g'(-1)=3$, $D_0f$ has the eigenvalues $2$ and $3$: Theorem~B(1) holds, but the conclusion of Theorem~B(2) fails at both points.
\item \emph{Triple point.} For the curve $X$ with $\cO(X)=\{h\in\CC[s]:h(0)=h(1)=h(-1)\}$, normalized by $\CC$, and $g(s)=2s^3-3s$, which fixes $0$ and exchanges $\pm1$, $f^*$ has the eigenvalues $-3,-3,3$ on $\fm_p/\fm_p^2=\bigoplus_{y\in\{0,\pm1\}}\fm_y/\fm_y^2$, and $3$ is not a monomial in $g'(0)=-3$: the assumption of Theorem~B(1) that the periodic points of $\pi^{-1}(p)$ are fixed is needed.
\item \emph{Symmetric products.} At the image in $\CC^n/\mathfrak S_n$ of a fixed point $(z,\dots,z)$ of $h^{\times n}$, the eigenvalues are $\alpha,\dots,\alpha^n$, $\alpha=h'(z)$: the degree $n$ occurs although $X$ and $Y$ are smooth.
\end{enumerate}
\end{example}

\section{Questions}\label{sec:questions}

\begin{question}\label{q:cluster-set}
Let $g$ be a polarized endomorphism of a normal projective variety $X$ with $\dim X\geq2$, and let $\varphi\colon D(0,R)\to X$ be the maximal rotation disc of Theorem~A for $E=\PC(g)$. Is the cluster set $\bigcap_{r<R}\overline{\varphi(\{r<|z|<R\})}$ contained in $\PC(g)$?
\end{question}

By Theorem~A(2), almost every radial limit of $\varphi$ lies in $\PC(g)$.

\begin{question}\label{q:strictness}
Does every polarized endomorphism of a normal projective variety which admits a branch boundary arise from a PCF endomorphism of some $\PP^N$ by finitely many operations of the following two kinds, in any order: restriction to an invariant subvariety followed by normalization, and quotient by a finite group of automorphisms commuting with the map? A test case is the map induced by $[2]$ on the Kummer surface of a simple abelian surface (Example~\ref{ex:kummer}).
\end{question}

\begin{question}\label{q:chebyshev-rank4}
For which root systems of the types $B_n$ \textup($n\geq3$\textup), $D_n$ \textup($n\geq4$\textup), $E_6$, $E_7$, $E_8$, $F_4$ and $G_2$, and which integers $m\geq2$, does the map $P_m$, in the fundamental-character coordinates, extend to an endomorphism of projective space?
\end{question}

By Remark~\ref{rem:foldings}, it does not for $G_2$ and $B_3$ with $m=2$.

\appendix
\section{Arcs with centred lifts}\label{app:arcs}

This appendix constructs, at a fixed point of a map with a branch boundary, a family of arcs which lift through every iterate to arcs centred at the point. It is used in the proof of Proposition~\ref{prop:local-aut}. The construction is that of Sun and Zhao \cite[Lemmas~3.2.5, 4.2.1 and~A.0.4]{SZ26b}; the proofs given here are different.

\subsection{The setting}\label{ssec:arcs}

Let $X$ be an irreducible normal projective variety of dimension $m\geq1$, $g\colon X\to X$ a finite surjective morphism, $B\subsetneq X$ a closed algebraic subset containing $\Sing X$ and $\Br(g^n)$ for every $n\geq1$, and $q\in X$ a fixed point of $g$. A \emph{formal arc} of a complete local $\CC$-algebra $R$ with residue field $\CC$ is a local homomorphism of $\CC$-algebras $R\to\CC[[t]]$, and a \emph{formal arc at $q$} is a formal arc of $\hcO_{X,q}$. Arcs are holomorphic maps defined on discs $D(0,\varepsilon)$, centred at $q$ if they send $0$ to $q$.

Choose a finite surjective morphism $\Pi\colon X\to\PP^m$ \cite[Chapter~1, \S5.4, Theorem~1.17]{Sha13}. The points of $\Reg X$ at which $D\Pi$ is not invertible form a proper closed algebraic subset of $\Reg X$, by generic smoothness in characteristic zero \cite[Chapter~2, \S6.2, Lemma~2.4]{Sha13}. Let $\Delta\subset\PP^m$ be a hypersurface containing $\Pi(q)$, $\Pi(B)$ and the image under $\Pi$ of the closure of this subset; it exists since $\Pi$ is finite. Then $B_1:=\Pi^{-1}(\Delta)$ is a proper closed algebraic subset of $X$ containing $B$, and $\Pi$ is locally biholomorphic on $X\setminus B_1$.

Choose affine coordinates $w$ on $\CC^m\subset\PP^m$ centred at $\Pi(q)$, and let $\delta$ be a polynomial equation of $\Delta\cap\CC^m$. As $\Pi$ is a finite surjective morphism, $\Delta\neq\PP^m$ and $X$ is irreducible, the regular function $\delta\circ\Pi$ is a nonzero element of $\cOa_{X,q}$.

Let $\Omega_0\subset\CC^m$ be the neighbourhood of $0$ given by \cite[Lemma~3.8]{Ued98} for $(\CC^m,\Delta\cap\CC^m,0)$, so that Lemma~\ref{lem:ueda-local} holds with $\Omega=\Omega_0$; it is connected, as noted after that lemma. Let $U$ be the connected component of $\Pi^{-1}(\Omega_0)$ containing $q$, and, for $n\geq0$, let $V_n$ be the connected component of $g^{-n}(U)$ containing $q$. Set
\[
\Pi_n:=\Pi\circ g^n\colon V_n\to\Omega_0,\qquad U^\circ:=U\setminus B_1\qquad\text{and}\qquad V_n^\circ:=V_n\setminus\Pi_n^{-1}(\Delta).
\]

\begin{lemma}\label{lem:projection}
For every $n\geq0$, the map $\Pi_n$ is proper and surjective with finite fibres, and $\Pi_n^{-1}(0)=\{q\}$. The maps $\Pi_n\colon V_n^\circ\to\Omega_0\setminus\Delta$ and $g^n\colon V_n^\circ\to U^\circ$ are finite coverings of connected complex manifolds. Moreover $g(V_{n+1})\subset V_n$.
\end{lemma}

\begin{proof}
The map $g^n\colon g^{-n}(U)\to U$ is proper with finite fibres, as the restriction of the finite map $g^n$ over an open set, and so is its restriction to the open and closed subset $V_n$. The same holds for $\Pi|_U\colon U\to\Omega_0$, hence for $\Pi_n$. Its image is closed, and open by the open mapping theorem \cite[Chapter~5, \S4.3]{GrRe84}, so $\Pi_n$ is surjective. Let $x\in V_n$ with $\Pi_n(x)\notin\Delta$. Then $g^n(x)\notin B_1$, and $B_1$ contains $\Sing X$ and $\Br(g^n)$. So $g^n$ is locally biholomorphic at $x$, the point $x$ is regular, and $\Pi$ is locally biholomorphic at $g^n(x)$; thus $\Pi_n$ is locally biholomorphic at $x$. As $V_n$ is connected, open in the normal space $X$, it is irreducible and normal, and $\Pi_n$ is a $(\Delta\cap\Omega_0)$-branched analytic covering. By Lemma~\ref{lem:ueda-local}, $\Pi_n^{-1}(0)$ is a single point, which is $q$.

We have shown that $V_n^\circ\subset\Reg X$. The manifold $\Reg V_n$ is connected \cite[Chapter~7, \S4.2]{GrRe84}, and $V_n^\circ$ is the complement in it of a proper analytic subset, so $V_n^\circ$ is connected; so is $\Omega_0\setminus\Delta$. The map $\Pi_n\colon V_n^\circ=\Pi_n^{-1}(\Omega_0\setminus\Delta)\to\Omega_0\setminus\Delta$ is a proper local homeomorphism, hence a finite covering, as in the proof of Lemma~\ref{lem:covering}(2). Since $U^\circ=V_0^\circ$ and $V_n^\circ=V_n\cap(g^n)^{-1}(U^\circ)$, the same argument applies to $g^n\colon V_n^\circ\to U^\circ$. Finally, $g(V_{n+1})$ is connected, contains $q$ and lies in $g^{-n}(U)$, hence in $V_n$.
\end{proof}

The next lemma plays the role of \cite[Lemma~A.0.4]{SZ26b}, which transfers the lifting property from a formal arc to the holomorphic arcs with the same contact orders along a normal crossing divisor.

\begin{lemma}\label{lem:arc-lift}
Let $\beta_0\in\CC[[t]]^m$ with $\beta_0(0)=0$ and $k:=\operatorname{ord}_t\delta(\beta_0(t))<\infty$. Let $n\geq0$, and let $\hat\alpha$ be a formal arc at $q$ with $\hat\alpha(w_i\circ\Pi_n)=\beta_{0,i}$ for $1\leq i\leq m$. Let $\beta\colon D(0,\varepsilon)\to\Omega_0$ be holomorphic, with $\beta\equiv\beta_0\bmod t^{k+1}$ and $\beta(D^*(0,\varepsilon))\cap\Delta=\emptyset$. Then there is a holomorphic map $\alpha\colon D(0,\varepsilon)\to V_n$ with $\Pi_n\circ\alpha=\beta$ and $\alpha(0)=q$.
\end{lemma}

\begin{proof}
\emph{Step 1: a convergent arc.} Embed an affine open neighbourhood of $q$ in $X$ into $\CC^\varkappa$, with $q=0$, and let $G_1,\dots,G_l$ be polynomials generating its ideal. Then $\hcO_{X,q}=\CC[[z]]/(G_1,\dots,G_l)$, since completion is exact on finitely generated modules. So $\hat\alpha$ is given by a formal solution $\hat z(t)\in\CC[[t]]^\varkappa$, with $\hat z(0)=0$, of the equations $G(z)=0$, which do not involve $t$. Artin's approximation theorem \cite[Theorem~(1.2), p.~277]{Art68} applies to systems of convergent power series over $\CC$ and their formal solutions without constant term. Applied here in the single variable $t$, it gives a convergent solution $z(t)$ with $z(t)\equiv\hat z(t)\bmod t^{k+1}$. It defines a holomorphic arc $\alpha''\colon D(0,\varepsilon'')\to X$ with $\alpha''(0)=q$, and we may assume that $\alpha''(D(0,\varepsilon''))\subset V_n$. The components of $w\circ\Pi_n$ are regular functions near $q$, hence restrictions of convergent power series in $z$, and modulo $t^{k+1}$ their values on $z(t)$ only depend on $z(t)\bmod t^{k+1}$. Hence $\beta'':=\Pi_n\circ\alpha''\equiv\beta_0\bmod t^{k+1}$.

\emph{Step 2: a homotopy avoiding $\Delta$.} Let $\varepsilon_1:=\min(\varepsilon,\varepsilon'')$, and $\beta_\sigma:=(1-\sigma)\beta''+\sigma\beta$ on $D(0,\varepsilon_1)$, for $\sigma\in[0,1]$. Then $\beta_\sigma\equiv\beta_0\bmod t^{k+1}$, so that
\[
\delta(\beta_\sigma(t))
 =t^k\bigl(\mathsf c+t\,\mathsf e(\sigma,t)\bigr),
\]
where $\mathsf c\neq0$ is the coefficient of $t^k$ in $\delta(\beta_0(t))$ and $\mathsf e$ is continuous on $[0,1]\times\overline{D}(0,\varepsilon_1/2)$. Choose $\varepsilon_2\leq\varepsilon_1/2$ with $\varepsilon_2\sup|\mathsf e|<|\mathsf c|$ and $\beta_\sigma(D(0,\varepsilon_2))\subset\Omega_0$ for all $\sigma$. Then $\beta_\sigma(D^*(0,\varepsilon_2))\subset\Omega_0\setminus\Delta$ for all $\sigma$.

\emph{Step 3: monodromy.} Let $\ell(\tau):=\frac{\varepsilon_2}2\exp(2\pi i\tau)$, $\tau\in[0,1]$. By Lemma~\ref{lem:projection}, the loop $\beta''\circ\ell$ has the closed lift $\alpha''\circ\ell$ through the covering $\Pi_n\colon V_n^\circ\to\Omega_0\setminus\Delta$. Lift the homotopy $(\sigma,\tau)\mapsto\beta_\sigma(\ell(\tau))$ to $H\colon[0,1]^2\to V_n^\circ$ with $H(0,\cdot)=\alpha''\circ\ell$. The paths $\sigma\mapsto H(\sigma,0)$ and $\sigma\mapsto H(\sigma,1)$ lift the same path and have the same origin, so they coincide, and $H(1,\cdot)$ is a closed lift of $\beta\circ\ell$.

\emph{Step 4: the lift.} As $\pi_1(D^*(0,\varepsilon))$ is generated by the class of $\ell$, the lifting criterion gives a continuous, hence holomorphic, map $\alpha\colon D^*(0,\varepsilon)\to V_n^\circ$ with $\Pi_n\circ\alpha=\beta$. As $t\to0$, the points $\alpha(t)$ stay in a compact set, since $\Pi_n$ is proper, and accumulate only in $\Pi_n^{-1}(0)=\{q\}$. In a local embedding of $(X,q)$, the components of $\alpha$ are therefore bounded near $0$, and $\alpha$ extends holomorphically to $D(0,\varepsilon)$ with $\alpha(0)=q$, by Riemann's theorem.
\end{proof}

The next proposition is due to Sun and Zhao \cite[Lemma~4.2.1]{SZ26b} (see also \cite[Lemma~3.2.5]{SZ26b}), who state it with generic rank $m$ and with centres which may vary with $n$. Their proof uses a log resolution of $(X,B)$; we give one which uses the finite projection $\Pi$, Lemma~\ref{lem:ueda-local} and Artin's approximation theorem.

\begin{proposition}[Sun--Zhao {\cite[Lemma~4.2.1]{SZ26b}}]\label{prop:family}
Let $\hat\alpha_0$ be a formal arc at $q$ with $\hat\alpha_0(\delta\circ\Pi)\neq0$. Assume that for every $n\geq1$ there is a formal arc $\hat\alpha_n$ at $q$ with $\hat\alpha_n\circ(g^*)^n=\hat\alpha_0$. Then there are $\varepsilon>0$, a ball $\BB\subset\CC^m$ centred at $0$ and a holomorphic map $\Xi\colon D(0,\varepsilon)\times\BB\to X$ with the following properties.
\begin{enumerate}[label=\textup{(\arabic*)}]
\item $\Xi(0,s)=q$ for every $s\in\BB$.
\item $\Xi$ maps $D^*(0,\varepsilon)\times\BB$ into $X\setminus B_1$, and has rank $m$ at every point of $D^*(0,\varepsilon)\times\BB$.
\item For every $s\in\BB$ and $n\geq1$, there is a holomorphic map $\eta_{s,n}\colon D(0,\varepsilon)\to X$ with $\eta_{s,n}(0)=q$ and $g^n\circ\eta_{s,n}=\Xi(\cdot,s)$.
\end{enumerate}
\end{proposition}

\begin{proof}
Let $\beta_0\in\CC[[t]]^m$ be given by $\beta_{0,i}:=\hat\alpha_0(w_i\circ\Pi)$. By assumption, $k:=\operatorname{ord}_t\delta(\beta_0(t))$ is finite. Let $\beta_0^{\leq k}$ be the sum of the terms of degree at most $k$ of $\beta_0$, and $\beta_s(t):=\beta_0^{\leq k}(t)+t^{k+1}s$ for $s\in\CC^m$. As in Step~2 of the proof of Lemma~\ref{lem:arc-lift}, $\delta(\beta_s(t))=t^k(\mathsf c+t\,\mathsf e(s,t))$ with $\mathsf e$ bounded for $s$ bounded. So there are $\varepsilon>0$ and a ball $\BB$ centred at $0$ such that $\beta_s(D(0,\varepsilon))\subset\Omega_0$ and $\beta_s(D^*(0,\varepsilon))\cap\Delta=\emptyset$ for $s\in\BB$. We set $b(t,s):=\beta_s(t)$.

\emph{Step 1: lifts for each $n$.} For $n\geq0$ we have $\hat\alpha_n(w_i\circ\Pi_n)=\hat\alpha_n\circ(g^*)^n(w_i\circ\Pi)=\beta_{0,i}$. By Lemma~\ref{lem:arc-lift}, for every $s\in\BB$ there is a holomorphic map $\eta'\colon D(0,\varepsilon)\to V_n$ with $\Pi_n\circ\eta'=\beta_s$ and $\eta'(0)=q$. The domain $D(0,\varepsilon)$ and the ball $\BB$ were chosen from $\beta_0$ and $\delta$ alone, so they are independent of $n$.

\emph{Step 2: a lift with lifts of all orders.} Let $\mathcal E$ be the set of holomorphic maps $\eta\colon D(0,\varepsilon)\to U$ with $\Pi\circ\eta=\beta_0^{\leq k}$. Such a map sends $D^*(0,\varepsilon)$ into $U^\circ$, and a lift through the finite covering $U^\circ\to\Omega_0\setminus\Delta$ is determined by its value at one point, so $\mathcal E$ is finite.

For $n\geq0$, let $\mathcal E_n\subset\mathcal E$ be the set of maps $g^n\circ\eta'$, where $\eta'\colon D(0,\varepsilon)\to V_n$ is holomorphic with $\Pi_n\circ\eta'=\beta_0^{\leq k}$. By Step~1 with $s=0$, the set $\mathcal E_n$ is nonempty. If $\eta'\colon D(0,\varepsilon)\to V_{n+1}$ is as in the definition of $\mathcal E_{n+1}$, then $g\circ\eta'$ maps $D(0,\varepsilon)$ into $V_n$ by Lemma~\ref{lem:projection}, and $\Pi_n\circ g\circ\eta'=\beta_0^{\leq k}$. So $\mathcal E_{n+1}\subset\mathcal E_n$, and we may fix $\eta_*\in\bigcap_n\mathcal E_n$.

\emph{Step 3: the family.} Let $\ell(\tau):=\frac\varepsilon2\exp(2\pi i\tau)$. The group $\pi_1(D^*(0,\varepsilon)\times\BB)$ is generated by the class of $\tau\mapsto(\ell(\tau),0)$, whose image under $b$ is $\Pi\circ\eta_*\circ\ell$. By the lifting criterion for the covering $\Pi\colon U^\circ\to\Omega_0\setminus\Delta$, there is a holomorphic map $\Xi\colon D^*(0,\varepsilon)\times\BB\to U^\circ$ with $\Pi\circ\Xi=b$ and $\Xi(\cdot,0)=\eta_*$. As $(t,s)\to(0,s_1)$, the accumulation points of $\Xi(t,s)$ lie in $U\cap\Pi^{-1}(0)=\{q\}$, since $\Pi|_U$ is proper and $b(0,s)=0$. So $\Xi$ extends continuously by $\Xi(0,s):=q$, and holomorphically by Riemann's extension theorem, applied in a local embedding of $(X,q)$. This gives (1). Since $\partial b/\partial s=t^{k+1}\,\Id$ is invertible for $t\neq0$, and $\Pi$ is locally biholomorphic on $U^\circ$, we get (2).

\emph{Step 4: lifts for every $s$.} Fix $n\geq1$ and $s\in\BB$, and write $\eta_*=g^n\circ\eta'$ as in the definition of $\mathcal E_n$. The loop $\Xi(\ell(\cdot),0)=g^n\circ\eta'\circ\ell$ has the closed lift $\eta'\circ\ell$ through the covering $g^n\colon V_n^\circ\to U^\circ$. Lifting the homotopy $\sigma\mapsto\Xi(\ell(\cdot),\sigma s)$, $\sigma\in[0,1]$, which takes values in $U^\circ$, we see as in Step~3 of the proof of Lemma~\ref{lem:arc-lift} that the loop $\Xi(\ell(\cdot),s)$ has a closed lift. By the lifting criterion, there is a holomorphic map $\eta_{s,n}\colon D^*(0,\varepsilon)\to V_n^\circ$ with $g^n\circ\eta_{s,n}=\Xi(\cdot,s)$. As $g^n|_{V_n}$ is proper, the accumulation points of $\eta_{s,n}(t)$ as $t\to0$ lie in $V_n\cap g^{-n}(q)\subset\Pi_n^{-1}(0)=\{q\}$, and $\eta_{s,n}$ extends holomorphically with $\eta_{s,n}(0)=q$, as in Step~4 of that proof.
\end{proof}

\section{Normalization boundaries on normal varieties}\label{app:whitney}

We prove Proposition~\ref{prop:normalization-boundary} for a polarized endomorphism $g$ of a normal projective variety $X$, using a conical neighbourhood along a stratum on which the pair $(X,B)$ is topologically trivial.

\subsection{Conical neighbourhoods}

Following Mather \cite[Section~5]{Mat12}, a \emph{Whitney stratification} of a smooth manifold $\mathcal X$ is a locally finite partition into smooth submanifolds, the \emph{strata}, satisfying the frontier condition and Whitney's condition~(b): a stratum which meets the closure of another stratum is contained in that closure, and every pair of strata satisfies condition~(b).

For a nonempty topological space $\mathcal K$, its \emph{open cone} is the quotient
\[
\operatorname{cone}(\mathcal K):=(\mathcal K\times[0,1))/(\mathcal K\times\{0\}).
\]
We denote its vertex by $*$ and the class of $(l,\tau)$ by $[l,\tau]$. We also set $\operatorname{cone}(\emptyset):=\{*\}$. The following conical form of Thom's first isotopy lemma follows from Mather's construction \cite{Mat12}; see also \cite[(4.14)]{Ver76}.

\begin{lemma}\label{lem:cone}
Let $\mathcal X$ be a smooth manifold with a Whitney stratification with finitely many strata, let $S$ be a stratum, of dimension $i$, and let $y\in S$. There are an open neighbourhood $U'$ of $y$ in $\mathcal X$, a compact space $\mathcal K$ which is the disjoint union of subsets $\mathcal K_{S'}$ indexed by the strata $S'\neq S$, and a homeomorphism
\begin{equation}\label{eq:cone}
\Upsilon\colon\RR^i\times \operatorname{cone}(\mathcal K)\to U'
\end{equation}
such that $\Upsilon(\RR^i\times\{*\})=S\cap U'$ and $\Upsilon(\RR^i\times\mathcal K_{S'}\times(0,1))=S'\cap U'$ for every stratum $S'\neq S$.
\end{lemma}

\begin{proof}
Set $k:=\dim\mathcal X-i$. If $k=0$, then $S$ is open in $\mathcal X$, and we take $\mathcal K:=\emptyset$ and a coordinate ball $U'\subset S$ centred at $y$. Assume that $k\geq1$.

\emph{Step 1: a tube on which the strata are transverse to the spheres.} As $S$ is a submanifold, a chart of $\mathcal X$ at $y$ identifies an open neighbourhood $\mathcal X_0$ of $y$ with $b_0\times\{v\in\RR^k:|v|<\delta_0\}$, where $b_0\subset\RR^i$ is an open ball centred at $0$ and $\delta_0>0$, so that $y=(0,0)$ and $S\cap\mathcal X_0=b_0\times\{0\}$. In the coordinates $(x,v)$ we set $\mathcal T(x,v):=(x,|v|^2)$. The traces of the strata form a Whitney stratification of $\mathcal X_0$, and the chart is a tubular neighbourhood of $S\cap\mathcal X_0$ in the sense of \cite[Section~6]{Mat12}, with tubular projection $(x,v)\mapsto x$ and tubular function $|v|^2$. By \cite[Lemma~7.3]{Mat12}, for every stratum $S'\neq S$ there is a positive continuous function $\varepsilon_{S'}$ on $b_0$ such that $\mathcal T$ restricted to $\{(x,v)\in S'\cap\mathcal X_0:|v|<\varepsilon_{S'}(x)\}$ is a submersion into $\RR^i\times\RR$. Let $b\subset b_0$ be the concentric ball of half the radius. As there are finitely many strata, there is $\delta\in(0,\delta_0)$ with $\delta<\varepsilon_{S'}$ on $\overline b$ for every $S'$. Set
\[
U':=\{(x,v):x\in b,\ |v|<\delta\},\quad \mathcal A:=\{(x,v):x\in b,\ 0<|v|<\delta\}\quad\text{and}\quad P:=b\times(0,\delta^2).
\]
Then $U'$ is an open neighbourhood of $y$ with $S\cap U'=b\times\{0\}$ and $U'\setminus S=\mathcal A$, the map $\mathcal T$ sends $\mathcal A$ into $P$, and its restriction to $S'\cap\mathcal A$ is a submersion into $P$ for every stratum $S'\neq S$. Moreover $\mathcal T\colon\mathcal A\to P$ is proper, as $\mathcal T^{-1}(K)\cap\mathcal A$ is closed and bounded in $\RR^i\times\RR^k$ for every compact set $K\subset P$.

\emph{Step 2: a trivialization of $\mathcal A$ over $P$.} The sets $S'\cap\mathcal A$ form a Whitney stratification of the manifold $\mathcal A$, and $\mathcal T$ is a smooth map whose restriction to each stratum is a submersion into $P$; it is also a submersion on $\mathcal A$, as $v\neq0$ there. By \cite[Proposition~7.1]{Mat12}, there is a family of control data for this stratification compatible with $\mathcal T$, which makes $\mathcal A$ an abstract stratified set and $\mathcal T$ a proper controlled submersion, as in the proof of \cite[Proposition~11.1]{Mat12}. Set
\[
\mathcal K:=\mathcal T^{-1}\bigl(0,\delta^2/2\bigr)=\bigl\{(0,v):|v|^2=\delta^2/2\bigr\},
\]
a sphere. Identify $P$ with $\RR^{i+1}$ by a diffeomorphism mapping $(0,\delta^2/2)$ to $0$; then $\mathcal T$ remains a proper controlled submersion, with fibre $\mathcal K$ over $0$. The proof of \cite[Corollary~10.2]{Mat12} gives a global product trivialization over $\RR^{i+1}$. It respects the strata by \cite[Corollary~10.3]{Mat12}. Thus there is a homeomorphism $\mathcal H\colon\mathcal A\to\mathcal K\times P$ with
\[
\mathrm{pr}_P\circ\mathcal H=\mathcal T
\qquad\text{and}\qquad
\mathcal H(S'\cap\mathcal A)=(S'\cap\mathcal K)\times P
\]
for every stratum $S'\neq S$. Set $\mathcal K_{S'}:=S'\cap\mathcal K$; as $\mathcal K$ does not meet $S$, it is the disjoint union of these sets.

\emph{Step 3: the homeomorphism $\Upsilon$.} Fix a diffeomorphism $\vartheta\colon\RR^i\to b$, and set
\[
\Upsilon(a,*):=\bigl(\vartheta(a),0\bigr)\qquad\text{and}\qquad\Upsilon(a,[l,\tau]):=\mathcal H^{-1}\bigl(l,(\vartheta(a),\delta^2\tau)\bigr)\quad\text{for }l\in\mathcal K,\ 0<\tau<1.
\]
The map $\Upsilon$ is a bijection of $\RR^i\times\operatorname{cone}(\mathcal K)$ onto $(b\times\{0\})\sqcup\mathcal A=U'$. It satisfies
\[
\begin{aligned}
\Upsilon(\RR^i\times\{*\})&=S\cap U',\\
\Upsilon(\RR^i\times\mathcal K_{S'}\times(0,1))
 &=\mathcal H^{-1}(\mathcal K_{S'}\times P)=S'\cap U'.
\end{aligned}
\]
Its restriction to the open set $\RR^i\times\mathcal K\times(0,1)$ is a homeomorphism onto $\mathcal A$. Let $a_0\in\RR^i$. As $\mathcal K$ is compact, the sets $N\times(\{*\}\sqcup\mathcal K\times(0,t))$, with $N$ a neighbourhood of $a_0$ and $0<t<1$, form a basis of neighbourhoods of $(a_0,*)$. Since $\mathcal T\circ\Upsilon(a,[l,\tau])=(\vartheta(a),\delta^2\tau)$, the map $\Upsilon$ sends such a set onto $\{(x,v):x\in\vartheta(N),\ |v|<\delta\sqrt t\}$, and these sets form a basis of neighbourhoods of $\Upsilon(a_0,*)$ in $U'$. Hence $\Upsilon$ is a homeomorphism.
\end{proof}

\subsection{Proof of Proposition~\ref{prop:normalization-boundary}}

Let $B$ be a branch boundary for $g$, let $Z\subset X$ be an irreducible invariant subvariety of positive dimension, and let $\nu\colon X_Z\to Z$ be its normalization. We write $g_Z$ for the induced map, and fix a polarization $g^*H\simeq H^{\otimes M}$.

\begin{proof}
The line bundle $\nu^*H$, the pull-back of the ample line bundle $H|_Z$, is ample, since $\nu$ is finite and therefore affine \cite[Remark~13.61(1), Propositions~13.65 and~13.50(1)]{GW10}, \cite[Tag~0B5V]{Stacks}, and $g_Z^*\nu^*H\simeq\nu^*g^*H\simeq(\nu^*H)^{\otimes M}$. So $g_Z$ is polarized. Set $s:=\dim Z$.

\emph{Step 1: a stratum of $Z$ and its neighbourhoods.} Embed $X$ in $\PP^L$. By \cite[D\'efinition~(2.1) and Th\'eor\`eme~(2.2)]{Ver76}, there is a finite stratification of $\PP^L$ into connected smooth Zariski locally closed subvarieties, compatible with $X$, $B$ and $Z$, which satisfies Verdier's condition~(w). Its strata are locally subanalytic, and condition~(w) implies Whitney's condition~(b) \cite[Th\'eor\`eme~(1.5)]{Ver76}, so Lemma~\ref{lem:cone} applies. Since $Z$ is irreducible, one of its strata, denoted by $S$, is Zariski open and dense in $Z$, and $S\subset\Reg Z$. Let $i$ be the real dimension of $S$.

Fix $y\in S$, and let $U'$, $\mathcal K$ and $\Upsilon$ be as in Lemma~\ref{lem:cone}. If $S'\neq S$ is a stratum contained in $Z$ with $\mathcal K_{S'}\neq\emptyset$, then $S'$ meets every neighbourhood of $y$, which is impossible because $S$ is open in $Z$; hence $U'\cap Z=U'\cap S$. Let $\mathcal K_X$ (resp.\ $\mathcal K_B$) be the union of the sets $\mathcal K_{S'}$ over the strata $S'\neq S$ contained in $X$ (resp.\ in $B$), and set
\[
U:=U'\cap X=\Upsilon\bigl(\RR^i\times(\{*\}\sqcup\mathcal K_X\times(0,1))\bigr).
\]
It is an open neighbourhood of $y$ in $X$ with $U\cap Z=U\cap S$. It is connected, since the path $\sigma\mapsto\Upsilon(a,[l,\sigma\tau])$, $\sigma\in[0,1]$, joins $\Upsilon(a,*)$ to $\Upsilon(a,[l,\tau])$ in $U$. As $U$ is a connected normal complex space, $\Reg U$ is a connected manifold \cite[Chapter~7, \S4.2]{GrRe84}, and $U\setminus B$ is the complement in it of a proper analytic subset, hence a connected manifold.

Let $u=\Upsilon(a,*)\in U\cap S$. For $r>0$ and $0<t\leq1$, set $U_u(r,t):=\Upsilon(b(a,r)\times(\{*\}\sqcup\mathcal K_X\times(0,t)))$, where $b(a,r)\subset\RR^i$ is the open ball. As $\mathcal K$ is compact, these sets form a fundamental system of open neighbourhoods of $u$ in $U$. We claim that the inclusion
\begin{equation}\label{eq:small}
U_u(r,t)\setminus B\hookrightarrow U\setminus B
\end{equation}
is a homotopy equivalence. Indeed, as $B$ is a union of strata, $\Upsilon$ identifies $U\setminus B$ with $\RR^i\times C_1$ and $U_u(r,t)\setminus B$ with $b(a,r)\times C_t$, where $C_t$ is the union of $(\mathcal K_X\setminus\mathcal K_B)\times(0,t)$ and, if $S\not\subset B$, of the vertex. Let $\varsigma\colon\RR^i\to b(a,r)$ be a radial homeomorphism centred at $a$. The map $(x,[l,\tau])\mapsto(\varsigma(x),[l,t\tau])$ sends $\RR^i\times C_1$ into $b(a,r)\times C_t$, and the maps
\[
(x,[l,\tau])\longmapsto\bigl((1-\sigma)x+\sigma\varsigma(x),[l,(1-\sigma+\sigma t)\tau]\bigr),\qquad\sigma\in[0,1],
\]
form a homotopy from the identity to this map, which preserves $\RR^i\times C_1$ and $b(a,r)\times C_t$. This proves the claim. In particular, $U_u(r,t)\setminus B$ is path connected, and $\pi_1(U_u(r,t)\setminus B)\to\pi_1(U\setminus B)$ is onto.

\emph{Step 2: fibres over $U\cap S$.} Let $y$ and $U$ be as in Step~1, let $n\geq1$, let $V$ be a connected component of $g^{-n}(U)$, and set $\eta:=g^n|_V\colon V\to U$. We claim that every point of $U\cap S$ has exactly one preimage under $\eta$.

The map $\eta$ is proper with finite fibres, as $V$ is open and closed in $g^{-n}(U)$. Since $V$ is a connected open subset of the normal space $X$, the manifold $\Reg V$ is connected \cite[Chapter~7, \S4.2]{GrRe84}. Set $V^\circ:=V\setminus\eta^{-1}(B)$. If $x\in V^\circ$, then $g^n$ is locally biholomorphic at $x$, because $\Br(g^n)\subset B$, and $\eta(x)\in\Reg X$, so $x\in\Reg X$. Moreover $\eta^{-1}(B)$ is an analytic subset of $V$ of dimension less than $\dim X$, as $\eta$ has finite fibres. Thus $V^\circ$ is the complement of a proper analytic subset of $\Reg V$, hence a connected manifold, and $\eta\colon V^\circ\to U\setminus B$ is a proper local homeomorphism, hence a finite covering, as in the proof of Lemma~\ref{lem:covering}(2). Its image $U\setminus B$ is dense in $U$, and $\eta(V)$ is closed, so $\eta(V)=U$.

Let $u\in U\cap S$, and write $\eta^{-1}(u)=\{x_1,\dots,x_k\}$, with $k\geq1$. Choose pairwise disjoint open sets $N_j\ni x_j$ in $V$. As $\eta$ is proper, there is a neighbourhood $U_u=U_u(r,t)$ of $u$ with $\eta^{-1}(U_u)\subset\bigcup_jN_j$. We show that $Y:=\eta^{-1}(U_u\setminus B)$ is connected. Fix $u'\in U_u\setminus B$. Every point of $Y$ is joined in $Y$ to a point of $\eta^{-1}(u')$, by lifting a path in $U_u\setminus B$. Let $v,v'\in\eta^{-1}(u')$, and join $v$ to $v'$ by a path in the connected manifold $V^\circ$. Its image is a loop in $U\setminus B$ based at $u'$, which is homotopic, with fixed endpoints, to a loop $\gamma$ in $U_u\setminus B$, since the map induced on fundamental groups by \eqref{eq:small} is onto. By the homotopy lifting property, the lift of $\gamma$ starting at $v$ ends at $v'$, and it lies in $Y$. So $Y$ is connected. On the other hand, $Y$ is contained in the disjoint union of the $N_j$ and meets each of them, as $N_j\cap\eta^{-1}(U_u)$ is a nonempty open subset of $V$ and $\eta^{-1}(B)$ is nowhere dense in $V$. Hence $k=1$.

\emph{Step 3: $g_Z^n$ is locally biholomorphic over $\nu^{-1}(S)$.} Let $n\geq1$, and let $\tilde x\in X_Z$ with $g_Z^n(\tilde x)\in\nu^{-1}(S)$. Set $x:=\nu(\tilde x)$ and $u:=g^n(x)=\nu(g_Z^n(\tilde x))\in S$. Let $U$ be as in Step~1 for $y:=u$, let $V$ be the connected component of $g^{-n}(U)$ containing $x$, and $\eta:=g^n|_V$. The set $S_U:=U\cap S=U\cap Z$ is a closed submanifold of $U$ of dimension $s$. As $g^n(Z)\subset Z$, we have $\eta(Z\cap V)\subset S_U$, and $\eta$ is injective on $Z\cap V$ by Step~2.

We claim that the germ of $Z$ at every point $x'\in Z\cap V$ is irreducible. Otherwise it has two distinct irreducible components, with representatives $Z_1$ and $Z_2$ near $x'$, of dimension $s$. The map $\eta|_{Z_2}\colon Z_2\to S_U$ has finite fibres, so it is open by the open mapping theorem \cite[Chapter~5, \S4.3]{GrRe84}, and $\eta(Z_2)$ is a neighbourhood of $\eta(x')$ in $S_U$. A point $x_1\in Z_1\setminus Z_2$ close to $x'$ then satisfies $\eta(x_1)=\eta(x_2)$ for some $x_2\in Z_2$, which contradicts the injectivity of $\eta$ on $Z\cap V$.

The morphism $\nu$ is an isomorphism over the normal locus of $Z$, a dense Zariski open subset containing $\Reg Z$, hence $S$. We deduce that $\nu$ is injective on the open neighbourhood $N:=\nu^{-1}(Z\cap V)$ of $\tilde x$. Indeed, $N$ and $Z\cap V$ have pure dimension $s$, and $Z\cap V$ is locally irreducible by the claim, so $\nu|_N\colon N\to Z\cap V$ is open by the open mapping theorem. If two distinct points of $N$ had the same image $x'$, the images of two disjoint neighbourhoods of these points would both be neighbourhoods of $x'$ in $Z$, and the points of the normal locus close to $x'$ would have two preimages. Moreover $N':=\nu^{-1}(S_U)$ is an open subset of $\Reg X_Z$, mapped biholomorphically onto $S_U$ by $\nu$, and $\nu\circ g_Z^n=g^n\circ\nu$ gives $g_Z^n|_N=(\nu|_{N'})^{-1}\circ\eta\circ\nu|_N$. This map is injective, hence, by the open mapping theorem, a homeomorphism of $N$ onto an open subset $N''$ of the manifold $N'$. Its restriction to $N\cap\Reg X_Z$ is an injective holomorphic map between manifolds of the same dimension, hence biholomorphic onto its image by Osgood's theorem \cite[Chapter~I, Theorem~2.14]{Ran86}. The inverse of the proper map $N\to N''$ is continuous, and holomorphic outside the image of $N\cap\Sing X_Z$, which is a nowhere dense analytic subset of $N''$ by Remmert's proper mapping theorem \cite[Chapter~10, \S6.1]{GrRe84}. By the Riemann extension theorem in local embeddings of $X_Z$, the inverse is holomorphic. Thus $g_Z^n$ is locally biholomorphic at $\tilde x$, and $\Br(g_Z^n)\cap\nu^{-1}(S)=\emptyset$ for every $n\geq1$.

\emph{Step 4: the branch boundary.} The set $\nu^{-1}(S)$ is a nonempty Zariski open subset of $X_Z$. By Step~3, its complement is a proper closed algebraic subset of $X_Z$ which contains $\Br(g_Z^n)$ for every $n\geq1$. By Lemma~\ref{lem:boundary}(3), $g_Z$ admits a branch boundary.
\end{proof}

\providecommand{\bysame}{\leavevmode\hbox to3em{\hrulefill}\thinspace}
\providecommand{\MR}{\relax\ifhmode\unskip\space\fi MR }
\providecommand{\MRhref}[2]{%
  \href{http://www.ams.org/mathscinet-getitem?mr=#1}{#2}
}
\providecommand{\href}[2]{#2}

\end{document}